# Disproof of the Continuum Hypothesis and Determination of the Cardinality of Continuum by Approximations of Sets

Slavko Rede

**Abstract**. A set theory is developed based on the *approximations* of sets and denoted by AS. In AS the set of all sets exists but the argument for Russell's and Cantor's paradox fail. The Axioms of Separation, Replacement and Foundation are not valid. All the other axioms of ZF are valid and all the basic sets, such as complement, intersection and cartesian product, exist although complement is not quite the same set as in ZF.

The set of all sets can be equipped with the *topology of approximations* ($\mathcal{T}_a$) which is analogical to the interval topology of real numbers in ZF. Every set is closed and every function is continuous in $\mathcal{T}_a$. This implies that the Continuum Hypothesis is false. The sets containing a subset which is perfect in $\mathcal{T}_a$ are of the greatest cardinality. A simple observation shows that the concept of well-ordering must be defined in AS in a slightly different way than in ZF. We prove that a set can be well-ordered if and only if it does not contain a perfect subset. Therefore the cardinalities of arbitrary sets are always comparable without assuming the Axiom of Choice. The cardinals following the smallest infinite cardinal $\omega$ are $2\times\omega$, $3\times\omega$, ..., $\omega^2$, ... $2\times\omega^2$, $3\times\omega^2$, ..., $\omega^3$.... each being of greater cardinality than the previous one, which is not the case in ZF. Immediately after these cardinals does not follow $\omega^\omega$ which is not a well-orderable set but some well-ordered cardinal $\kappa$, and this one is followed by the cardinals $2\times\kappa$, $3\times\kappa$, ..., $\kappa^2$, ... $2\times\kappa^2$, $3\times\kappa^2$, ..., $\kappa^3$...., etc. The greatest cardinal is $\mathcal{P}(\omega)$ and is not a well-orderable set. The cofinality of a well-ordered set is either 2 or $\omega$. The only regular cardinals are 0, 1, 2, $\omega$ and $\mathcal{P}(\omega)$. All other cardinals are singular. The only strong limit cardinal is $\omega$. The only inaccessible cardinal is $\mathcal{P}(\omega)$. Strongly inaccessible cardinals do not exist.

# Introduction

In the 1. section we define a recursive function $F(n, A)$ for an arbitrary natural number $n$ and an arbitrary set $A$. We put $F(0, A) = \emptyset$ and $F(n + 1, A) = \{F(n, x) \mid x \in A\}$. The value of $F(n, A)$ is always a hereditarily finite set (a HF set) of the rank not greater than $n$. If we write { } instead of Ø, we can interpret $F(n, A)$ as a finite word made only of pairs of curly braces { } which are syntactically correctly inserted into each other. When $n$ approaches to infinity, the value of $F(n, A)$ approaches to $A$. It even becomes equal to $A$, if $n$ becomes greater than the rank of $A$. And we prove that $F(\infty, A) = A$ for any set $A$. Therefore, it seems to be natural that the infinite sequence $F(0, A)$, $F(1, A)$, $F(2, A)$, ... should uniquely determine any set $A$. It also seems to be natural that $x \in A$ should hold true for arbitrary sets $x$ and $A$ if $F(n, x) \in F(n + 1, A)$ holds true for every $n$. However, none of these two assumptions is true in the classical set theory, that is in ZF.

In the 2. section we begin the construction of a set theory which is based on both assumptions. We do not adopt any axiom of ZF. We adopt only the logical frame of predicate calculus, and even this one in a limited range. We allow only the existence of arbitrary finite words of some chosen finite alphabet, and the existence of only those predicates which are computable on these words by some algorithm, that is, by a Turing machine. As the only objects in our universe are words, we sketch how a Turing machine can be defined only by the concept of word and without the concept of set, since sets will have in our set theory not the same intuitive meaning as in ZF. Then we define HF sets as special strings made only of curly braces. We also define natural numbers as a special type of HF sets: $0 = \{\}$ and $n + 1 = \{n\}$. Then we define basic operations and relations on natural numbers: sum, difference, <, >, ≤ and ≥. We also prove the validity of inference by induction. We explicitly define when some HF set is an element of another HF set and when two HF sets are equal (Axiom of Extensionality). We calculate the rank of a HF set and by its help we define a computable function $F(n, A)$. This function is a restriction of the function $F(n, A)$ introduced already in the 1. section. Its domain are here not arbitrary sets but only HF sets, and its values are HF sets called *approximations* of $A$. We prove some basic properties of $F$ as, for example, that $F(n, A) = F(n, F(n + m, A))$ holds for every $n$, $m$ and every HF set $A$, that every HF set $A$ is uniquely determined by the infinite sequence $F(0, A)$, $F(1, A)$, $F(2, A)$, ..., and that a HF set $x$ is an element of a HF set $A$ precisely when $F(n, x) \in F(n + 1, A)$ is true for every $n$. Therefore, both assumptions from the 1. section are true for any HF set $A$.

In the 3. section we widen our universe from HF sets to arbitrary sets. These are, in general, not words but are implicitly defined by an axiom, called the Axiom of Set Existence. By this axiom any set $A$ is uniquely determined by the sequence of its approximations $F(0, A)$, $F(1, A)$, $F(2, A)$, ..., and every such sequence determines some set $A$. Since we can choose the terms of this sequence randomly, with the only restriction that $F(n, F(n + 1, A) = F(n, A)$ holds true for every $n$, we can construct sequences which do not correspond to any set from ZF. In our set theory, which we denote by AS (*Approximations of Sets*), there also exist sets which are elements of themselves. Such a set is, for example, the set determined by $F(n + 1, A) = \{F(n, A)\}$ for every $n$. Hence the Axiom of Foundation is not in force in AS. We explicitly define the relation $x \in A$ for arbitrary sets $x$ and $A$ as: $x \in A$ is true precisely when $F(n, x) \in F(n + 1, A)$ is true for every $n$. This relation, as we have already told, holds for HF sets anyway. (As it is known, the relation $\in$ is not defined explicitly in ZF.) We also show that all properties of the function $F(n, A)$, where $A$ is a HF set, are in force if $A$ is an arbitrary set.

All possible sequences $F(0, A)$, $F(1, A)$, $F(2, A)$, ... form an infinite tree. This tree is in fact a computable predicate on words, called the universal tree **U**. The nodes on the $n$th level of **U** are the approximations $F(n, A)$, where $A$ is an arbitrary set. Thus the root of **U** is $F(0, A) = \{\}$, and every node $F(n, A)$ has finitely many successors. The form of these successors is $F(n + 1, A)$ and they fulfil the condition $F(n, F(n + 1, A)) = F(n, A)$. Every infinite path in **U** that starts in the root of **U** represents some set. If such paths are different they represent different sets. Also every infinite subtree of **U** that has the same root as **U** represents some set. The paths in such subtree represent the elements of the corresponding set. This implies that **U** represents some set $U$ as well. Since to every set corresponds a path in **U**, the set $U$ is the universal set.

In 4. section we show that $U$ is the set of all sets. The existence of the set of all sets and the validity of the Axiom of Separation would imply the existence of the set of all sets which are not elements of themselves and this would imply Russell's paradox. However, we show that such a set does not exist. Hence the Axiom of Separation is not in force in AS. Neither is the Axiom of Replacement in force in AS, since the latter implies the former. There arises the question which sets from ZF do exist in AS? In particular, which axioms of ZF are valid in AS? In order to answer this question we first determine some sufficient conditions for the existence of a predicate on sets. Namely, in the 2. section we allowed only computable predicates on HF sets and the application of predicate calculus. And in the 3. section we extended the function $F$ from HF sets to arbitrary sets. Therefore, a predicate on sets certainly exists if it is expressible in predicate calculus by computable predicates on HF sets and by the function $F$ on arbitrary existent sets. Then we prove that for every predicate **P** which has one free variable with domain $U$ there exists the smallest set $P$ whose elements are all $x \in U$ for which **P**$(x)$ is true. Of course, this does not mean that the elements of $P$ are only those $x \in U$ for which **P**$(x)$ is true. For example, the set $\{x \mid x \notin A\}$ does not exist for arbitrary set $A$. That is, the complement of a set, as defined in ZF, does not always exist. But this is an exception, since we can prove the existence of: the union and the intersection of an arbitrary family of sets, unordered and ordered pair of arbitrary two sets, cartesian product of two arbitrary sets, relations and functions on arbitrary sets, power set of an arbitrary set, and the set of all natural numbers. This set, which we denote by $\boldsymbol{N}$, is somewhat different from the one which exists in ZF. It contains an additional element, called the infinite number and denoted by $\infty$, for which holds $F(n + 1, \infty) = \{F(n, \infty)\}$ for every $n$. This implies that $\infty$ is not a regular set. In the same way as the set $\{x \mid x \notin A\}$ does not always exist, the set $\{x \mid x \notin f(x)\}$ does not exist for an arbitrary function $f : A \to \mathcal{P}(A)$ where $\mathcal{P}(A)$ is the power set of some set $A$. For this reason

Cantor's diagonal procedure cannot always be carried out, and the argument for Cantor's paradox fails in AS since $\mathcal{P}(U) = U$.

In the 5. section we introduce a topology into $U$ on the base of which we disprove the Continuum Hypothesis. Since $U$ can be represented by **U**, which is an infinite tree of finite branching degree, $U$ can be interpreted as a set of points in the interval [0, 1]. A neighborhood of some $x \in U$ is the set of all $y \in U$ which have the same $n$th approximation as $x$ for some fixed $n$. Therefore, a neighborhood of $x$ is the set $\{y \mid F(n, x) = F(n, y)\}$. In this topology, called the *topology of approximations* and denoted by $\mathcal{T}_a$, every set is closed and every function is continuous. Therefore, every bijection is a homemorphism. Hence the bijections do not preserve only the sizes of sets but also their topological structure in $\mathcal{T}_a$. We prove that $\mathcal{P}(\boldsymbol{N})$ includes an infinite set whose every element is an accumulation point in it. That is, $\mathcal{P}(\boldsymbol{N})$ includes a set which is perfect in $\mathcal{T}_a$. We also prove that $\boldsymbol{N}$ has only one accumulation point in $\mathcal{T}_a$, namely $\infty$. And we prove that for every natural number $i$ there exists a set $A_i$ which consists of pairs $\{j, k\}$, where $j, k \in \boldsymbol{N}$ and $j < i \leq k$. The set $A_i$ has precisely $i$ accumulation points in $\mathcal{T}_a$ for every $i$. Therefore, there exists a bijection between $\boldsymbol{N}$ and $A_1$ and an injection from $A_i$ to $A_{i+1}$ and to $\mathcal{P}(\boldsymbol{N})$ for every $i$. But there exists no bijection between $A_i$ and $A_{i+1}$ for any $i$, since $A_i$ and $A_{i+1}$ have different number of accumulation points in $\mathcal{T}_a$. All this implies that the Continuum Hypothesis is false in AS.

In sections 6 to 11 we develop cardinal arithmetic and describe the universe of all cardinals giving the position of the cardinality of continuum in this universe. It is known that in ZF the size of sets is measured by the existence of injections between these sets. In AS an injection is a homeomorphic embedding of one set into another in the topology of approximations $\mathcal{T}_a$. This topology is equivalent to the interval topology in ZF. Thus injections are not a proper set size measure in AS. Nevertheless we define that a set $A$ *is not of greater cardinality* then a set $B$ if there exists an injection from $A$ into $B$, denoting this by $A \leq_c B$. If $A \leq_c B$ and $B \leq_c A$ then we say that $A$ and $B$ are of *equal cardinality* and write it $A =_c B$. It is interesting that Cantor-Bernstein Theorem does not hold in AS. It turns out that sets containing a subset which is perfect in $\mathcal{T}_a$ are of the greatest cardinality. Comparison of the cardinalities of two arbitrary sets is in ZFC always possible because any set can be well-ordered. Investigation shows that a sensible definition of well-ordering in AS must be slightly different from the one in ZF. With regard to this definition we prove that a set can be well-ordered in AS if and only if it does not contain a subset which is perfect in $\mathcal{T}_a$. These facts are sufficient to conclude that the cardinalities of two sets in AS are always comparable. Since the interval topology in ZF is equivalent to $\mathcal{T}_a$ in AS, and since every set in AS is closed in $\mathcal{T}_a$, the existence of an injection from one set into another is in AS determined by Cantor-Bendixson rank and Cantor-Bendixson degree of both sets. It turns out that Cantor-Bendixson rank ($r_{CB}$) cannot be a von Neumann ordinal since $\omega$ is the greatest such ordinal in AS. We can overcome this difficulty by defining that any well-ordered set is an ordinal.

The cardinality of the sums of infinite sets in AS has not quite the same properties as in ZFC. If $A \leq_c B$ then $A+B =_c B$ holds only if $r_{CB}(A) <_o r_{CB}(B)$ where the subscript $_o$ denotes comparison of ordinals. However, if $r_{CB}(A) =_o r_{CB}(B)$ then $A+B >_c B$. The cardinality of the product of sets also differs from the one in ZFC. By [H&J] the problem of determining the cardinality of the exponential set $A^B$ where $B$ is infinite has not yet been solved in ZFC. In AS, however, the cardinality of $A^B$ is the same as the cardinality of a perfect set if only $A >_c 1$. We have $A+B <_c A \times B$ for any two infinite $A$ and $B$. And we have $n \times A <_c (n+1) \times A$ for any well-ordered set $A$ and any finite ordinal $n$, which is not the case in ZFC if $A$ is infinite. From this we derive that if $\kappa$ is a cardinal then $n \times \kappa$ is a cardinal for any finite ordinal $n$. We also prove that there are no cardinals between $n \times \kappa$ and $(n+1) \times \kappa$. This implies that the cardinals following immediately after $\omega$ are $2 \times \omega$, $3 \times \omega$, $4 \times \omega$,.... Further we prove that if $\kappa$ is a cardinal whose Cantor-Bendixson degree is 1, then $\omega \times \kappa$ is a cardinal following immediately after the cardinals $\kappa$, $2 \times \kappa$, $3 \times \kappa$, .... Therefore the cardinals $\omega$, $2 \times \omega$, $3 \times \omega$, ... are followed by the cardinal $\omega \times \omega = \omega^2$. Then follow the cardinals $2 \times \omega^2$, $3 \times \omega^2$, $4 \times \omega^2$,.... and then follows the cardinal $\omega \times \omega^2 = \omega^3$. We can see how this goes on. It is interesting that $\omega^\omega$, which is a non-well-orderable set, does not follow immediately after the cardinals $\omega$, $\omega^2$, $\omega^3$, $\omega^4$, .... Immediately after all these cardinals follows some well-ordered cardinal $\kappa_1 <_c \omega^\omega$. And $\kappa_1$ is followed by the cardinals $2\times\kappa_1$, $3\times\kappa_1$, ..., and so on without an end.

A definition of cofinality which gives interesting results in AS is, for example, the one from [G], and can be defined only for well-ordered sets. It turns out that the cofinality of any well-ordered cardinal is not greater than $\omega$. And the only other possible cofinality is 2. The only regular cardinals are 0, 1, 2, $\omega$ and $\mathcal{P}(\omega)$. All other cardinals are singular. The only strong limit cardinal is $\omega$. The only inaccessible cardinal is $\mathcal{P}(\omega)$. Strongly inaccessible cardinals do not exist.

In section 7 we make some intuitive observations about measuring the sizes of infinite sets in AS

# 1. The Membership Structure of Sets

Since in ZF (Zermelo-Fraenkel set theory) all objects are sets, we can define the following function

**1.1 Definition** $F(n, A)$ is a function for which holds

$$F(0, A) = \emptyset = \{\}$$

$$F(n + 1, A) = \{F(n, x) \mid x \in A\}$$

for every set $A$ and every natural number $n$, that is for every $n \in \boldsymbol{N}$.

**1.2 Examples**

1) $F(1, \{\{\}, \{\{\}\}\}) = \{F(0, \{\}), F(0, \{\{\}\})\} = \{\{\}\}$

2) For every $n$ we have $F(n + 1, \emptyset) = \{F(n, x) \mid x \in \emptyset\} = \emptyset = \{\}$

Let's denote the rank of a set $A$ by $r(A)$ assuming that $r(\emptyset) = 0$. The rank of a set is a natural number precisely when this set is hereditarily finite.

**1.3 Theorem** *If $A$ is a hereditarily finite set, then $F(r(A) + m, A) = A$ for every $m \in \boldsymbol{N}$.*

*Proof.* Since {} is the only set of the rank 0, the theorem holds for $r(A) = 0$ by Definition 1.1 and Example 1.2.2. Let the theorem hold for every set of the rank not greater than $n$ and let $r(A) = n + 1$. Then $F(r(A) + m, A) = F(n + 1 + m, A) = \{F(n + m, x) \mid x \in A\}$ by Definition 1.1. Since $r(x) \leq n$ for every $x \in A$, we have by assumption $F(n + m, x) = x$ for every $x \in A$ and every $m \in \boldsymbol{N}$, and this implies $F(r(A) + m, A) = \{x \mid x \in A\} = A$ for every $m \in \boldsymbol{N}$. Thus, by induction on $n$, the theorem holds for any set whose rank is a natural number. That is, it holds for every hereditarily finite set. □

**1.4 Examples**

1) Let natural numbers be written in Zermelo form. That is, let $0 = \{\}$ and let $n + 1 = \{n\}$ for every $n$. Then $F(n, m) = \min\{m, n\}$. Namely, $r(m) = m$ and by Theorem 1.3 we have $F(n, m) = m$ for $n \geq m$. Now, let $n < m$. By Definition 1.1 we have $F(0, m) = \{\} = 0$ for every $m$. If $F(n, m) = n$ for some $n$ and every $m > n$, then we have for every $m > n$

$$F(n + 1, m + 1) = F(n + 1, \{m\}) = \{F(n, m)\} = \{n\} = n + 1$$

Thus $F(n + 1, m) = n + 1$ holds for any $m \geq n + 1$, which by induction on $n$ means that it holds for every $n$.

2) Let $\boldsymbol{N}$ be the set of all natural numbers written in Zermelo form. Then by the previous example we have for $n > 0$

$$F(n, \boldsymbol{N}) = F(n, \{i \mid i \in \boldsymbol{N}\}) = \{F(n-1, i) \mid i \in \boldsymbol{N}\} = \{\min\{n-1, i\} \mid i \in \boldsymbol{N}\} = \{i \mid i \leq n-1\}$$

3) Let $\mathcal{P}(\boldsymbol{N})$ be the power set of $\boldsymbol{N}$. Then $F(0, \mathcal{P}(\boldsymbol{N})) = \{\}$, $F(1, \mathcal{P}(\boldsymbol{N})) = \{\{\}\}$, and for $n > 1$ we have

$$F(n, \mathcal{P}(\boldsymbol{N})) = F(n, \{\{i \mid i \in A\} \mid A \subseteq \boldsymbol{N}\}) = \{\{F(n-2, i) \mid i \in A\} \mid A \subseteq \boldsymbol{N}\} =$$

$$\{\{\min\{n-2, i\} \mid i \in A\} \mid A \subseteq \boldsymbol{N}\} = \{\{i \mid i \leq n-2\} \mid i \in A\} \mid A \subseteq \boldsymbol{N}\}\} = \mathcal{P}(\{i \mid i \leq n-2\})$$

$F(n, A)$ is a hereditarily finite set for every natural number $n$ and every set $A$. The greater the number $n$ is the better approximation of $A$ should $F(n, A)$ be, which is evident from the following observation. We can think of a hereditarily finite set $A$ as of a string of brackets { and }. The brackets at each end of $A$ are the *outside* brackets of $A$. To each of the brackets $A$ is made of we can assign its *depth*. The outside brackets of $A$ have depth 0 in $A$ and if $A \ni x_1 \ni x_2 \ni x_3 \ldots \ni x_k$ is a *membership chain* of the length $k$ in $A$, then the outside brackets of $x_k$ have depth $k$ in $A$. The function $F(n, A)$ operates on $A$ in such a way that it erases from $A$ all brackets of the depth greater than $n$. Therefore, the greater that $n$ is the better approximation of $A$ is $F(n, A)$, since the less brackets are erased from $A$. If $n > r(A)$, then no membership chain in $A$ is of the length greater than $n$ and $F(n, A)$ erases from $A$ no brackets at all. This is also evident from Theorem 1.3. Even if $A$ is not a hereditarily finite set no membership chain in $A$ is of the length greater than $\infty$. In fact every such chain is finite if the Axiom of Foundation holds. Therefore $F(\infty, A) = A$ holds for every set $A$. We can prove this formally as well if we extend Definition 1.1 for $n = \infty$ and if we take into account that $\infty + 1 = \infty$.

**1.5 Theorem** *For every set A holds* $F(\infty, A) = A$.

*Proof.* We shall prove the theorem by transfinite induction on the rank of $A$. The theorem holds for all sets of the rank 0, since the only set of the rank 0 is { }, and for this set we have

$$F(\infty, \{\}) = F(\infty + 1, \{\}) = \{F(\infty, x) \mid x \in \{\}\} = \{\}.$$

Let $F(\infty, x) = x$ hold for every $x$ whose rank smaller than some ordinal $\alpha$ and let the rank $A$ of be $\alpha$. Then we have $F(\infty, A) = F(\infty + 1, A) = \{F(\infty, x) \mid x \in A\} = \{x \mid x \in A\} = A$ which implies that the theorem holds for every set of every rank. □

According to the previous observation we could say that the sequence $F(0, A)$, $F(1, A)$, $F(2, A)$, ... approaches to $A$. If $A_1$ and $A_2$ are two sets and $F(n, A_1) = F(n, A_2)$ holds for every $n \in \boldsymbol{N}$, then we would expect that $F(\infty, A_1) = F(\infty, A_2)$ holds as well. That is, we would expect that $A_1$ and $A_2$ are equal sets. In other words, the sequence $\{F(n, A)\}_{n\in \boldsymbol{N}}$ should uniquely determine every set $A$. However, this is not true in ZF. Let $S(n) = \{i \mid i \leq n-1\}$ where $i$ is a natural number in Zermelo form. Let $A = \{S(i) \mid i \in \boldsymbol{N}\}$. Then it is not difficult to show that the same sequence $\{F(n, A)\}_{n\in \boldsymbol{N}}$ corresponds to different sets $A$ and $A \cup \{\boldsymbol{N}\}$. But if we ignore ZF axioms for a moment, are the sets $A$ and $A \cup \{\boldsymbol{N}\}$ really different? Does it really hold $\boldsymbol{N} \notin A$? Namely, for every $n$ we have $F(n, \boldsymbol{N}) = S(n)$ (See Example 1.4.2) and $F(n+1, A) = \{S(i) \mid i \leq n\}$. Therefore $F(n, \boldsymbol{N}) \in F(n+1, A)$ for every $n$. Thus we should have $F(\infty, \boldsymbol{N}) \in F(\infty, A)$, that is $\boldsymbol{N} \in A$. Let's assume that what should hold by intuition actually holds:

***Principle of Set Equality*** *Let F be the function from Definition* 1.1 *and let A and B be two arbitrary sets. If* $F(n, A) = F(n, B)$ *holds for every* $n \in \boldsymbol{N}$, *then* $A = B$.

The previous example implies that the Principle of Set Equality is not consistent with ZF. We can express our observations also by the following principle:

***Principle of Membership*** *Let F be the function from Definition* 1.1 *and let A and x be two arbitrary sets. If* $F(n, x) \in F(n + 1, A)$ *holds for every* $n \in \boldsymbol{N}$*, then* $x \in A$.

In the sequel we shall construct a theory of sets which will be founded on both principles. Of course this theory will not be compatible with ZF.

# 2. Hereditarily Finite Sets and their Approximations

The construction of a new set theory must begin at the very basis of mathematics. Our set theory will be founded on the concept of string or word. In this section we shall define hereditarily finite sets as a special type of words, and natural numbers as a special type of hereditarily finite sets. We shall assume the predicate calculus with equality but without the concept of set which is used in this calculus in the intuitive meaning as a domain. Neither shall we allow arbitrary predicates. For the time being we shall assume only the existence of those predicates which are decidable for every word in finitely many steps by some algorithm, that is, by a Turing machine.

*Character* is an elementary concept which we shall not define by more elementary concepts as we do not define points in geometry. Characters are letters (capital and small, written in ordinary, italic, bold, subscript or superscript style), numbers, logical, arithmetical and set theory signs, comma, apostrophe, bar and braces. Every character is a *word*. If we write a word $y$ next to a word $x$ we obtain a word denoted by $xy$. By repeating this procedure we can obtain any word. Two words $x$ and $y$ are *equal* ($x = y$) if they are two equal characters, or if $x = x'x''$ and $y = y'y''$ where $x' = y'$ and $x'' = y''$. If words $x$ and $y$ are not equal, then they are *different* ($x \neq y$).

**2.1 Definition**

1) The word {} is a *hereditarily finite set* (HF *set*).
2) If $x$ is a HF set, then $x$ is a *list of* HF *sets* and $x$ is a *term* of this list.
3) If $L$ is a list of HF sets and $x$ is a HF set, then the word $Lx$ is a list of HF sets and $x$ is a term of $Lx$.
4) If $L$ is a list of HF sets, then the word $\{L\}$ is a HF *set*.

This definition implies that the words {}, {{}}, {{}{{}}}, and {{{}}{}} are HF sets while the words {}{{}} and {{}}{} are lists of HF sets. For the sake of greater clearness we shall usually (but not always) place commas between consecutive terms of a list of HF sets. Thus we shall write {{}, {{}, {{}}}} instead of {{}{{}{{}}}}.

HF *sets* are intuitively finite sets whose elements are again HF sets. We shall prove that the Principle of Set Equality and the Principle of Membership, introduced in the previous section, hold for HF sets without any additional axioms. The main reason why we have introduced HF sets is because we shall interpret them as approximations of arbitrary sets. Now we shall define natural numbers as a special type of HF sets. Natural numbers will be represented in Zermelo form (See Example 1.4.1).

**2.2 Definition** A HF set $n$ is a *natural number* if $n = \{\}$ (natural number 0) or if $n = \{m\}$ where $m$ is a natural number. In this case $n$ is the *successor* of $m$ denoted also by $m'$.

For the names of natural numbers we shall almost always use small italic letters $n$, $m$, $i$, $j$, and $k$. We shall denote the constants {}, {{}}, {{{}}}, {{{{}}}}, ... by the usual decimal names 0, 1, 2, 3, .... Instead of „for every natural number $n$“ we shall write „for every $n$“ and instead of „there exists a natural number $n$“ we shall write „there exists $n$“.

**2.3 Definition** The *sum* of natural numbers $n$ and $m$ is a natural number denoted by $n + m$. For every $n$ and every $m$ holds $n + 0 = n$ and $n + m' = (n + m)'$. If $i = j + k$, then $k$ is the *difference* of the numbers $i$ and $j$ ($k = i - j$). In this case $j$ is *smaller or equal* than $i$ ($j \le i$), and $i$ is *greater or equal* than $j$ ($i \ge j$). If $k \ne 0$, then $j$ is *smaller* than $i$ ($j < i$), and $i$ is *greater* than $j$ ($i > j$).

The number $n + 1$ is the successor of the number $n$, since $n + 1 = n + 0' = (n + 0)' = n'$. In the sequel we shall often use inference by induction. For this reason we must prove the following.

**2.4 Theorem** *If a predicate* **P** *holds for* 0 *and if* **P** *holds for n' as soon as it holds for n, then* **P** *holds for every n.*

*Proof.* Every natural number is a word composed only of pairs of curly braces {} inserted into each other. Every such word can be obtained from 0 by the repetitions of the construction of successor. □

It is also obvious that there is no natural number whose successor is 0 and that $n' = m'$ implies $n = m$. Namely, if $\{n\} = \{m\}$, then we put $\{n\} = x\}$ where $x = \{n$ and we put $\{m\} = y\}$ where $y = \{m$. Therefore $x\} = y\}$ and hence $x = y$, that is $\{n = \{m$. From this we obtain $n = m$. Thus all the axioms of the first order theory called number theory are in force so that all theorems of number theory concerning the sum hold. Hence the sum is commutative and associative. Since the relations ≥, ≤, >, < and the function difference are defined in the same way as in number theory, the number $n$ is the successor of the number $n - 1$ for $n > 0$, and $i < j$ and $j < k$ implies $i < k$.

If we have two natural numbers, we can decide by some algorithm in a finite number of steps whether the relations <, >, ≤ and ≥ hold between them or not. In the same way we can decide whether some number is the sum or the difference of two given numbers. We shall say that these relations and functions are *computable*. Of course, using the concept of algorithm demands its definition by the concepts we have introduced so far. An algorithm is formally a Turing machine. And this is defined by a list of characters (an alphabet), list of states, list of accepting states, starting state, tape divided into cells, control unit, reading-writing head, ordered quintuples of characters, states and left-right moves, list of the rules made of these quintuples, and input and output words. All these concepts can be defined by words of different types. For example, the $n$-tuple of words $a_1, a_2, ..., a_n$ is the word $(a_1, a_2, ..., a_n)$ where commas and brackets are part of the word. In order to avoid the term set we have replaced it with the term list which can be defined in a similar way as it was in Definition 2.1.

The relations „is the successor of a natural number ", <, >, ≤, and ≥ are 2-*place* computable relations between natural numbers, while the sum and the difference of two natural numbers are 3-place computable relations or 2-place computable functions between natural numbers. The 2-place relation = is computable even for two arbitrary words. The relations "is a HF set" and "is a natural number" are 1-*place* computable relations for arbitrary word, or computable 1-place *predicates* on words.

Let's emphasize that a predicate on words does not determine the *set* of all those words for which this predicate is true, as it is assumed in the predicate calculus. Namely, in our set theory a group of arbitrary words will *not* be a set, and what is even more important, for every predicate **P** there will in general not exist such a set $P$ that the statement **P**($a$) will be true precisely when $a$ will be an element of $P$. This is not at all surprising if we take into account that the Principle of Set Equality and the Principle of Membership, which we have introduced in the previous section, are not compatible with ZF.

**2.5 Definition** A HF set $x$ is an *element* of a HF set $A$ ($x \in A$) precisely when $A = \{L\}$, where $L$ is a list of HF sets and $x$ is a term of $L$. The denotation $x \notin A$ is an abbreviation for $\neg (x \in A)$.

HF sets $a$, $b$, $c$ and $d$ are the only elements of the HF set $\{a, b, c, d\}$. The HF set {} is the *empty set*, since it has no elements. The 2-place predicate "is an element of a HF set" is computable.

**2.6 Theorem**. *Every nonempty* HF *set is of the form* $\{a_1, a_2, ..., a_k\}$ *where* $a_1, a_2, ..., a_k$ *are precisely all the elements of this* HF *set. Every list of* HF *sets is of the form* $a_1, a_2, ..., a_k$ *where* $a_1, a_2, ..., a_k$ *are precisely all the terms of this list. Here k can be equal to* 1 *or* 2 *as well.*

*Proof.* Every HF set is obtained by the rules of Definition 2.1. If $a_1$ is a HF set and we apply Rule 2.1.2 we obtain a list $L$ of HF sets which is of the form $a_1, a_2, ..., a_k$ where $k = 1$. Let $L$ be a list of HF sets of the form $a_1, a_2, ..., a_{k-1}$ and let $a_k$ be a HF set. If we apply Rule 2.1.3 we obtain a list of HF sets $La_k$ which is of the form $a_1, a_2, ... , a_{k-1}, a_k$. And finally, if $L$ is a list of HF sets of the form $a_1, a_2, ..., a_k$ and we apply Rule 2.1.4 we obtain a HF set $\{L\}$ which is of the form $\{a_1, a_2, ..., a_k\}$. □

The elements of a HF set are not necessarily all different from each other. However, a HF set should be independent of the number of repetitions of the same element. Therefore the words $\{a, a\}$ and $\{a\}$ should be equal HF sets, although these HF sets are two different words. A HF set should also be independent of the order of appearance of its elements. Hence the words $\{a, b\}$ and $\{b, a\}$ should be equal HF sets. Therefore

**2.7 Definition** Let $L_1$ and $L_2$ be two lists of HF sets. If for every HF set $x$ holds that $x$ is a term of $L_1$ precisely when $x$ is a term of $L_2$, then $L_1$ and $L_2$ are *equal lists* of HF sets. HF sets $A_1$ and $A_2$ are *equal* HF sets if $A_1 = A_2 = \{\}$, or if $A_1 = \{L_1\}$ and $A_2 = \{L_2\}$ where $L_1$ and $L_2$ are equal lists of HF sets.

The relations "are equal lists of HF sets" and "are equal HF sets" are 2-place computable relations. It is an easy to prove theorem that these two relations are equivalence relations. Now we can show that the Axiom of Extensionality holds for HF sets.

**2.8 Theorem** *Let* $A_1$ *and* $A_1$ *be two* HF *sets. If for every* HF *set x holds that x is an element of* $A_1$ *precisely when x is an element of* $A_2$*, then* $A_1$ *and* $A_2$ *are equal* HF *sets.*

*Proof.* If $A_1 = A_2 = \{\}$, then the theorem obviously holds. Otherwise, by Definition 2.7, we have $A_1 = \{L_1\}$ and $A_2 = \{L_2\}$ where $L_1$ and $L_2$ are *equal* lists of HF sets. In this case the theorem is a consequence of Definitions 2.5 and 2.7. □

We shall denote the equality of HF sets and the equality of words by the same sign =. We can afford this, since every HF set is a uniquely determined word if we lexicographically order it, which we do as follows. First it is obvious that for every word $x$ holds $x = x_1x_2...x_k$ where $x_i$ is a character for $i = 1, 2, ..., k$. Here $k$ can be equal to 1 or 2 as well. The number $k$ is the *length* of the word $x$ denoted by $l(x)$. The length of a word is a computable function on words. In the sequel we shall limit ourselves only to words made of characters { and }, since we shall have no need for any other words.

**2.9 Definition** A word $x$ is *smaller* than a word $y$ ($x < y$) if $l(x) < l(y)$, or if $x = x_1x_2...x_m$ and $y = y_1y_2...y_m$, where $x_i$ and $y_i$ are characters { or } for $i = 1, 2, ..., m$, and if there exists such $k$ that $x_i = y_i$ for $i < k$ and $x_k < y_k$ where { is smaller than }. In this case $y$ is *greater* than $x$ ($y > x$). Let $L = a_1, a_2, ... , a_k$ be a list of HF sets. The terms *of L follow each other lexicographically* if $i < j$ implies $a_i < a_j$ for $i, j = 1, 2, ..., k$. A list $L$ is *lexicographically ordered* (l.o.) if every term of $L$ is l.o. and if the terms of $L$ follow each other lexicographically. A HF set $A$ is l.o. if $A = \{\}$ or if $A = \{L\}$ where $L$ is a l.o. list of HF sets.

The relations < and > defined on arbitrary words are extensions of the relations < and > defined on natural numbers. It is not difficult to show that for two arbitrary words $x$ and $y$ precisely one of the following possibilities holds: $x < y$ or $x > y$ or $x = y$.

**2.10 Theorem** *Every HF set can be lexicographically ordered.*

*Proof*. Every HF set $A$ can be constructed by one of the rules of Definition 2.1. This can be done in such a way that $A$ is l.o.. If we apply Rule 2.1.1 we obtain the HF set {} which is l.o.. If we apply Rule 2.1.2 and $x$ is a l.o. HF set, then $x$ is a l.o. list of HF sets. If we apply Rule 2.1.3 and $L$ is a l.o. list of HF sets and $x$ is a l.o. HF set, then there are three possibilities. If $x$ is a term of $L$, then $Lx$ and $L$ are equal lists by Definition 2.7. Otherwise, if no term of $L$ is greater than $x$, then $Lx$ is l.o.. Otherwise, if no term of $L$ is smaller than $x$, then $Lx$ and $xL$ are equal lists and $xL$ is a l.o. list. The only remaining case is $L = L_1L_2$ where $L_1$ and $L_2$ are l.o. lists, each term of $L_1$ is smaller than $x$ and each term of $L_2$ is greater than $x$. Then the lists $Lx$ and $L_1xL_2$ are equal and $L_1xL_2$ is a l.o. list. Thus we always obtain a l.o. list. If we apply Rule 2.1.4 and $L$ is a l.o. list, we obtain the HF set {L} which is l.o.. Since any HF set can be obtained in the described way, any HF set can be l.o.. □

**2.11 Theorem** *Two equal lexicographically ordered HF sets are equal words.*

*Proof*. Let $A$ be the smallest l.o. HF set which is equal to but not the same word as another l.o. HF set $A'$. Let $A = \{a_1, a_2, ..., a_m)$, $A' = \{a_1', a_2', ..., a_n')$, and $m \leq n$. Let $k$ be such natural number that $a_i = a_i'$ for $0 < i < k$ but $a_k \neq a_k'$. If $k > m$, then $k = m + 1$, since there is no $(m + 1)$th element in $A$. In this case $a_{m+1}'$ is not an element of $A$ and therefore $A$ and $A'$ are not equal sets. Hence $k \leq m$. Let $a_k < a_k'$. Since $A$ is l.o., we have $a_k > a_i$ for $0 < i < k$. Hence $a_k > a_i'$ for $0 < i < k$ and $a_k < a_i'$ for $k \leq i \leq n$. Therefore $a_k \notin A'$. If, however $a_k' < a_k$, then we conclude in the same way that $a_k' \notin A$. In both cases $A$ and $A'$ are not equal HF sets, which contradicts the assumption. Hence $k$ does not exist, $m = n$, and $a_i = a_i'$ for $0 < i \leq m$. Since $A$ is the smallest HF set which can be expressed in l.o. form by two different words, and since every element of $A$ is smaller than $A$, the word $a_i$ is the same as the word $a_i'$ for $0 < i \leq m$. Therefore $A$ and $A'$ are equal words which contradicts the assumption. □

If $\{a_1, a_2, ..., a_m\}$ is a l.o. HF set, then $m$ is the *number of the elements* of this HF set. If two HF sets are equal then they have equal number of elements. Now we shall introduce the concept of *rank* of a HF set. For this reason we need to define first the function max. By the way we shall define also the function min which we shall need later.

**2.12 Definition** Let $A$ be a HF set. Then min($A$) denotes such element $m$ of $A$ for which holds $m \leq a$ for every $a \in A$ and max($A$) denotes such element $M$ of $A$ for which holds $a \leq M$ for every $a \in A$. The definition of min($L$) and max($L$), where $L$ is a list of HF sets, is almost exactly the same. We only replace the expression "element of a HF set" with the expression "term of a list of HF sets". The *rank* of $A$, which we denote by $r(A)$, is the natural number for which the following holds:

1) $r(\{\}) = 0$
2) If $L$ is a list of HF sets and $x$ is a HF set, then $r(Lx) = \max\{r(L), r(x)\}$
3) If $L$ is a list of HF sets, then $r(\{L\}) = r(L) + 1$

The functions min and max are 1-place computable functions, and rank is a 1-place computable function which transforms every HF set and every list of HF sets to a natural number. Namely, if we apply Rules 2.1.1, 2.1.3 or 2.1.4, then we compute the rank by Rules 2.12.1, 2.12.2 or 2.12.3 respectively. The application of Rule 2.1.2 preserves the rank. Rank, defined in 2.12, corresponds to the rank of hereditarily finite sets which we used in the first section.

Let the expression $A = \{x \mid \mathbf{P}(x)\}$, for the time being, mean that the elements of some HF set $A$ are precisely all those HF sets $x$ for which some computable predicate **P** is true.

**2.13 Theorem** *If A is a nonempty* HF *set, then* $r(A) = \max\{r(a) \mid a \in A\} + 1$.

*Proof.* Since every nonempty HF set $A$ can be expressed as $\{a_1a_2...a_k\}$, where $a_i$ is a HF set for $i = 1, 2, ... k$ (Theorem 2.6), and since $r(\{a_1a_2...a_k\}) = r(a_1a_2...a_k) + 1$ by Rule 2.12.3, the theorem is equivalent to the assertion that $r(a_1a_2...a_k) = \max(r(a_1)r(a_1)...r(a_k))$. We shall prove this assertion by induction on the *length* of the list $a_1a_2...a_k$. The assertion holds for every list of the length 1. Let the assertion hold for every list of the length $k$. Then $r(a_1a_2...a_ka_{k+1}) = \max(r(a_1a_2...a_k), r(a_{k+1}))$ by Rule 2.12.2. By assertion this is equal to

$$\max(\max(r(a_1)r(a_2)...r(a_k)), r(a_{k+1}))) = \max(r(a_1)r(a_2)...r(a_k)r(a_{k+1})).$$

Hence the assertion holds for $k + 1$ if it holds for $k$. Since it holds for $k = 1$, it holds for every $k > 1$. □

Now we shall introduce the concept of *approximation* of a HF set and show some of its properties.

**2.14 Definition** Let $F(n, x)$ be the following 2-place computable function where $n$ is a natural number and $x$ is a HF set or a list of HF sets.

1) $F(n, \{\}) = \{\}$ for every $n$.
2) If $L$ is a list of HF sets and $x$ is a HF set, then $F(n, Lx) = F(n, L)F(n, x)$.
3) If $L$ is a list of HF sets, then $F(n + 1, \{L\}) = \{F(n, L)\}$.
4) $F(0, x) = \{\}$ for every HF set $x$.

**2.15 Example** $F(2, \{\{\}, \{\{\}\}, \{\{\{\}\}\}, \{\{\{\{\}\}\}\}\})$ =

$$\{F(1, \{\}), F(1, \{\{\}\}), F(1, \{\{\{\}\}\}), F(1, \{\{\{\{\}\}\}\})\} =$$

$$\{\{\}, \{F(0, \{\})\}, \{F(0, \{\{\}\})\}, \{F(0, \{\{\{\}\}\})\}\} = \{\{\}, \{\{\}\}, \{\{\}\}, \{\{\}\}\} = \{\{\}, \{\{\}\}\}$$

**2.16 Theorem** *If A is a* HF *set, then* $F(n, A)$ *is a* HF *set and* $r(F(n, A)) \leq n$ *for every n.*

*Proof.* We shall prove the theorem by induction on $r(A)$. The theorem holds for $r(A) = 0$, that is for $A = \{\}$, by Rule 2.14.1. Let the theorem hold for all HF sets of the rank not greater than $\rho$ and let $A$ be a HF set of the rank $\rho + 1$. If $n = 0$, then the theorem holds because of Rule 2.14.4. If $n = m + 1$, then

$$F(n, A) = F(m + 1, A) = F(m + 1, \{a_1a_2...a_k\}) = \{F(m, a_1a_2...a_k)\} = \{F(m, a_1)F(m, a_2)...F(m, a_k)\}$$

by Theorem 2.6, by Rule 2.14.3 and by $k - 1$ times applying Rule 2.14.2. The rank of every $a_i$, where $i = 1, 2, ..., k$, is by Theorem 2.13 not greater than $\rho$. Hence every $F(m, a_i)$ is a HF set and by assumption $r(F(m, a_i)) \leq m$. Therefore $F(n, A)$ is a HF set by Definition 2.1 and

$$r(F(n, A)) = r(\{F(m, a_1)F(m, a_2)...F(m, a_k)\}) \leq m + 1 = n$$

by Theorem 2.13. Thus, by induction on $n$, we have $r(F(n, A)) \leq n$ for every $n$ and every $A$ of the rank $\rho + 1$. Therefore the theorem holds for a HF set $A$ of every rank $\rho$ and for every $n$. □

**2.17 Theorem** *For every* HF *set A and every n holds* $F(n + 1, A) = \{F(n, a) \mid a \in A\}$.

*Proof.* The theorem obviously holds for $A = \{\}$. So let $A = \{a_1a_2...a_k\}$. Then we have $F(n + 1, A) = \{F(n, a_1)F(n, a_2)...F(n, a_k)\}$ (See the proof of the previous theorem.). Since $F(n, a)$ is a HF set for every HF set $a$, by the previous theorem, and since every element of $F(n + 1, A)$ is of the form $F(n, a)$ for some $a \in A$, the theorem holds. □

**2.18 Theorem** *$F(r(A) + m, A) = A$ for every* HF *set A and every m.*

The function $F$ defined by 2.14 operates on HF sets in the same way as the function $F$ defined by 1.1 because of Rule 2.14.4 and Theorem 2.17. And the rank defined by 2.12 operates on HF sets in the same way as the rank we have used in the first section. Therefore the proof of the theorem is the same as the proof of Theorem 1.3. □

The function $F(n, A)$ is recursive. It operates on $A$ in such a way that it operates on each of it elements. If no membership chain in $A$ is longer than $n$, then $F$ does not transform $A$ at all. If, however, there exists some $x_{n+1}$ for which holds $x_{n+1} \in x_n \in .... \in x_2 \in x_1 \in A$ then $F$ deletes this $x_{n+1}$ from the word $A$. The greater that $n$ is the less characters deletes $F$ from $A$ and the more precise approximation of $A$ is $F(n, A)$. Because of this we call $F(n, A)$ the *nth approximation* of the HF set $A$. Let's prove that an approximation of an approximation equals to the less precise of both approximations.

**2.19 Theorem** *For every n, m and every* HF *set A holds $F(n, F(n+m, A)) = F(n+m, F(n, A)) = F(n, A)$.*

*Proof.* Since $r(F(n, A)) \leq n$ by Theorem 2.16, we have $F(n + m, F(n, A)) = F(n, A)$ by Theorem 2.18. We shall prove the validity of $F(n, F(n + m, A)) = F(n, A)$ by induction on $n$. If $n = 0$ then, by Rule 2.14.4, we have $F(n, F(m, A)) = F(n, A)$ for every $A$ and $m$. Let the equation hold for some $n$. Then by Theorem 2.17 we have

$$F(n + 1, F(n + 1 + m, A)) = F(n + 1, \{F(n + m, x) \mid x \in A\}) =$$

$$\{F(n, F(n + m, x)) \mid x \in A\} = \{F(n, x) \mid x \in A\} = F(n + 1, A).$$

Thus, if the equation holds for $n$, then it holds for $n + 1$. Since it holds for $n = 0$, it holds for every $n$. □

**2.20 Corollary** *For every n and every* HF *set A holds $F(n, F(n + 1, A)) = F(n, F(n, A)) = F(n, A)$.* □

By Theorem 2.16 we have $r(F(n, A)) \leq n$. Now we shall prove a stronger result.

**2.21 Theorem** *If $r(F(n, A)) < n$, then $F(n, A) = A$.*

*Proof.* We shall prove the theorem by induction on $n$. The smallest $n$ for which the condition of the theorem can be fulfilled is 1 so let $r(F(1, A)) < 1$. Then $r(F(1, A)) = 0$ and $F(1, A) = \{\}$. Since $A \neq \{\}$ implies $F(1, A) = \{\{\}\}$, we have $A = \{\}$. Thus $r(F(1, A)) < 1$ implies $F(1, A) = A$ and the theorem holds for $n = 1$. Let the theorem hold for $n$ and let $r(F(n + 1, A)) < n + 1$ hold for some $A$. Then $r(\{F(n, x) \mid x \in A\}) < n + 1$ by Theorem 2.17. Therefore $r(F(n, x)) < n$ for every $x \in A$ by Theorem 2.13 and thus $F(n, x) = x$ for every $x \in A$ by assumption. This implies $F(n + 1, A) = \{F(n, x) \mid x \in A\} = \{x \mid x \in A\} = A$. Thus the theorem holds for every $n$. □

Now we shall prove that the Principle of Set Equality and the Principle of Membership hold for HF sets without any additional axioms.

**2.22 Theorem** *A* HF *set A is equal to a* HF set *B precisely when $F(n, A) = F(n, B)$ holds for every n.*

*Proof.* Let $A$ and $B$ be two HF sets. If $A = B$, then for every $n$ certainly holds $F(n, A) = F(n, B)$. Let now for every $n$ hold $F(n, A) = F(n, B)$ and let $r(A) \geq r(B)$. Then $A = F(r(A), A) = F(r(A), B) = B$ by Theorem 2.18. □

**2.23 Theorem** *A* HF *set x is an element of a* HF *set A precisely when $F(n, x) \in F(n + 1, A)$ holds for every n.*

*Proof.* Let $A$ be a HF set and $x \in A$. Since $F(n + 1, A) = \{F(n, y) \mid y \in A\}$ for every $n$ by Theorem 2.17, we have $F(n, x) \in F(n + 1, A)$ for every $n$. Let now $F(n, x) \in F(n + 1, A)$ hold for every $n$. Then we have $A = F(r(A), A)$ by Theorem 2.18, and we have $F(r(A) - 1, x) = x$ by Theorem 2.13. Since $F(r(A) - 1, x) \in F(r(A), A)$, we have $x \in A$. □

# 3. Set theory based on Approximations of Sets

If the $n$th approximation of a HF set is known for some $n > 0$, we can compute its $(n - 1)$th approximation by Corollary 2.20. However, in general we cannot compute its $(n + 1)$th approximation. If, for example, $F(1, A) = \{\{\}\}$, then $F(0, A) = \{\}$, while $F(2, A)$ can be either $\{\{\}\}$ or $\{\{\}\{\{\}\}\}$. We could say that an approximation has precisely one predecessor, but it can have several successors. By Theorems 2.16 and 2.21 a $n$th approximation can be any HF set of the rank not greater than $n$. It is known that the number of HF sets of the rank not greater than $n$ grows very rapidly by $n$. Therefore at least one $n$th approximation has as many successors as we wish if only $n$ is large enough. On the other hand for every $n$ a lot of $n$th approximations have only one successor. By Theorem 2.21 this is true for every $n$th approximation whose rank not greater than $n - 1$.

If we connect every $n$th approximation with its successors, we obtain a tree whose levels from 0 to 3 are drawn in the following picture. In order that the root { } will be at the top of the tree, as it is usual, the picture has to be turned round for 90 degrees.

```
{}
├─{}
│ └─{}
│   └─{}
└─{{}}
  ├─{{}}
  │ └─{{}}
  ├─{{{}}}
  │ ├─{{{}}}
  │ ├─{{{{}}}}
  │ ├─{{{},{{}}}}
  │ ├─{{{}},{{{}}}}
  │ ├─{{{}},{{},{{}}}}
  │ ├─{{{{}}},{{},{{}}}}
  │ └─{{{}},{{{}}},{{},{{}}}}
  └─{{},{{}}}
    ├─{{},{{}}}
    ├─{{},{{{}}}}
    ├─{{},{{},{{}}}}
    ├─{{},{{}},{{{}}}}
    ├─{{},{{}},{{},{{}}}}
    ├─{{},{{{}}},{{},{{}}}}
    └─{{},{{}},{{{}}},{{},{{}}}}
```

We shall call this tree the *universal tree* **U**. The reasons for the name will be explained in the following section. The nodes on the $n$th level of **U** are precisely all HF sets of the rank not greater than $n$. Now we shall define **U** formally. The definition will be based on the concept of node. We must take into account that different nodes on different levels of **U** can represent the same HF set. This is also evident from the picture.

**3.1 Definition** The *ordered pair* of HF sets $a$ and $b$ is the HF set $\{\{a\}, \{a, b\}\}$ denoted by $(a, b)$. A *node* is every ordered pair $(n, A_n)$, where $n$ is a natural number called *level* and $A_n$ is a HF set of the rank not greater than $n$. The *successor* of the node $(n, A_n)$ is a node $(n + 1, A_{n+1})$ where $F(n, A_{n+1}) = A_n$. A node $v$ is a *predecessor* of a node $v'$ precisely when $v'$ is a successor of $v$. The *universal tree* is the predicate $\mathbf{U}(v, v')$ which is a true statement precisely when the node $v'$ is a successor of the node $v$. The node $(0, \{\})$ is the *root* of **U**. Every successor of a node $v$ is a *descendant* of $v$. A descendant of a descendant of $v$ is a descendant of $v$. The predecessor of $v$ is an *ancestor* of $v$. An ancestor of an ancestor of $v$ is an ancestor of $v$.

The definition implies that every node $(n, A_n)$ has at least one successor, namely $(n + 1, A_n)$. If a node $(n, A_n)$ is not the root of **U**, it has precisely one predecessor, namely $(n - 1, F(n - 1, A_n))$. The predicate $\mathbf{U}(v, v')$ is computable, since there certainly exists an algorithm which checks whether $v$ and $v'$ are nodes and whether $v$ is the predecessor of $v'$.

Let's take an arbitrary HF set $A$ and let's calculate $F(n, A)$ for every $n$. This calculation determines some path in **U** which consists of nodes $(n, F(n, A))$ where $n$ is any natural number. The path starts in the root of **U** and continues without an end. We formally define path as follows.

**3.2 Definition** A *path in* **U** is a predicate $\mathbf{P}(v, v')$ for which the following holds

1) $\mathbf{P}(v, v')$ implies $\mathbf{U}(v, v')$. (**P** consists only of nodes and their successors in **U**.)
2) There exists such $v$ that $\mathbf{P}((0, \{\}), v)$ is true. (**P** starts in the root of **U**.)
3) If $\mathbf{P}(v, v')$, then there is only one such $v''$ that $\mathbf{P}(v', v'')$. (Every node of **P** has only one successor.)

If $v = (n, A_n)$ and $\mathbf{P}(v, v')$ is true for some path **P** then $v$ is the $n$th node of **P**.

This definition implies that no path ends. The approximations represented by the nodes of some path become larger and larger words as the level of the nodes increases. If these words stop expanding and are from some level on all equal to each other, then the path represents some HF set. Immediately there arises the question what represent all other paths. We obtain such a path if we put, for example, $A_{n+1} = \{A_n\}$, that is $A_n = n$, where $v_n = (n, A_n)$ is the $n$th node of the path and $n$ is a natural number in Zermelo form.

A path is in general not a computable predicate as is the universal tree. We can determine the nodes of some path with an algorithm or we can choose them in a completely random way. The only important thing is that these nodes satisfy Rules 3.2.

To every set which exists in ZF corresponds some path in **U**. Namely, for every such set we can calculate the approximations $F(0, A)$, $F(1, A)$, $F(2, A)$, $F(3, A)$, ... by Definition 1.1. Every of these approximations is a HF set, and $F(n, F(n + 1, A)) = F(n, A)$ must hold for every $n$ by Corollary 2.20. However, as we have already found out in the first section, one such path can correspond to different sets in ZF. We shall prevent this on the base of the Principle of Set Equality. In our set theory one path will correspond only to one set. On the other hand, some paths in **U** do not correspond to any set in ZF. One of these is the described path where $A_{n+1} = \{A_n\}$. If there exists a set $A$, which is determined by this path, then this set is not regular, since there must exists an infinite membership chain $A \ni x_1 \ni x_2 \ni x_3 \ni x_4 \ni x_5$ .... In spite of that we shall extend the concept of set in such a way that some set will correspond to such a path as well.

**Axiom of Set Existence** *Every set can be represented by a unique path in* **U** *and every path in* **U** *represents a unique set.*

Now we can extend the concept of approximation of a set and the domain of the function $F$ to arbitrary sets $A$.

**3.3 Definition** If $(n, A_n)$ is the $n$th node of the path that represents a set $A$, then the HF set $A_n$ is the $n$th *approximation* of the set $A$ and $F(n, A) = A_n$.

The function $F$ is here defined in the environment of predicate calculus as a rule and not in the environment of set theory as a set of ordered pairs of elements of some domain and range. However, when we shall prove the existence of the ordered pair of two arbitrary sets and the existence of the cartesian product of two arbitrary sets, we shall extend the domain and the range of $F$ once again and then $F$ will really be a set. The function $F$ is in general not computable, since it is defined for every set, that is, for every path in **U** and there exist paths which are not computable by any algorithm.

**3.4 Theorem** *Sets A and B are equal precisely when* $F(n, A) = F(n, B)$ *holds for every n.*

*Proof.* If $A = B$, then $F(n, A) = F(n, B)$ holds for every $n$, since $A$ and $B$ are, by the Axiom of Set Existence, the same set represented by the same path. Let now $F(n, A) = F(n, B)$ hold for every $n$. Then the corresponding paths for $A$ and $B$ in **U** are equal. Therefore $A$ and $B$ are equal by the Axiom of Set Existence. □

We shall denote the set theory in which every set is determined by its approximations by AS (Approximations of Sets). From now on this theory will be our working environment. If we shall work in the environment of ZF, we shall stress it explicitly.

**3.5 Theorem**. *A set A exists precisely when for every n holds* $F(n, F(n + 1, A)) = F(n, A)$.

*Proof.* If $A$ exists then the condition in the theorem is certainly fulfilled, since the path for the set $A$ and the path for the HF set $F(n + 1, A)$ coincide to the $(n + 1)$th node, and since we can apply Corollary 2.20 on the HF set $F(n + 1, A)$. On the other hand, if this condition is fulfilled, then the node $(n + 1, F(n + 1, A))$ is a successor of the node $(n, F(n, A))$ for every $n$ and hence the nodes, determined by the approximations of $A$, form a path in **U**. Thus $A$ exists by the Axiom of Set Existence. □

**3.6 Definition** The function $f_A$: $n \rightarrow F(n, A)$, which maps every $n$ to the $n$th approximation of $A$, is the *developing sequence* of $A$. We write this sequence symbolically as $F(0, A)$, $F(1, A)$, $F(2, A)$, ....

The sequence $F(0, A)$, $F(1, A)$, $F(2, A)$, ... is called the developing sequence of the set $A$ because it shows the development of the set $A$ from its most rough approximation $F(0, A)$, which is equal to $\{\}$, through better and better approximations which tend towards the set $A$ if $n$ approaches infinity, or even before that if $A$ is a HF set. We shall give a rigorous proof of this at the end of the article.

**3.7 Theorem** *The developing sequence of any set determines a unique path in* **U** *and any path in* **U** *determines the developing sequence of a unique set.*

*Proof.* A simple consequence of the Axiom of Set Existence and Theorem 3.4. □

The concept of element of a set can be introduced in AS by the Principle of Membership which holds for HF sets anyway. We could do this in a form of a new axiom, but since we want to keep the number of axioms as small as possible, and since the concept of element of a set has not yet been defined, we can introduce this concept by a definition.

**3.8 Definition** A set $x$ is an *element* of a set $A$ $(x \in A)$ if $F(n, x) \in F(n + 1, A)$ holds for every $n$. The denotation $x \notin A$ is an abbreviation for $\neg (x \in A)$.

We must emphasize that there is an essential difference between the relations $A = B$ and $x \in A$ defined on HF sets and the same two relations defined on arbitrary sets although we used the same denotations in both cases. The first two relations are decidable in a finite number of steps and so they are computable 2-place predicates on HF sets. But there exists no algorithm which would in a finite number of steps decide whether any of the second two relations is true for two arbitrary sets. Namely, if we want to prove $A = B$, then we must decide whether $F(n, A) = F(n, B)$ holds for every $n$, and if we want to prove $x \in A$, then we must decide whether $F(n, x) \in F(n + 1, A)$ holds for every $n$.

Now we shall show that the extended function $F(n, A)$ ddefined in 3.3 has the same properties as the function $F(n, A)$ defined in 2.14. We must prove these properties once again, since we have proved them in the previous section only for HF sets. Maybe there is a second thought that in this way the properties of $F$ will be proved twice and that we should have avoided this redundancy by postponing the proofs till the present section. But in this case we could not have constructed the universal tree **U** which is the base of AS and we could not have proved that the Principle of Set Equality and the Principle of Membership hold for HF sets without any additional axiom.

**3.9 Theorem** *For every n and every set A holds* $r(F(n, A)) \leq n$.

*Proof*. By Theorem 3.5 we have $F(n, A) = F(n, F(n + 1, A))$ for every $n$ and every set $A$. Therefore $r(F(n, A)) = r(F(n, F(n + 1, A))) \leq n$ by Theorem 2.16. □

**3.10 Theorem** *For every set A and every m holds* $F(n, F(n + m, A)) = F(n + m, F(n, A)) = F(n, A)$.

*Proof*. Since $r(F(n, A)) \leq n$ by Theorem 3.9, we have $F(n + m, F(n, A)) = F(n, A)$ by Theorem 2.18. In order to prove $F(n, F(n + m, A)) = F(n, A)$ we infer as follows. If sets $A$ and $B$ have the same $n$th approximation, then the paths that represent them in **U** coincide from the root of **U** up to the $n$th node. Thus their $m$th approximations for $m < n$ must also be the same. Hence $F(m, A) = F(m, B)$ holds for every $m \leq n$. The sets $F(n + m, A)$ and $A$ have the same $(n + m)$th approximation. Therefore we have $F(n, F(n + m, A)) = F(n, A)$. □

**3.11 Corollary** *For every n and every set A holds* $F(n, F(n, A)) = F(n, A)$. □

Now we shall prove a very important property of the function $F(n, A)$, namely that for every set $A$ and every $n$ holds $F(n + 1, A) = \{F(n, x) \mid x \in A\}$. First we shall prove:

**3.12 Lemma** *If* $F(n, a) \in F(n + 1, A)$, *then* $F(m, a) \in F(m + 1, A)$ *for* $m \leq n$.

*Proof*. Let $F(n, a) \in F(n + 1, A)$. By Theorem 2.6 we can write $F(n + 1, A) = \{x_1x_2...F(n, a)...x_k\}$. Since $m \leq n$, we have, by Theorem 3.10, by Rule 2.14.3 and by $k - 1$ times application of Rule 2.14.2,

$$F(m + 1, A) = F(m + 1, F(n + 1, A)) = F(m + 1, \{x_1x_2...F(n, a)...x_k\}) =$$

$\{F(m, x_1)F(m, x_2)...F(m, F(n, a))...F(m, x_k)\}$. Thus $F(m, a) = F(m, F(n, a)) \in F(m + 1, A)$. □

**3.13 Theorem** *For every set A and every n holds* $F(0, A) = \{\}$ *and* $F(n + 1, A) = \{F(n, a) \mid a \in A\}$.

*Proof*. The fact $F(0, A) = \{\}$ follows from Theorem 3.9. Certainly $\{F(n, a) \mid a \in A\} \subseteq F(n + 1, A)$ holds for every $n$, since $a \in A$ implies $F(n, a) \in F(n + 1, A)$ for every $n$. Does the set $F(n + 1, A)$ contain something else? Let some HF set $a_n$ be an element of $F(n + 1, A)$. Then there exists $a_{n+1} \in F(n + 2, A)$ such that $F(n, a_{n+1}) = a_n$. Namely, let $F(n + 1, A) = \{x_1x_2...a_n...x_i\}$ and let $F(n + 2, A) = \{y_1y_2...y_j\}$. By Theorem 3.5 and by Rules 2.14.2 and 2.14.3 we have

$$F(n + 1, A) = F(n + 1, F(n + 2, A)) = F(n + 1, \{y_1y_2...y_j\}) = \{F(n, y_1)F(n, y_2)...F(n, y_j)\}.$$

Hence some $F(n, y_k)$, $1 \leq k \leq j$, must be equal to $a_n$. By putting $a_{n+1} = y_k$ the assertion is proved. In a similar way we conclude that there exists such $a_{n+2} \in F(n + 3, A)$ that $F(n + 1, a_{n+2}) = a_{n+1}$ holds. This implies that there exists a sequence $a_n, a_{n+1}, a_{n+2}, a_{n+3}, ...$, for which holds $a_{n+m} \in F(n + m + 1, A)$ and $F(n + m, a_{n+m+1}) = a_{n+m}$ for every $m$. Since $a_n \in F(n + 1, A)$, we have $F(m, a_n) \in F(m + 1, A)$ for $m < n$ by Lemma 3.12. Let's put $a_m = F(m, a_n)$ for $m < n$. Then $a_0, a_1, a_2, ..., a_n, a_{n+1}, a_{n+2}, a_{n+3}, ...$ is the developing sequence for some set $a$. Since $F(m, a) = a_m \in F(m + 1, A)$ holds for every $m$, we have $a \in A$. Therefore, if $a_n \in F(n + 1, A)$, then there exists such $a \in A$ that $F(n, A) = a_n$ holds for every $n$. This implies that $F(n + 1, A) \subseteq \{F(n, a) \mid a \in A\}$ holds for every $n$. Hence we can conclude that for every $n$ holds $F(n + 1, A) = \{F(n, a) \mid a \in A\}$. □

The following theorem is an extension of Theorem 2.8 from HF sets to arbitrary sets. It proves the validity of the Axiom of Extensionality for arbitrary sets.

**3.14 Theorem** *Let A and B be two sets. If for every set x holds that x is an element of A precisely when x is an element of B, then A and B are equal sets.*

*Proof.* Let $x \in A \Leftrightarrow x \in B$ for every set $x$. Then by Definition 3.8 we have $F(n, x) \in F(n + 1, A)$ precisely when $F(n, x) \in F(n + 1, B)$ for every $n$ and every $x$. Since, by Theorem 3.13, the elements of $(n + 1)$th approximation of any set are precisely all the $n$th approximations of the elements of this set, it must be true that $F(n + 1, A) = F(n + 1, B)$ for every $n$. And since $F(0, A) = \{\} = F(0, B)$, the sets $A$ and $B$ are equal by Theorem 3.4. □

Now we shall show how we can find out from the developing sequence of some set whether this set is a HF set. First we have

**3.15 Theorem** *If $r(F(n, A)) < n$, then $F(n, A) = A$.*

*Proof.* The proof of this theorem is the same as the proof of Theorem 2.21 with the only difference that we use Theorem 3.13 instead of Theorem 2.17 which holds only for HF sets. □

**3.16 Theorem** *Some set A is a* HF *set precisely when there exists such n that $F(n, A) = F(n + 1, A)$.*

*Proof.* If $A$ is a HF set, we put $n = r(A)$ and then $F(n, A) = F(n + 1, A)$ holds by Theorem 2.18. If, however, there exists such $n$ that $F(n, A) = F(n + 1, A)$, then $r(F(n + 1, A)) \leq n < n + 1$ by Theorem 3.9. Hence $F(n + 1, A) = A$ by Theorem 3.15 and consequently $A$ is a HF set. □

Thus we have proved all important properties of the function $F(n, A)$ where $A$ is an arbitrary set.

It is interesting that every set can not only be interpreted as a path in **U** but also as a subtree of **U**.

**3.17 Definition** A *subtree* of **U** or a *tree* in **U** is a predicate **T**($v$, $v'$) which satisfies the same conditions as the predicate path in Definition 3.2 with the only difference that the word *only* in the Condition 3.2.3 is replaced with the words *at least*. That is, every node of a path has only one successor while every node of a subtree of **U** has at least one successor.

**3.18 Theorem** *Every set is represented by a unique (possibly empty) subtree of* **U** *and every subtree of* **U** *represents a unique set. The paths in a subtree of* **U** *represent the elements of the corresponding set.*

*Proof.* If $A$ is the empty set, then the corresponding tree is empty as well. Let $A$ be a nonempty set. By the Axiom of Set Existence every $x \in A$ corresponds to a uniquely determined path. All these paths form some tree **T**. Namely, let **T**($v$, $v'$) hold precisely when **P**($v$, $v'$) holds for some path **P** which represents some $x \in A$. We shall now prove that each of the three conditions for the existence of **T** is fulfilled.

1) If $\mathbf{T}(v, v')$ is true, there exists such a path $\mathbf{P}$ that $\mathbf{P}(v, v')$ is true. This implies that $\mathbf{U}(v, v')$ is true.
2) Since $A$ is not empty, there exists $x \in A$. Therefore there exists the path $\mathbf{P}$ for this $x$. For $\mathbf{P}$ holds $\mathbf{P}((0, \{\}), F(1, x))$ by Condition 3.2.2 and hence $\mathbf{T}((0, \{\}), F(1, x))$ is true.
3) If $\mathbf{T}(v, v')$ is true, then there exists such a path $\mathbf{P}$ that $\mathbf{P}(v, v')$ is true. Then there exists such a node $v''$ that $\mathbf{P}(v', v'')$ is true. Hence $\mathbf{T}(v', v'')$ is true.

The tree $\mathbf{T}$ uniquely determines the set $A$. Namely, if two sets are different, then by the Axiom of Extensionality (Theorem 3.14) there exists at least one element of one of these two sets which is not an element of the other set. Therefore there exists at least one path in the tree of one set which is not in the tree of the other set. Hence these two trees are different. And vice versa is true as well. If there exists a tree, then a uniquely determined set corresponds to this tree. The elements of this set are precisely those sets whose paths are in this tree. □

**3.19 Definition** The tree that represents a set $A$ is the *developing tree* of $A$.

**3.20 Corollary** *A set exists precisely when its developing tree exists.* □

**3.21 Theorem** *The nodes on the nth level of the developing tree of a set A represent the elements of the* $(n + 1)$*th approximation of A. That is,* $(n, a)$ *is a node of the developing tree of A precisely when* $a \in F(n + 1, A)$.

*Proof.* Let $\mathbf{T}$ be the developing tree of some set $A$. A node on the $n$th level of $\mathbf{T}$ is the $n$th node of some path in $\mathbf{T}$. This path represents some set $x \in A$ (Theorem 3.18). Therefore the node represents the approximation $F(n, x)$. Since for every $x \in A$ there exists a path in $\mathbf{T}$, the $n$th approximation of any $x \in A$ is represented by some node on the $n$th level of $\mathbf{T}$. Hence the nodes on the $n$th level of $\mathbf{T}$ represent precisely all the elements of the HF set $\{F(n, x) \mid x \in A\}$. And this set is by Theorem 3.13 equal to $F(n + 1, A)$. □

# 4. A Comparison between ZF and AS

First we shall show some fundamental differences between ZF and AS. In ZF the set of all sets does not exist, since in the opposite case the Axiom of Separation would imply the existence of the set of all those sets which are not elements of themselves and this would imply *Russell's paradox*. Namely, such a set is an element of itself precisely when it is not an element of itself. However, in AS the set of all sets (the universal set) exists but in spite of that the argument for Russell's paradox fails.

**4.1 Theorem** *The set of all sets exists.*

*Proof.* By Theorem 3.18 and by Definition 3.19 every subtree of $\mathbf{U}$ is a developing tree of some set. Therefore also $\mathbf{U}$ is the developing tree of some set $U$. Since every set is represented by some path in $\mathbf{U}$ by the Axiom of Set Existence, every set is an element of $U$ by Theorem 3.18. And since every element of $U$ is some set, $U$ is the set of all sets. □

**4.2 Theorem** *The Axiom of Foundation is not consistent with AS.*

*Proof.* By the Axiom of Foundation no set is an element of itself. However, the set of all sets is already by definition an element of itself. As this set exists by the previous theorem, the Axiom of Foundation is not consistent with AS. □

Invalidity of the Axiom of Foundation is, however, not a restriction but an expansion of the universe of all sets.

**4.3 Theorem** *The argument for Russell's paradox fails in AS.*

*Proof*. Let the elements of some set $R$ be all sets which are not elements of themselves. Every HF set is certainly an element of $R$, since none of HF sets is an element of itself. If it were, then, by Theorem 2.13, its rank would not be a natural number. However, this is not the case. (See the observations following Definition 2.12.) But if every HF set is an element of $R$, then every set is an element of $R$. Namely, let $A$ be an arbitrary set. Since every HF set is an element of $R$, we have $F(n, A) \in R$ for every $n$. By Definition 3.8. we have $F(n, A) \in R$ precisely when $F(m, F(n, A)) \in F(m + 1, R)$ for some $n$ and every $m$. If we put $n = m$ and apply Corollary 3.11, then we obtain $F(n, A) \in F(n + 1, R)$. And this relation holds for every $n$, which, by Definition 3.8, means that $A \in R$. This implies $R = U$. As $U$ is the set of all sets, we have $U \in U$. Consequently, $R$ includes itself as an element and does not meet the conditions of Russell's paradox. □

Now we see again what we have already seen in section 1. In AS we cannot arbitrarily decide which objects are elements of some set and which are not. For example, there is no set whose elements would be precisely all HF set. Every set whose elements are all HF sets must be equal to the set of all sets. Generally, if a set $A$ exists, then the set of precisely those $x \in A$ for which some property **P** holds does not necessarily exist. Hence:

**4.4 Theorem** *The Axiom of Separation is not consistent with AS.* □

Since the Axiom of Replacement implies the Axiom of Separation, the Axiom of Replacement is not consistent with AS either. When exactly does the set of all $x \in U$, for which some predicate **P** is true, exist? In order to answer this question we must first find out which predicates exist in AS. As we have already told in the beginning of the 2. section, we allow only the computable predicates on words and we work in the environment of predicate calculus without the concept of domain, since domains are sets and all words do not form a set in AS. But there is no longer any need for this restriction. The domain of AS is, of course, the universal set $U$. However, this narrows our universe of all possible words to those words which are elements of $U$, that is to HF sets. Apart from sets we also introduced in the 3. section the function $F(n, A)$ which is a 3-place predicate and whose value exists for every natural number $n$ and every existent set $A$. Let's resume the sufficient conditions for the existence of a predicate in AS.

**4.5 Theorem** *If a predicate* **P** *is expressible in predicate calculus by computable predicates on* HF *sets, by the function F defined in* 3.3, *and by predicates which exist, then* **P** *exists as well.* □

**4.6 Example** Let $\mathbf{N}(x)$ be the predicate "$x$ is a natural number". This predicate is a computable predicate on HF sets. Let $A$ be a fixed set. Then the predicate $x \in A \Leftrightarrow \forall n\ (\mathbf{N}(n) \wedge F(n, x) \in F(n, A))$ and the predicate $x = A \Leftrightarrow \forall n\ (\mathbf{N}(n) \wedge F(n, x) = F(n, A))$ exist, since they are expressible in the described way. Notice that $F(n, x) \in F(n, A)$ and $F(n, x) \in F(n, A)$ are computable predicates on HF sets.

**4.7 Theorem**. *If* **P** *is a predicate with one free variable over U, then*

*a) the* HF *set* $\{F(n, x) \mid \mathbf{P}(x)\}$ *exists for every n*
*b) there exists a set P whose* $(n + 1)$*th approximation is equal to* $\{F(n, x) \mid \mathbf{P}(x)\}$ *for every n*
*c) if* $x \in P$, *then for every n there exists such* $y_n$ *that* $F(n, x) = F(n, y_n)$ *and* $\mathbf{P}(y_n)$ *is true*
*d) the set P is the smallest set of all* $x \in U$ *for which* $\mathbf{P}(x)$ *is true*
*e) the set* $\{x \mid \mathbf{P}(x)\}$ *exists precisely when* $x \in P$ *implies* $\mathbf{P}(x)$

*Proof.*

a) Since the set $U$ exists, its approximation $F(n, U)$ exists for every $n$. By Theorem 3.13 we have $F(n + 1, U) = \{F(n, x) \mid x \in U\}$. Let $n$ be fixed. As **P** has one free variable with domain $U$, the same is true for the predicate $\mathbf{P}(x) \wedge F(n, x) \in F(n+1, U)$. If this predicate is a true statement for some $x \in U$, then $F(n, x)$ is an element of the HF set $\{F(n, x) \mid \mathbf{P}(x)\}$ otherwise it is not. Therefore the HF set $\{F(n, x) \mid \mathbf{P}(x)\}$ exists.

b) By a) the set $P_n = \{F(n, x) \mid \mathbf{P}(x)\}$ exists for every $n$. We have $F(n, P_{n+1}) = F(n, \{F(n, x) \mid \mathbf{P}(x)\}) = \{F(n-1, F(n, x)) \mid \mathbf{P}(x)\} = \{F(n-1, x) \mid \mathbf{P}(x)\} = F(n, P_n)$. Hence $\{\}$, $P_1$, $P_2$, ... is the developing sequence of some set $P$ by Theorem 3.5 and hence $F(n+1, P) = P_{n+1} = \{F(n, x) \mid \mathbf{P}(x)\}$ for every $n$.

c) If $x \in P$, then $F(n, x) \in F(n+1, P)$ for every $n$ which, because of b), means that for every $n$ there exists such $y_n$, for which $\mathbf{P}(y_n)$ is true, that $F(n, y_n) = F(n, x)$ and $\mathbf{P}(y_n)$ is true.

d) If $\mathbf{P}(x)$ is true, then $F(n, x) \in F(n+1, P)$ for every $n$ by b). Consequently, $x \in P$. Let every $x$ for which $\mathbf{P}(x)$ is true be an element of a set $A$, and let $y \in P$. Then, by c), for every $n$ holds $F(n, y) = F(n, y_n)$ for some $y_n$ where $\mathbf{P}(y_n)$ is true. This implies that $y_n \in A_n$ and consequently $F(n, y) = F(n, y_n) \in F(n+1, A)$ for every $n$. Therefore $y \in A$ which implies $P \subseteq A$.

e) a consequence of d). □

Thus equipped we shall prove that all the axioms of ZF, with the exception of the Axioms of Separation, Replacement and Foundation, are in force in AS. We have already proved this for the Axiom of Extensionality (Theorem 3.14) and the Axiom of the Existence of empty set, since the set $\{\}$ exists. For the rest of ZF axioms we shall prove it now. We shall also prove that intersection, complement, and cartesian product exist in AS although complement has not quite the same properties as in ZF.

If we want to emphasize that some set is a set of sets, then we talk about a *family* of sets (of course, in AS there is no other alternative). The *intersection* of a family of sets $\mathcal{F}$ is the set of all those sets that are elements of every $A \in \mathcal{F}$. The intersection of $\mathcal{F}$ is denoted by $\cap\mathcal{F}$. If $\mathcal{F} = \{A_1, A_2, ..., A_k\}$ we denote $\cap\mathcal{F}$ by $A_1 \cap A_2 \cap ... \cap A_k$.

**4.8 Theorem** *The intersection of an arbitrary family of sets exists.*

*Proof*. The corresponding predicate for $\cap\mathcal{F}$ is $\mathbf{P}(x) \Leftrightarrow \forall A\ (A \in \mathcal{F} \Rightarrow x \in A)$. The predicate $\in$ exists by Theorem 4.5 and Example 4.6. Hence by Theorem 4.7.b there exists a set $P$ determined for ever $n$ by $F(n+1, P) = \{F(n, x) \mid \mathbf{P}(x)\}$. Let $x \in P$. Then by Theorem 4.7.c there exists for every $n$ such $y_n$, not necessarily different from $x$, that $F(n, x) = F(n, y_n)$ and $y_n \in A$ for every $A \in \mathcal{F}$. Hence $F(n, x) \in F(n+1, A)$ for every $A \in \mathcal{F}$ and every $n$. Therefore $x \in A$ for every $A \in \mathcal{F}$. Hence $\mathbf{P}(x)$ is true and by Theorem 4.7.e the set $\cap\mathcal{F}$ exists. □

Now we shall prove that the Axiom of Union is a theorem in AS. The *union* of a family of sets $\mathcal{F}$ is the set of all those sets that are elements of at least one $A \in \mathcal{F}$. The union of $\mathcal{F}$ is denoted by $\cup\mathcal{F}$. If $\mathcal{F} = \{A_1, A_2, ..., A_k\}$ we denote $\cup\mathcal{F}$ by $A_1 \cup A_2 \cup ... \cup A_k$.

**4.9 Theorem** *The union of an arbitrary family of sets exists.*

*Proof*. Let $\mathcal{F}$ be a family of sets. The predicate for $\cup\mathcal{F}$ is $\mathbf{P}(x) \Leftrightarrow \exists A\ (x \in A \wedge A \in \mathcal{F})$. This predicate exists by Theorem 4.5. Hence, by Theorem 4.7.b, there exists the set $P$ determined by $F(n+1, P) = \{F(n, x) \mid \mathbf{P}(x)\}$ for every $n$. If $x \in P$, then, by Theorem 4.7.c, there exists for every $n$ such $y_n \in A_n$ for some $A_n \in \mathcal{F}$ that $F(n, x) = F(n, y_n)$. This implies that for every $n$ we have $F(n, x) \in F(n+1, A_n)$ for some $A_n \in \mathcal{F}$. In other words, $F(n, x)$ is an element of at least one HF set represented by a node on the $(n+1)$th level of the developing tree of $\mathcal{F}$. Lemma 3.12 implies that $F(m, x)$ is an element of the HF sets represented by the ancestors of these nodes on the $(m+1)$th level for every $m < n$. Therefore, all the nodes of **U**, of any level $m + 1$, which represent a HF set whose element is $F(m, x)$, for $m \leq n$, form a tree structure. But this structure is not a tree by Definition 3.17, as every node in this structure does not necessarily have a successor inside this structure. However, since the nodes of this structure are on every level, the structure must contain at least one path $\mathbf{Q}_A$ in the sense of Definition 3.2. That is, at least one of the branches which starts at the root and continues to level $n$ can be prolonged to the level $n + 1$ inside this structure (König's Lemma). Since this structure is inside the developing tree of $\mathcal{F}$, the path $\mathbf{Q}_A$ represents some $A \in \mathcal{F}$. Thus $F(n, x) \in F(n + 1, A)$ holds for every $n$. Hence $x \in A$ and $\mathbf{P}(x)$ is true which, by Theorem 4.7.e, means that $\cup\mathcal{F}$ exists. □

Does the complement of a set exist? That is, does the set $\{x \mid x \notin A\}$ exist for arbitrary set $A$? In general this is not true. Let every element of some set $A$ be a set which is not a HF set. Then every HF set must be an element of $\{x \mid x \notin A\}$. But then $\{x \mid x \notin A\} = U$ (See the proof of Theorem 4.3.), and every element of $A$ is an element of $\{x \mid x \notin A\}$ which is a contradiction. Consequently, according to Theorem 4.7, the only sensible definition of complement would be:

**4.10 Definition** The *relative complement* of a set $A$ in a set $B$ is the smallest set containing all elements of $B$ that are not elements of $A$. This set is denoted by $B \setminus A$. The set $U \setminus A$ is the *complement* of $A$ denoted by $\mathsf{C}(A)$.

**4.11 Theorem** *The set $B \setminus A$ exists for arbitrary sets $A$ and $B$.*

*Proof.* If $A$ and B are two sets, then the predicates $x \in B$, $x \notin A$ exist by Theorem 4.5. Hence, by Theorem 4.7.d, there exists the set $B \setminus A$. □

It is not difficult to see that $A \cup \mathsf{C}(A) = U$, $\mathsf{C}(U) = \emptyset$ and $\mathsf{C}(\emptyset) = U$. However, we have already found out that $A \cap \mathsf{C}(A) = \emptyset$ is in general not true. In the following section we shall show when precisely this is true.

Now we shall prove that the Power set Axiom is valid in AS.

**4.12 Definition** A set $B$ is a *subset* of a set $A$ ($B \subseteq A$) if $x \in B$ implies $x \in A$. The *power set* of a set $A$ is the set whose elements are precisely all subsets of $A$. The power set of a set $A$ is denoted by $\mathcal{P}(A)$. A tree $\mathbf{T}_1$ is a *subtree of* a tree $\mathbf{T}_2$, if $\mathbf{T}_1(v, v')$ implies $\mathbf{T}_2(v, v')$ for every pair if nodes $(v, v')$ (See Definition 3.17).

**4.13 Lemma** $A_1 \subseteq A_2$ *if and only if* $F(n, A_1) \subseteq F(n, A_2)$ *for every n.*

*Proof.* Let $\mathbf{T}_1$ and $\mathbf{T}_2$ be the developing trees of $A_1$ and $A_2$ respectively. Then $A_1 \subseteq A_2$ precisely when $\mathbf{T}_1$ is a subtree of $\mathbf{T}_2$, since every path in $\mathbf{T}_1$ must be in $\mathbf{T}_2$ by Theorem 3.18. By Theorem 3.21 the nodes on the $n$th level of $\mathbf{T}_1$ and $\mathbf{T}_2$ represent the elements of $F(n + 1, A_1)$ and $F(n + 1, A_2)$ respectively. Thus $F(n + 1, A_1) \subseteq F(n + 1, A_2)$ must hold for every $n$. And since $F(0, A_1) \subseteq F(0, A_2)$, we have $F(n, A_1) \subseteq F(n, A_2)$ for every $n$. It is not difficult to see that this holds only if $A_1 \subseteq A_2$. In the opposite case some node on some level $n$ of $\mathbf{T}_1$ is not in $\mathbf{T}_2$ and consequently $F(n + 1, A_1)$ is not a subset of $F(n + 1, A_2)$. □

**4.14 Theorem** *The power set of an arbitrary set exists.*

*Proof.* The corresponding predicate for the power set of $A$ is $\mathbf{P}(x) \Leftrightarrow x \subseteq A$. It is not difficult to see that this predicate satisfies the conditions of Theorem 4.5 and that it therefore exists. Hence by Theorem 4.7.b there exists the set $P$ for which holds $F(n + 1, P) = \{F(n, x) \mid \mathbf{P}(x)\}$ for every $n$. Let $x \in P$. By Theorem 4.7.c there exists such $y_n \subseteq A$ for every $n$ that $F(n, x) = F(n, y_n)$. Therefore $F(n, x) \subseteq F(n, A)$ for every $n$ by Lemma 4.13. Hence $x \subseteq A$ by Lemma 4.13 and $\mathbf{P}(x)$ is a true statement which by Theorem 4.7.e means that $\mathcal{P}(A)$ exists. □

Till now we have come across only one function in AS, namely $F(n, A)$. If we want to investigate the status of the Continuum Hypothesis in AS and other important issues we must be able to decide which sets in AS are functions, bijections, etc. Since functions from a set $A$ to a set $B$ are special subsets of the set of all ordered pairs $(a, b)$, where $a \in A$ and $b \in B$, we must first define the ordered pair of two sets and then prove its existence. In other words, we must prove that the Axiom of Pair holds in AS.

**4.15 Definition** A *pair of sets a* and *b* is the set denoted by $\{a, b\}$ whose only elements are $a$ and $b$. The set whose only element is $a$ is denoted by $\{a\}$.The *ordered pair* of sets $a$ and $b$ is the set $\{\{a\},\{a, b\}\}$ denoted by $(a, b)$.

**4.16 Theorem** *If a and b are arbitrary sets, then the pair* $\{a, b\}$ *and the ordered pair* $(a, b)$ *exist. For* $(a, b)$ *holds* $F(n + 2, (a, b)) = (F(n, a), F(n, b))$ *for every n.*

*Proof.* The corresponding predicate for the unordered pair of sets $a$ and $b$ is $\mathbf{P}(x) \Leftrightarrow x = a \vee x = b$. This predicate exists by Theorem 4.5. Hence by Theorem 4.7.b there exists the set $P$ for which $F(n + 1, P) = \{F(n, x) \mid \mathbf{P}(x)\}$ holds for every $n$. If $x \in P$ then $F(n, x) = F(n, a)$ or $F(n, x) = F(n, b)$ for every $n$. Let $a \neq b$. Then there exists such $i$ that $F(i, a) \neq F(i, b)$. If $F(i, x) = F(i, a)$, then $F(j, x) = F(j, a)$ for every $j > i$. Namely, if for some $j > i$ we have $F(j, x) = F(j, b)$, then $F(i, x) = F(i, b) \neq F(i, a)$, which is not true. Since $F(i, x) = F(i, a)$, we have $F(j, x) = F(j, a)$ also for $j < i$. Hence $F(n, x) = F(n, a)$ for every $n$, which means $x = a$. In the same way we conclude that $F(i, x) = F(i, b)$ implies $x = b$. Hence $\mathbf{P}(x)$ is true and the pair $\{a, b\}$ exists. In a similar way we conclude that $\{a, a\}$ exists and hence $\{a\}$ exists. This implies that the ordered pair $(a, b)$ exists.

Let's prove now the second part of the theorem. We have

$$F(n + 2, (a, b)) = F(n + 2, \{\{a\}, \{a, b\}\}) = \{F(n + 1, \{a\}), F(n + 1, \{a, b\})\} =$$

$$\{\{F(n, a)\}, \{F(n, a), F(n, b)\}\} = (F(n, a), F(n, b))$$

by definition of the ordered pair of two HF sets in 3.1. In words, the $(n + 2)$th approximation of the ordered pair $(a, b)$ is the ordered pair of the $n$th approximations of sets $a$ and $b$. □

The *cartesian product of sets A* and *B* is the set of all ordered pairs $(a, b)$ where $a \in A$ and $b \in B$.

**4.17 Theorem** *The cartesian product* $A \times B$ *exists for arbitrary sets A and B. For every n holds* $F(n + 3, A \times B) = \{(F(n, a), F(n, b)) \mid a \in A, b \in B\}$.

*Proof.* The corresponding predicate is $\mathbf{P}(x) \Leftrightarrow \exists a\, \exists b\, (x = (a, b) \wedge a \in A \wedge b \in B)$. Since the ordered pair $(a, b)$ exists by the previous theorem, the predicate $\mathbf{P}$ is expressible by existent predicates. Hence, by 4.7.b, there exists the set $P$ determined by $F(n + 1, P) = \{F(n, x) \mid \mathbf{P}(x)\}$. Let $x \in P$. By Theorem 4.7.c we have $F(n, x) = F(n, (a_n, b_n))$ for every $n$ where $a_n \in A$ and $b_n \in B$. And by Theorem 3.5 we have

$$F(n, F(n + 1, (a_{n+1}, b_{n+1}))) = F(n, F(n + 1, x)) = F(n, x) = F(n, (a_n, b_n)) \quad (1)$$

By the previous theorem we have

$$F(n, F(n + 1, (a_{n+1}, b_{n+1}))) = (F(n - 2, F(n - 1, a_{n+1})), F(n - 2, F(n - 1, b_{n+1}))) \quad (2)$$

and $$F(n, (a_n, b_n)) = (F(n - 2, a_n), F(n - 2, b_n)) \quad (3)$$

Items (1) to (3) imply $F(n - 2, F(n - 1, a_{n+1})) = F(n - 2, a_n)$ and $F(n - 2, F(n - 1, b_{n+1})) = F(n - 2, b_n)$ for every $n \geq 2$. Hence $F(0, a_2)$, $F(1, a_3)$, $F(2, a_4)$ is the developing sequence of some set $a$. Since $a_{n+2} \in A$ for every $n$, we have $F(n, a) = F(n, a_{n+2}) \in F(n + 1, A)$ for every $n$. Hence $a \in A$. In a similar way we conclude that $F(0, b_2)$, $F(1, b_3)$, $F(2, b_4)$ is the developing sequence of some $b \in B$. Therefore the approximations $F(n, (a_n, b_n))$ form the developing sequence of $(a, b)$ which means that $x = (a, b)$ where $a \in A$ and $b \in B$. Hence $\mathbf{P}(x)$ is true by Theorem 4.7.e and $A \times B$ exists.

For the second part of the theorem we infer as follows. By the previous theorem we have

$$F(n + 3, A \times B) = F(n + 3, \{(a, b) \mid a \in A, b \in B\}) = \{F(n + 2, (a, b)) \mid a \in A, b \in B\} =$$

$$\{(F(n, a), F(n, b)) \mid a \in A, b \in B\}.$$

In words, the $(n + 3)$th approximation of the product $A \times B$ is the set of ordered pairs of all $n$th approximations of the elements of $A$ and $B$. □

A subset of $A \times B$ is a *relation between the elements of A and B*. The set of all $(y, x)$, where $(x, y) \in R$ for some relation $R$, is the *inverse* relation of $R$ denoted by $R^{-1}$.

**4.18 Theorem** *If a relation R exists, then the relation $R^{-1}$ exists.*

*Proof.* Let $R \subseteq A \times B$. The predicate for the inverse relation is $\mathbf{P}(x) \Leftrightarrow \exists a\ \exists b\ (x = (b, a) \wedge (a, b) \in R)$. Since $R$ exists, $\mathbf{P}$ is expressible by existent predicates and therefore $\mathbf{P}$ exists. Hence the set $P$, determined by $F(n + 1, P) = \{F(n, x) \mid \mathbf{P}(x)\}$, exists. If $x \in P$, then $F(n, x) = F(n, (b_n, a_n))$ for every $n$ where $(a_n, b_n) \in R$. In the same way as in the proof for the existence of cartesian product we find out that the developing sequence for $x$ is equal to the developing sequence of some ordered pair $(b, a)$ where $(a, b) \in R$. Therefore $x = (b, a)$. Hence $\mathbf{P}(x)$ is true and $R^{-1}$ exists on the base of Theorem 4.7. □

**4.19 Definition** A relation $R$ is
1) *total* if for every $a \in A$ there exists $b \in B$ such that $(a, b) \in R$.
2) *surjective* if for every $b \in B$ there exists $a \in A$ such that $(a, b) \in R$.
3) *functional* if $(a, b) \in R$ and $(a, b') \in R$ imply $b = b'$.
4) *injective* if $(a, b) \in R$ and $(a', b) \in R$ imply $a = a'$.

A *function f* which *maps* a set $A$ to a set $B$ is a functional total relation between the elements of $A$ and $B$. Such function is denoted by $f : A \to B$. Instead of $(a, b) \in f$ we write also $f(a) = b$. The set $A$ is the *domain* of $f$ and the set $B$ is the *range* of $f$. If $f$ is injective or surjective, then $f$ is an *injection* from $A$ into $B$ or a *surjection* from $A$ onto $B$ respectively. If $f$ is an injection and a surjection from $A$ to $B$, then $f$ is a *bijection* from $A$ to $B$. The set of all ordered pairs $(a, a)$, $a \in A$, is the *identity function on A*.

If $f$ is a bijection from $A$ to $B$, then $f^{-1}$ is a relation between the elements of $B$ and $A$. This relation is total because $f$ is surjective, it is surjective because $f$ is total, it is functional because $f$ is injective, and it is injective because $f$ is functional. Therefore $f^{-1}$ is a bijection as well. The identity function on $A$ is an example of a bijection. This functions exists for every set $A$. The proof is a modification of the proof of Theorem 4.17.

In ZF it is possible to prove that no surjection from $A$ onto $\mathcal{P}(A)$ exists for any set $A$. Namely, let $f$ be an arbitrary map from $A$ to $\mathcal{P}(A)$ and let $D \in \mathcal{P}(A)$ consist of all such $x \in A$ for which $x \notin f(x)$. Then we prove by diagonalization that $D$ is not an $f$-image of any $x \in A$. Hence $f$ is not a surjection from $A$ onto $\mathcal{P}(A)$ and therefore $\mathcal{P}(A)$ has *more* elements than $A$ for any set $A$. But the existence of the universal set $U$ in AS would thus imply the so called *Cantor's paradox* by which $U$ being the set of all sets would have *less* elements than its power set $\mathcal{P}(U)$. But, as $\mathcal{P}(U)$ is certainly a subset of $U$ it has no more elements than $U$. In short, $U$ would have less elements than $U$. However, this paradox does not appear in AS, since in AS the diagonalization cannot be always carried out.

**4.20 Theorem** *The argument for Cantor's paradox fails in AS.*

*Proof.* Let $A$ be an arbitrary set. As every element of $A$ is certainly an element of $U$, we have $A \subseteq U$. In other words, $A \in U$ implies $A \in \mathcal{P}(U)$. This means that $U \subseteq \mathcal{P}(U)$. And, since $\mathcal{P}(U) \subseteq U$ certainly holds, we have $U = \mathcal{P}(U)$. Therefore, the identity map on $U$ is a bijection between $U$ and $\mathcal{P}(U)$, and so $U$ and $\mathcal{P}(U)$ have the *same number* of elements. □

Therefore, if $f$ is the identity function on $U$, then $f$ is a surjection from $U$ onto $\mathcal{P}(U)$. For this function we cannot find the set of all those $x \in U$ for which $x \notin f(x)$, since such a set would be equal to the set of precisely all those $x \in U$ that are not elements of themselves, and this set does not exist (See the proof of Theorem 4.3.). Thus the diagonalization is not always realizable in AS.

Now we shall prove that the Axiom of Infinity holds in AS. That is, we shall prove that the set of all natural numbers exists in AS and that it is infinite. We shall see that this set contains an additional element which is not a natural number. Let's describe this element. The first $n$ terms of the developing sequence of a natural number $n$ are equal to the first $n$ natural numbers, and all the other terms are equal to $n$ (See Example 1.4.1). For instance, the developing sequence of the natural number 3 is 0, 1, 2, 3, 3, 3, 3, ... This fact can be naturally extrapolated. We can say that the sequence of all natural numbers 0, 1, 2, 3, 4, 5, 6, ..., is the developing sequence of the so called infinite number.

**4.21 Definition** The set whose $n$th approximation is equal to $n$ is the *infinite number*. This set is denoted by $\infty$. For every $n$ we have $n < \infty$, $\infty > n$ and $\infty + n = n + \infty = \infty$.

By Example 1.4.1 we have $F(n, F(n + 1, \infty)) = F(n, n + 1) = n = F(n, \infty)$ for every $n$. So the set $\infty$ exists by Theorem 3.5. Since $F(n, \infty)$ is the only element of $F(n + 1, \infty)$ for every $n$, we have $\infty = \{\infty\}$ and the set $\infty$ is the only element of itself. The existence of this set proves again that in AS the Axiom of Foundation does not hold. Notice that we have not defined infinite number as a natural number. Every natural number is a HF set while infinite number is not a HF set.

**4.22 Theorem** *There exists the set of all natural numbers whose element is also the set* $\infty$.

*Proof.* The predicate $\mathbf{N}(x)$ for "to be a natural number" is computable on HF sets, as we have seen in the 2. section. Hence the set $\boldsymbol{N}$, determined by $F(n + 1, \boldsymbol{N}) = \{F(n, x) \mid \mathbf{N}(x)\}$, exists. Let $x \in \boldsymbol{N}$. Then for every $n$ there exist such natural number $y_n$ that $F(n, x) = F(n, y_n)$ which, by Example 1.4.1, means $F(n, x) = \min\{n, y_n\}$. Hence $F(n, x)$ is a natural number which is not larger than $n$. If $F(n, x) < n$ for some $n$, then $F(n, x) = x$ by Theorem 3.15, and $x$ is a natural number not larger than $n$. If, however, $F(n, x) = n$ for every $n$, then $x = \infty$ by Definition 4.21. □

**4.23 Definition**. If $B \subseteq A$ and $B \neq A$, then $B$ is a *proper* subset of $A$. A set is *infinite* if there exists a bijection between this set and its proper subset.

**4.24 Theorem** *The set* $\boldsymbol{N}$ *is infinite.*

Let $\mathbf{P}$(x) be the predicate "to be a natural number greater than 0". This predicate is computable on HF sets. Hence the set $\boldsymbol{N}'$, determined for every $n$ by $F(n + 1, \boldsymbol{N}') = \{F(n, x) \mid \mathbf{P}(x)\}$, exists. The elements of $\boldsymbol{N}'$ are all natural numbers greater than 0 and the set $\infty$. The proof of this is almost exactly the same as the proof for the existence of $\boldsymbol{N}$. The function which maps $i$ to $i + 1$ for every $i$ is a bijection between $\boldsymbol{N}$ and $\boldsymbol{N}'$. This function maps $\infty$ to $\infty$. Since $\boldsymbol{N}'$ is a proper subset of $\boldsymbol{N}$, the set $\boldsymbol{N}$ is infinite.□

# 5. Topology of Approximations and the Continuum Hypothesis

In this section we shall introduce *topology of approximations* into the universal set $U$. This topology, which we shall denote by $\mathcal{T}_a$, corresponds to the interval topology in ZF. We shall prove that every bijection in AS is a homeomorphism in $\mathcal{T}_a$. This fact has very important consequences for the truthfulness of the Continuum Hypothesis in AS. Namely, we shall prove that $\mathcal{P}(\boldsymbol{N})$ contains an infinite subset which is perfect in $\mathcal{T}_a$. This means that every element of this set is an accumulation point of this set in $\mathcal{T}_a$. On the other hand $\boldsymbol{N}$ has only one accumulation point in $\mathcal{T}_a$, namely $\infty$. We shall also prove that for every $n$ there exists an infinite set with precisely $n$ accumulation points in $\mathcal{T}_a$. It is not difficult to conclude that all these facts disprove the Continuum Hypothesis.

First we shall intuitively show that every set in AS can be interpreted as a set of points in the interval [0, 1]. Namely, we assign the interval [0, 1] to the root of **U**. Then we proceed recursively. Let's say that we have assigned a closed interval $[a, b]$ to some node $v$ and that the successors of this node are $v_1, v_2, ..., v_k$. We choose in $[a, b]$ non-degenerated mutually disjoint closed intervals $I_1, I_2, ..., I_k$ none of which is longer than the half of $[a, b]$. We assign each $I_i$ to the node $v_i$ for $i = 1, 2, ..., k$ and we repeat the procedure for every $I_i$. If $x$ is an arbitrary element of $U$, then a chain of closed intervals $J_0(x) \supset J_1(x) \supset J_2(x), \ldots$ is assigned by this recursive procedure to the nodes which represent the approximations $F(0, x), F(1, x), F(2, x)$ .... Since the length of $J_{n+1}(x)$ is at the most half of the length of $J_n(x)$, the intersection of all $J_n(x)$ contains precisely one point which we assign to $x$.

The described procedure assigns different points of [0, 1] to different elements of the universal set $U$. Thus $U$ can be interpreted as a set of points in [0, 1]. Since every set is a subset of $U$, every set can be interpreted as a subset of [0, 1]. This shows that there is no set in AS which has *more* than "continuum many" elements. Of course all these conclusions are only intuitive, since we have not yet defined the interval [0, 1], and even less have we proved that such an interval exists in AS. In spite of that we can, on the base of this observation, introduce into $U$ a topology which corresponds to the interval topology in [0, 1]. In order to introduce this topology we must prove the existence of a particular type of sets.

**5.1 Theorem** *The set* $\{x \mid F(i, x) = F(i, a)\}$ *exists for every i and every* $a \in U$.

*Proof.* The corresponding predicate is $\mathbf{P}(i, a, x) \Leftrightarrow F(i, x) = F(i, a)$. This predicate exists because it is expressible by existent predicates in predicate calculus. Therefore there exists a set $P$ determined by $F(n + 1, P) = \{F(n, x) \mid F(i, x) = F(i, a)\}$ for every $n$. Let $x \in P$. Then for every $n$ there exists such $y_n$ that $F(n, x) = F(n, y_n)$ and $\mathbf{P}(i, a, y_n)$ is true. This implies $F(i, x) = F(i, y_i) = F(i, a)$. Therefore $\mathbf{P}(i, a, x)$ is true and consequently the set $\{x \mid \mathbf{P}(x)\}$ exists by Theorem 4.7. □

**5.2 Definition** The elements of $U$ are *points* in $U$. The set of points $\{ y \mid F(n, y) = F(n, x)\}$ is a *neighborhood* of the point $x$ for every $n$. This set is denoted by $V(n, x)$. The topology determined in this way is the *topology of approximations* $\mathcal{T}_a$. The universal set $U$ equipped with $\mathcal{T}_a$ is the *topological space* $(U, \mathcal{T}_a)$. A point $x$ is *in* a set $A$

*a*) *internal* if there is a neighborhood of $x$ with all of its points from $A$
*b*) *external* if there is a neighborhood of $x$ with none of its points from $A$
*c*) *boundary* if in every neighborhood of $x$ there are points from $A$ and points which are not from $A$.

A set $A$ is *open* if every one of its points is internal in $A$. A point $x$ is an *accumulation point in* a set $A$ if in every neighborhood of $x$ there is at least one point $y \in A, y \neq x$. A point $x \in A$ which is not an accumulation point in $A$ is an *isolated point in A*. A set is *closed* if every of its accumulation points is also its element. The *closure* of a set $A$ is the smallest closed set whose subset is $A$ (here $X$ is *smaller* than $Y$ means $X \subset Y$ ).

Notice that the expression $(U, \mathcal{T}_a)$ is not an ordered pair of two sets but only a designation because $\mathcal{T}_a$ is not a set. Also notice that $n$ in $V(n, x)$ is a natural number, that is, $n \neq \infty$. At the end of the article we shall prove that $F(\infty, x) = x$ for every $x$ (Theorem 5.15). If we allow $V(\infty, x)$ to be a neighborhood of $x$ then the resulting topology is the trivial discrete topology where every point is isolated. This topology does not correspond to the interval topology in ZF and is without any use for our investigations.

Whenever we shall deal with topological structure of sets, without mentioning the topology explicitly, we shall have in mind $\mathcal{T}_a$. Now we shall prove a very important theorem of AS.

**5.3 Theorem** *Every set is closed in* $(U, \mathcal{T}_a)$.

*Proof.* Let $x$ be an accumulation point in a set $A$. Then for every $n$ there exists such $y \in A$, $y \neq x$, that $F(n, y) = F(n, x)$. This means that $F(n, x) \in F(n + 1, A)$ for every $n$. Therefore $x \in A$ by Definition 3.8. Thus $A$ includes all of its accumulation points and for this reason $A$ is closed in $(U, \mathcal{T}_a)$. □

When we were investigating Russell's and Cantor's paradox we saw that diagonalization is in general not realizable in AS. Theorem 5.3 explains why this is so. The reason is the fact that for every set $A$ and every function $f : A \rightarrow \mathcal{P}(A)$ the set $D = \{x \mid x \in A \text{ and } x \notin f(x)\}$ does not always exist. By Theorem 5.3 there exists at most the closure of $D$ in $(U, \mathcal{T}_a)$, say $D'$, whose elements can, however, be also such points $x \in A$ for which $x \in f(x)$ holds and which are accumulation points in $D'$. These points cannot be eliminated from $D'$. And for this reason we cannot make the conclusion that there is no such $x \in A$ which does not map to $D'$ by $f$. Of course, we can still make this conclusion when $A$ is finite, as in this case any subset of $A$ has no accumulation points.

As we have already found out, the equation $A \cap \mathsf{C}(A) = \emptyset$ does not hold in general. Now we see what is the reason. The sets $A$ and $\mathsf{C}(A)$ contain all their accumulation points and some points can be accumulation points in $A$ and in $\mathsf{C}(A)$.

**5.4 Theorem** *Every set has precisely all its boundary points in common with its complement.*

*Proof.* Every $x \in U$ is either a boundary, an internal or an external point in a set $A$. If $x$ is a boundary point in $A$, then in every neighborhood of $x$ there exist points from $A$ and points from $\mathsf{C}(A)$. Hence $x$ is an accumulation point in $A$ and in $\mathsf{C}(A)$. Since $A$ and $\mathsf{C}(A)$ contain all their accumulation points by Theorem 5.3, we have $x \in A \cap \mathsf{C}(A)$. If $x$ is an internal point in $A$, then there exists at least one neighborhood of $x$ in which no point from $\mathsf{C}(A)$ exists, hence $x \notin \mathsf{C}(A)$. However, if $x$ is an external point in $A$, then $x \notin A$. □

**5.5 Theorem** $A \cap \mathsf{C}(A) = \emptyset$ *precisely when A is an open set. The set* $\mathsf{C}(A)$ *is open precisely when the set A is open.*

*Proof.* If $A$ is open, then $A$ has no boundary points. Hence $A \cap \mathsf{C}(A) = \emptyset$ by the previous theorem. If, however, $A \cap \mathsf{C}(A) = \emptyset$, then $A$ has no boundary points by the same theorem. Hence every point in $A$ is internal. Because of the symmetric role of $A$ and $\mathsf{C}(A)$ in this proof the set $\mathsf{C}(A)$ is open precisely when the set $A$ is open. □

A function $f : A \rightarrow B$ is continuous in a point $a \in A$ if for every neighborhood $W$ of $f(a)$ in $B$ there exists such small enough neighborhood $V$ of $a$ in $A$ that $f(V) \subseteq W$. The function $f$ is continuous if it is continuous in every point of $A$. In $\mathcal{T}_a$ the neighborhood $V(j, x)$ is *smaller* than the neighborhood $V(i, x)$ if $j > i$, since every point $y$ which is in the neighborhood $V(j, x)$ is in the neighborhood $V(i, x)$ but not vice versa. Therefore the greater the number $i$ is the smaller is the neighborhood $V(i, x)$. Hence

**5.6 Definition** In $\mathcal{T}_a$ a function $f : A \to B$ is *continuous* in a point $a \in A$ if for every $i$ there exists such large enough $j$ that $F(j, x) = F(j, a)$ implies $F(i, f(x)) = F(i, f(a))$ for every $x \in A$. A function is *continuous* if it is continuous in every point of $A$.

Now we shall prove the second very important theorem of AS which is that every function is continuous. In the proof and in some subsequent proofs we shall not distinguish between nodes and the approximations they represent. We shall do this for the sake of greater simplicity and therefore greater clarity, since the approximations will be given in the form which will precisely determine the nodes that represent them. Thus we shall use the expressions "a successor of the approximation" or "the node $F(n, x)$", etc.

Let $f$ be a function which maps $A$ to $B$. The developing tree of $f$ is a subtree of the developing tree of $A \times B$, since every path from the first one must be in the second one by Theorem 3.18. Therefore, by Theorem 4.17, the nodes of the developing tree of $f$ are of the form $(F(n, a), F(n, b))$ where $a \in A$ and $b \in B$. These nodes are on the $(n + 2)$th level of the developing tree of $f$.

**5.7 Lemma** *Let $f : A \to B$ be a function and $a \in A$. If for every $i$ there exists such $j$ that the node $(F(i, a), F(i, f(a)))$ is an ancestor of every node $(F(j, a), F(j, b))$ of the developing tree of $f$, where $b$ can be any element of $B$, then $f$ is continuous in the point $a$ in $(U, \mathcal{T}_a)$.*

*Proof.* Let $a \in A$, let $i$ be some natural number and $j$ the corresponding natural number for which the condition of the theorem is fulfilled. Let $(x, f(x)) \in f$ and $F(j, x) = F(j, a)$. Then the $(j + 2)$th approximation of $(x, f(x))$ is equal to $(F(j, x), F(j, f(x)))$, and this is equal to $(F(j, a), F(j, b))$ where $b = f(x)$. Therefore, by assumption of the theorem, the node $(F(i, a), F(i, f(a)))$ is an ancestor of the node $(F(j, x), F(j, f(x)))$. This ancestor is in fact the approximation $(F(i, x), F(i, f(x)))$, which implies $F(i, f(x)) = F(i, f(a))$. In short, if $F(j, x) = F(j, a)$, then $F(i, f(x)) = F(i, f(a))$, which means that $f$ is continuous in the point $a$. □

In the proof of the following theorem every node is a node of the developing tree of a function $f$ and $H_i$ denotes a HF set for which there exists such $b \in B$ that $F(i, b) = H_i$ holds.

**5.8 Theorem** *Every function is continuous in $(U, \mathcal{T}_a)$.*

*Proof.* Let $f : A \to B$ be a function for which the condition from the previous lemma does not hold in some point $a$. Then there exists such $i$ that for every $j$ there exists a node $(F(j, a), H_j)$ which is not a descendant of the node $(F(i, a), F(i, f(a)))$. Then there must exist a node $(F(i, a), H_i)$, where $H_i \neq F(i, f(a))$, for which for every $j$ there exists such $k > j$ that a node $(F(k, a), H_k)$ is a descendant of the node $(F(i, a), H_i)$. In the opposite case for every node $(F(i, a), H_i')$, where $H_i' \neq f(a)$, there would exist such $j$ that $(F(i, a), H_i')$ would have no descendant of the form $(F(k, a), H_k)$ for any $k > j$. Then none of the nodes $(F(i, a), H_i')$, where $H_i' \neq f(a)$, would have a descendant of the form $(F(k, a), H_k)$ for some great enough $k$ and therefore every node of such form would be a descendant of the node $(F(i, a), F(i, f(a)))$ which contradicts the assumption. And then we conclude in the same way that there must exist such a successor of $(F(i, a), H_i)$, say $(F(i + 1, a), H_{i+1})$, that for every $j$ there exists such $k > j$ that a node of the form $(F(k, a), H_k)$ is a descendant of the node $(F(i + 1, a), H_{i+1})$. And in the same way we conclude that there exists a successor of $(F(i + 1, a), H_{i+1})$, say $(F(i + 2, a), H_{i+2})$, which has the same properties ..., etc. The sequence $(F(i, a), H_i)$, $(F(i + 1, a), H_{i+1})$, $(F(i + 2, a), H_{i+2})$, ... together with the ancestors of the node $(F(i, a), H_i)$ form the developing sequence of some $(a, b) \in f$, where $H_n = F(n, b)$ for every $n$. Since $F(i, b) = H_i \neq F(i, f(a))$, the point $a$ would be in the relation $f$ with $f(a)$ and with $b$ where $b \neq f(a)$. Thus $f$ would not be a function. This contradiction shows that $f$ is continuous in every point $a$. □

**5.9 Remark** It is known that the following theorem holds in ZF: *A function $f$ is continuous precisely when $f^{-1}(A)$ is a closed set for any closed set $A$.* All sets in AS are closed in $(U, \mathcal{T}_a)$ by Theorem 5.3.

Thus Theorem 5.8 should follow already from Theorem 5.3. However, such a conclusion is not valid in AS. In the proof of the stated theorem from ZF we use the fact that the complement of a closed set is an open set. But this fact is derived from the Axiom of Separation and does not hold in AS. Namely, there are closed sets in AS whose complements are not open sets. It is evident from Theorem 5.4 that these sets have nonempty boundary.

A bijection $f$ is a *homeomorphism* if $f$ and $f^{-1}$ are both continuous.

**5.10 Theorem** *Every bijection is a homeomorphism in* $(U, \mathcal{T}_a)$.

*Proof.* If $f$ is a bijection, then $f^{-1}$ is a bijection as well. The functions $f$ and $f^{-1}$ are both continuous in $(U, \mathcal{T}_a)$ by Theorem 5.8 and consequently $f$ is a homeomorphism in $(U, \mathcal{T}_a)$. □

Thus, if there exists a bijection between two sets, these two sets are homeomorphic in topology $\mathcal{T}_a$. This surprising fact gives to bijections in AS a different meaning from the one they have in ZF. The existence of a bijection between two sets in AS does not only show the same size of these two sets but also the same topological structure of them.

Now we shall investigate the topological structure of $\mathcal{P}(\boldsymbol{N})$. This set exists by Theorems 4.22 and Theorem 4.14. First we shall prove a lemma which is very well known in ZF. We must prove it again because in the ZF proof the Axiom of Separation might be used, and this axiom, as we know, is not valid in AS.

**5.11 Lemma** *A homeomorphism in* $(U, \mathcal{T}_a)$ *maps isolated points to isolated points and accumulation points to accumulation points.*

*Proof.* Let $f: A \to B$ be a homeomorphism in $(U, \mathcal{T}_a)$, $a \in A$ and $f(a)$ an isolated point in $B$. Then there exists such $n$ that $f(a)$ is the only point from $B$ in the neighborhood $V(n, f(a))$. Because of the continuity of $f$ there exists such $i$ that for every $x \in A$, which is in the neighborhood $V(i, a)$, holds that $f(x)$ is in the neighborhood $V(n, f(a))$. Therefore $f$ maps every such $x$ to $f(a)$, as $f(a)$ is the only point from $B$ which is in the neighborhood $V(n, f(a))$. As $f$ is bijective, there can be only one such $x$, which implies $x = a$. Thus $a$ is the only point from $A$ which is in the neighborhood $V(i, a)$. Hence $a$ is an isolated point in $A$. Since the function $f^{-1}$ is continuous and bijective as well, the point $a$ is isolated in $A$ precisely when the point $f(a)$ is isolated in $B$. □

**5.12 Theorem** *If $x$ is a* HF *set, then $x$ is an isolated point in every set. In other words, an accumulation point in some set is not a* HF *set.*

*Proof.* Let $x$ be a HF set and let $y$ be a point in the neighborhood $V(r(x) + 1, x)$ of $x$. Then we have $F(r(x) + 1, x) = F(r(x) + 1, y)$. As $F(r(x) + 1, x) = x$, by Theorem 2.18, we have $r(F(r(x) + 1, y)) = r(x)$. This and Theorem 3.15 imply $F(r(x) + 1, y) = y$. Therefore $y = x$. Hence, apart from the point $x$, there are no other points in the neighborhood $V(r(x) + 1, x)$. □

A *perfect set* is a nonempty closed set whose every element is an accumulation point in this set. Thus a perfect set is always infinite.

**5.13 Theorem** *There exists a subset of* $\mathcal{P}(\boldsymbol{N})$ *which is perfect in* $\mathcal{T}_a$.

*Proof.* By Theorem 3.13 we have

$$F(n + 2, \mathcal{P}(\boldsymbol{N})) = F(n + 2, \{A \mid A \subseteq \boldsymbol{N}\}) = \{F(n + 1, A) \mid A \subseteq \boldsymbol{N}\} = \{\{F(n, i) \mid i \in A\} \mid A \subseteq \boldsymbol{N}\}$$

Since $F(n, i) = \min\{n, i\}$, the elements of $F(n + 2, \mathcal{P}(\boldsymbol{N}))$ are all possible sets of natural numbers not greater than $n$. These sets are represented by the nodes on the $(n + 1)$th level of the developing tree of $\mathcal{P}(\boldsymbol{N})$ (Theorem 3.21). That is, if $v = (n + 1, S)$ is such node, then $S = \{F(n, i) \mid i \in A\}$ for some $A \subseteq \boldsymbol{N}$. From these nodes we select those for which $n \in S$. If we perform such a selection on every level $n > 0$, then the selected nodes form a complete binary tree. Namely, every successor of the previously described node $v$ is of the form $(n + 2, S')$ where $S' = \{F(n + 1, i) \mid i \in A\}$. If $k \in S$ and $k < n$, then $k \in A$, since such $k$ is the only natural number whose $n$th approximation is equal to $k$. And then $k \in S'$. Since a successor of $v$ must be one of the selected nodes, we have $n + 1 \in S'$. This means that $S'$ is almost completely determined. The only open question is whether $n \in S'$ or not. Hence $v$ has two successors. As $v$ was selected arbitrarily, any selected node has two successors. Therefore the selected nodes form a complete binary tree which starts at level 1. If we add the root of **U** to this tree, then we obtain the developing tree of some set $\mathcal{P}'(\boldsymbol{N}) \subseteq \mathcal{P}(\boldsymbol{N})$ (Corollary 3.20). Every $x \in \mathcal{P}'(\boldsymbol{N})$ is an accumulation point in $\mathcal{P}'(\boldsymbol{N})$. Namely, for every $x \in \mathcal{P}'(\boldsymbol{N})$ and every $n$ there exists such $y \in \mathcal{P}'(\boldsymbol{N})$, $y \neq x$, that $F(n + 1, x) = F(n + 1, y)$, and that one of $F(n + 2, x)$ and $F(n + 2, y)$ contains $n$ and the other does not, but both contain $n + 1$. Since $\mathcal{P}'(\boldsymbol{N})$ is closed by Theorem 5.3, it is a perfect set. □

Let $A <_c B$ denote that there exists an injection from $A$ into $B$ but not from $B$ into $A$, symbolizing that $A$ is of smaller *cardinality* than $B$. The Continuum Hypothesis asserts that there does not exist a set $S$ for which $\boldsymbol{N} <_c S <_c \mathcal{P}(\boldsymbol{N})$ would hold. The Continuum Hypothesis is not true in AS. Even more

**5.14 Theorem** *There exists a chain* $A_1 <_c A_2 <_c A_3 <_c \ldots <_c \mathcal{P}(\boldsymbol{N})$ *where* $A_1$ *is homeomorphic to* $\boldsymbol{N}$.

*Proof.* First we shall prove that the set of all unordered pairs $\{j, k\}$, $j, k \in \boldsymbol{N}$, $j < i$, $k \geq i$ exists for every $i > 0$. The corresponding predicate is $\mathbf{P}_i(x) \Leftrightarrow \exists j\, \exists k\, (j \in \boldsymbol{N} \wedge k \in \boldsymbol{N} \wedge j < i \wedge k \geq i \wedge x = \{j, k\})$. As $\mathbf{P}_i(x)$ is a conjunction of the predicates which exist, there exists a set $A_i$, by Theorem 4.7.b, determined by $F(n + 1, A_i) = \{F(n, x) \mid \mathbf{P}_i(x)\}$. Let $x \in A_i$. Then for every $n$ there exists a pair $\{j_n, k_n\}$ for which $\mathbf{P}_i(\{j_n, k_n\})$ is true and $F(n, x) = F(n, \{j_n, k_n\})$. By Theorem 3.5 we have

$$F(n, F(n + 1, \{j_{n+1}, k_{n+1}\})) = F(n, F(n + 1, x)) = F(n, x) = F(n, \{j_n, k_n\})$$

which implies $\{F(n-1, F(n, j_{n+1})), F(n-1, F(n, k_{n+1}))\} = \{F(n, j_n), F(n, k_n)\}$ (1)

Since $j_n < k_n$ for every $n$, we have $F(m, j_n) \leq F(m, k_n)$ for every $m$ and every $n$. This and (1) implies that $F(n - 1, F(n, j_{n+1})) = F(n, j_n)$ and $F(n - 1, F(n, k_{n+1})) = F(n, k_n)$ for every $n$. Thus, the approximations $F(n, j_n)$ and $F(n, k_n)$ form two developing sequences. It is not difficult to see that the first one is the developing sequence of some natural number $j$ smaller than $i$, and that the second one is the developing sequence of some natural number $k$ not smaller than $i$. Hence $F(n, x) = F(n, \{j, k\})$ for every $n$ where $\mathbf{P}_i(\{j, k\})$ is true. Consequently, $x = \{j, k\}$ and $\mathbf{P}_i(x)$ is true. Therefore $A_i = \{x \mid \mathbf{P}_i(x)\}$ for every $i > 0$.

The set $A_i$ is infinite for every $i > 0$. Namely, let $A_i'$ be the set of pairs $\{j, k\} \in A_i$ for which $k > i$. This set is a proper subset of $A_i$, since it does not contain the pairs $\{j, i\} \in A_i$. The existence of $A_i'$ can be proved in the same way as the existence of $A_i$. The function which maps $\{j, k\}$ to $\{j, k + 1\}$ for every $\{j, k\} \in A_i$ is a bijection between $A_i$ and $A_i'$. Therefore $A_i$ is infinite.

The pair $\{j, \infty\}$ is an accumulation point in $A_i$ for every $j < i$, since for every $\{j, k\} \in A_i$, where $k \geq n$, we have $F(n, \{j, k\}) = F(n, \{j, \infty\})$, but $\{j, k\}) \neq \{j, \infty\}$ for $k < \infty$. By Theorem 5.12 there are no other accumulation points in $A_i$. Hence $A_i$ has $i$ accumulation points, as $j \neq j'$ implies $(j, \infty) \neq (j', \infty)$.

The set $\infty$ is an accumulation point in $\boldsymbol{N}$, as $F(n, \infty) = F(n, n)$ for every $n$. The set $\infty$ is, by Theorems 5.12 and 4.22, the only accumulation point in $\boldsymbol{N}$. A function which maps $\{0, k\}$ to $k - 1$ for every $\{0, k\} \in A_1$ is a bijection from $A_1$ to $\boldsymbol{N}$. A function which maps $\{j, k\}$ to $\{j, k + 1\}$ for every $\{j, k\} \in A_i$ is an injection from $A_i$ to $A_{i+1}$. However, since $A_{i+1}$ has more accumulation points than $A_i$, no bijection between these two sets exists by Lemma 5.11 and Theorem 5.10. Thus $A_i <_c A_{i+1}$ holds for any $i$. Since $A_i \subseteq \mathcal{P}(\boldsymbol{N})$, the identity map on $A_i$ is an injection from $A_i$ into $\mathcal{P}(\boldsymbol{N})$. But an injection from $\mathcal{P}(\boldsymbol{N})$ into $A_i$ is not possible for any $i$. Namely, by Theorem 5.13, there exists a perfect subset of $\mathcal{P}(\boldsymbol{N})$ having infinitely many accumulation points, while there are only $i$ accumulation points in $A_i$. □

We have one more debt to settle. Till now we have treated $F$ as a function in predicate calculus. Now we shall show that $F$ is also a function by Definition 4.19, if we properly extend its domain and range. We shall also show that $F(\infty, A) = A$ for every set $A$.

**5.15 Theorem** *$F(n, A)$ is a function $F : \boldsymbol{N} \times U \rightarrow U$ where $F(\infty, A) = A$ for every set $A$.*

*Proof.* The function $F$ is determined for every element of $\boldsymbol{N} \times U$ except for the pairs $(\infty, A)$ where A is an arbitrary set. Since $F$ must be continuous by Theorem 5.8, there exists such $i$ for every $n$ that $F(n, F(x)) = F(n, F(\infty, A))$, if only $F(i, x) = F(i, (\infty, A))$. If we put $x = (j, A)$, where $j \geq i - 2$, then we obtain, by Theorem 4.16,

$$F(i, x) = F(i, (j, A)) = (F(i-2, j), F(i-2, A)) = (F(i-2, \infty), F(i-2, A)) = F(i, (\infty, A)).$$

Hence $F(n, F(j, A)) = F(n, F(\infty, A))$ must hold for $j \geq i - 2$. This certainly holds if $j \geq \max\{i-2, n\}$. Therefore we have $F(n, F(\infty, A)) = F(n, F(j, A)) = F(n, A)$ for every $n$. And then Theorem 3.4 implies $F(\infty, A) = A$ and $F$ is a surjection from $\boldsymbol{N} \times U$ onto $U$. □

# 6. Cardinality of Sets in AS

If there exists an injection from a set $A$ into a set $B$ we say that $A$ is smaller or at most of the same size as $B$. In ZFC (Zermelo-Fraenkel system of axioms with the addition of the Axiom of Choice) there always exists an injection from one of two arbitrary sets into the other. This follows from the fact that in ZFC every set can be well-ordered. Of course, we must realize that in AS an injection from a set $A$ into a set $B$ is a homeomorphism between $A$ and its image in $B$ in topology of approximations $\mathcal{T}_a$ (Theorem 5.10). Consequently, the existence of an injection from $A$ into $B$ means that there is a topological embedding of $A$ into $B$ in $\mathcal{T}_a$. If then an injection from $B$ into $A$ does not exist we could hardly say that this implies that $A$ is of smaller size than $B$. We can topologically embed a line into a circle, but not vice versa. However, a line surely has not less points than a circle. In spite of this we shall investigate the cardinality of sets based on injections. We shall see that in AS the cardinality is much less undetermined than in ZFC.

If there exists an injection from $A$ into $B$ and an injection from $B$ into $A$ then in ZF, by Cantor-Bernstein Theorem, there exists a bijection between $A$ and $B$ and we say that $A$ and $B$ are of equal size. As we shall see Cantor-Bernstein Theorem does not hold in AS in general since we cannot transfer the proof of this theorem from ZF to AS. Namely, this proof is based on the families of sets which exist only by the Axiom of Separation. And by this axiom the accumulation points of sets in $\mathcal{T}_a$ are not necessarily elements of these sets which must be the case in AS by Theorem 5.3. Before we disprove Cantor-Bernstein Theorem we must define some new concepts and prove one lemma.

**6.1 Definition**
a) A set $A$ is *in bijection with* a set $B$ if there exists a bijection between these two sets. A set is *infinite* if it is in bijection with its proper subset. A set is *finite* if it is not infinite.

b) The developing tree of a set $A$ will be denoted by $\mathbf{T}_A$ (See Theorem 3.18 and Definition 3.19). Let $v$ be a node of $\mathbf{T}_A$ and let $B$ be a subset of $A$ such that $\mathbf{T}_B$ consists of all infinite paths from $\mathbf{T}_A$ which go through $v$. Then *$B$ develops in $A$ from $v$*. The subsets $B_1$ and $B_2$ of $A$, which develop in $A$ from $v_1$ and $v_2$ respectively, are *separated in $A$* if neither $v_1$ is a descendant of $v_2$ nor $v_2$ is a descendant of $v_1$ (See Definition 3.1).

It is not difficult to see that separated sets have no common elements.

**6.2 Lemma** *Every set is in bijection with some subset of a perfect set.*

*Proof.* Let $A$ be an arbitrary set and let $P$ be a set which is perfect in $\mathcal{T}_a$. We shall construct a subtree $\mathbf{T}_{P'}$ of $\mathbf{T}_P$ such that different infinite paths from $\mathbf{T}_A$ will correspond to different infinite paths from $\mathbf{T}_{P'}$ and vice versa. In this way we shall prove that there is a bijection between $A$ and $P' \subseteq P$ whose developing tree is $\mathbf{T}_{P'}$ (See Theorem 3.18). Let $a_1, a_2, ..., a_k$ be all the successors of the root of $\mathbf{T}_A$. Since $P$ is an infinite set, there certainly exists such $n_1$ that $\mathbf{T}_P$ has at least $k$ nodes on the $n_1$th level. Let's choose from these nodes the nodes $p_1, p_2, ..., p_k$ and assign $a_1, a_2, ..., a_k$ to $p_1, p_2, ..., p_k$ respectively. From $a_1, a_2, ..., a_k$ develop in $A$ mutually disjoint subsets of $A$ (See Definition 6.1.b), and from $p_1, p_2, ..., p_k$ develop in $P$ mutually disjoint perfect subsets of $P$ otherwise $P$ would contain isolated points, which is not the case since $P$ is perfect. By the same argument as before we come to the conclusion that each node $p_i$ has on some level $n_2$ of $\mathbf{T}_P$ at least as many descendants as $a_i$ has successors in $\mathbf{T}_A$. Thus we can assign every successor of $a_i$ uniquely to the corresponding descendant of $p_i$ on level $n_2$. By using this argument recursively we construct the subtree $\mathbf{T}_{P'}$ of $\mathbf{T}_P$ whose infinite paths correspond to the infinite paths of $\mathbf{T}_A$. Hence $P' \subseteq P$ is in bijection with $A$. □

**6.3 Theorem** *Cantor-Bernstein Theorem does not hold in general.*

*Proof.* Let $P$ be a perfect set in $\mathcal{T}_a$ and let $A = P \cup \{a\}$, where $a$ is a hereditarily finite set. By Theorem 5.12 the point $a$ is isolated in every set. This implies $a \notin P$. By the previous lemma there exists an injection from $A$ into $P$. Of course, there also exists an injection from $P$ into $A$ as $P$ is a subset of $A$. However, since $a$ is an isolated point in every set, no bijection between $P$ and $A$ exists by Lemma 5.11 and Theorem 5.10. □

This theorem and the observations we made in the beginning of this section force us to the following definition of the comparison of cardinality.

**6.4 Definition** If there exists an injection from a set $A$ into a set $B$, then $A$ is *not of greater cardinality than B*. We write this $A \leq_c B$. If $A \leq_c B$ and $B \leq_c A$, then $A$ *is of the same cardinality as* $B$ which we write $A =_c B$. If $A \leq_c B$ but not $A =_c B$, then $A$ *is of smaller cardinality than* $B$ which we write $A <_c B$.

The relation 'being of the same cardinality as' must be an *equivalence* relation. We shall prove now that this really is the case. In the proof we shall use the concept of composition of two functions. So we must first define this concept and prove that it exists.

Let $R_1$ and $R_2$ be two relations. The *composition* of $R_1$ and $R_2$ is the relation $R_3$ whose elements are all ordered pairs $(a, c)$ for which there exists such $b$ that $(a, b) \in R_1$ and $(b, c) \in R_2$. We denote the composition of $R_1$ and $R_2$ by $R_1R_2$.

**6.5 Lemma** *Every infinite set has at least one accumulation point in* $\mathcal{T}_a$.

*Proof.* Let $A$ be an infinite set. Every node on the $n$th level of $\mathbf{U}$ can be interpreted as the $n$th approximation of all those $x \in U$ whose paths go through this node. (See Definition 3.1 and Theorem 4.1.) Hence the root of $\mathbf{T}_A$ is the 0th approximation of infinitely many points from $A$. Let's denote the root of $\mathbf{T}_A$, which is also the root of $\mathbf{U}$, by $v_0$. At least one of the successors of $v_0$ is the 1st approximation of infinitely many points from $A$ otherwise $A$ would be finite. Namely, in this case $A$ would be the union of a finite family of finite sets and as such would be finite. The proof for this is the same as in ZF. So let a successor $v_1$ of $v_0$ be the approximation of infinitely many points from $A$. Applying the same argument we conclude that among the successors of $v_1$ there must be a node $v_2$ which is the approximation of infinitely many points from $A$, etc. The infinite path $v_0, v_1, v_2, ...$ is the developing sequence of some $x \in A$. Obviously every of the nodes $v_0$, $v_1$, $v_2$, ... determines a neighborhood of $x$ in which there is at least one $y \in A$, $y \neq x$. Thus $x$ is an accumulation point in $A$. □

**6.6 Theorem** *If the relations $R_1$ and $R_2$ exist then the composition $R_1R_2$ exists.*

*Proof.* The predicate for $x$ to be an element of $R_1R_2$ is

$$\mathbf{P}(x) \Leftrightarrow \exists a\exists b\exists c((a, b) \in R_1 \wedge (b, c) \in R_2 \wedge x = (a, c)).$$

Since $R_1$ and $R_2$ exist, $\mathbf{P}$ is expressible by existent predicates and therefore $\mathbf{P}$ exists by Theorem 4.5. Hence the set $P$, determined for every $n$ by $F(n + 1, P) = \{F(n, x) \mid \mathbf{P}(x)\}$, exists by Theorem 4.7.b. Let $x \in P$. Then for every $n$ there exist $a_n$, $b_n$ and $c_n$ such that $(a_n, b_n) \in R_1$, $(b_n, c_n) \in R_2$ and $F(n + 2, x) = F(n + 2, (a_n, c_n))$. Since $F(n + 2, F(n + 3, x)) = F(n + 2, x)$ we have for every $n$

$$F(n + 2, F(n + 3, (a_{n+1}, c_{n+1}))) = F(n + 2, F(n + 3, x)) = F(n + 2, x) = F(n + 2, (a_n, c_n)). \quad (1)$$

By Theorem 4.16 $$F(n + 2, F(n + 3, (a_{n+1}, c_{n+1}))) = (F(n, a_{n+1}), F(n, c_{n+1})) \quad (2)$$

and $$F(n + 2, (a_n, c_n)) = (F(n, a_{n+1}), F(n, c_{n+1})) \quad (3)$$

for every $n$. By Combining (1), (2) and (3) we obtain $F(n, a_{n+1}) = F(n, a_n)$ for every $n$. Hence the approximations $F(0, a_0)$, $F(1, a_1)$, $F(2, a_2)$, ... form the developing sequences of some $a$. Namely

$$F(n, F(n + 1, a)) = F(n, F(n +1, a_{n+1}) = F(n, a_{n+1}) = F(n, a_n) = F(n, a)$$

which by Theorem 3.5 means that $a$ exists. As $F(n, c_{n+1}) = F(n, c_n)$ also holds for every $n$ because of (1), (2) and (3) the approximations $F(0, c_0)$, $F(1, c_1)$, $F(2, c_2)$, ... form the developing sequence of some $c$. For the sequence $b_0$, $b_1$, $b_2$, ... there exists at least one accumulation point $b$ by Lemma 6.5, or $b_n = b$ from some $n$ onwards. Hence $(a, b) \in R_1$ and $(b, c) \in R_2$. Therefore $x = (a, c) \in R_1R_2$. Hence $\mathbf{P}(x)$ is true and $R_1R_2$ exists on the base of Theorem 4.7.e. □

**6.7 Lemma** *If $f_1$ and $f_2$ are injections then the composition $f_1f_2$ is an injection.*

We shall omit the proof since it is the same as in ZF. □

**6.8 Theorem** *The relation $A =_c B$ is an equivalence relation in the universal set $U$.*

*Proof.* Reflexivity: the identity map on $A$ guarantees that $A =_c A$. Symmetry: if $A =_c B$ then there exists an injection from $A$ into $B$ and an injection from $B$ into $A$ by Definition 6.4. Thus $B =_c A$ as well. Transitivity: if $A =_c B$ and $B =_c C$ then there exists an injection $f$ from $A$ into $B$ and an injection $g$ from $B$ into $C$. The composition $fg$ exists by Theorem 6.6 and it is an injection from $A$ into $C$ by the previous two lemmas. This implies $A \leq_c C$. Since $=_c$ is symmetric, there exists also an injection $g'$ from $C$ into $B$ and an injection $f'$ from $B$ into $A$. Thus the composition $g'f'$ is an injection from $C$ into $A$ which implies $C \leq_c A$. Hence $A =_c C$. □

**6.9 Remark** An equivalence relation between the elements of some set separates this set in ZF into mutually disjoint equivalence classes which, of course, are sets. The relation $=_c$ is in ZF an equivalence relation. If the set of all sets existed in ZF, then $=_c$ would separate it into mutually disjoint equivalence classes. Every equivalence class would represent a precisely determined cardinal number. But the set of all sets does not exist in ZF. In AS the set of all sets does exist. It is the universal set $U$. However, $=_c$ does not separate $U$ into mutually disjoint equivalence classes. Here is an example. Let $A_n = \{\{^n\{\}\}^n,\{^n\{\},\{\{\}\}\}^n\}$, where $\{^n$ and $\}^n$ designate the concatenation of $n$ brackets { and } respectively. Obviously $A_n$ is a hereditarily finite set for every $n \geq 0$ and its cardinal number is 2. The sets $A_n$, $n \geq 0$, are therefore elements of the equivalence class $E_2 \subseteq U$ which is the smallest class containing all sets having cardinal number 2. The set $\infty$ is an accumulation point of the sets $A_n$ in $\mathcal{T}_a$ (See Definition 4.21, the text following it, and Theorem 5.12). Hence $\infty$ is an element of $E_2$ although its cardinal number is 1. Thus $E_2$ and the equivalence class $E_1$, which represents the cardinal number 1, are not mutually disjoint. This is no surprise since the set of all sets having some property contains in AS also other sets as its elements by Theorem 4.7.d.

# 7. The Cardinalities of any two Sets can be compared

If the relation $\leq_c$ is to make sense, then it must be an *ordering*. That is, it must be reflexive, antisymmetric with regard to $=_c$ ($A \leq_c B$ and $B \leq_c A$ implies $A =_c B$) and transitive. Besides it must be *total*, which means that for any two sets $A$ and $B$ holds either $A \leq_c B$ or $B \leq_c A$. The proof that $\leq_c$ is an ordering is a simplification of the proof of Theorem 6.8. But the proof that $\leq_c$ is total is a harder nut. First, Lemma 6.2 implies that $A \leq_c P$ holds for any set $A$ and any perfect set $P$. Hence if at least one of two sets contains a perfect subset, then the comparison of cardinality between these two sets is possible. However, if this is not the case, then, as we shall prove in the sequel, these two sets can be well-ordered in AS without the assumption of the Axiom of Choice. And it can be proved, in the same way as in ZF, that there always exists an injection from at least one of two well-orderable sets into the other. It is interesting that these two cases are complementary in AS. A set can be well-ordered precisely when it does not contain a perfect subset.

We shall see that a sensible definition of well-ordering in AS is slightly different from the one in ZF. In ZF a set is *well-ordered* by a relation $R$ if it is linearly ordered by $R$ and if every of its nonempty subsets has the minimal element with regard to $R$. As it is well known, a *linear ordering* of a set $A$ is an asymmetric, transitive and total relation on $A$. The *minimal* element of $A$ with regard to a linear ordering $<$ is such $a \in A$ that $a < x$ for every $x \in A$, $x \neq a$. Analogously, the *maximal* element of $A$ is such $a' \in A$ that $x < a'$ for every $x \in A$, $x \neq a'$. We shall denote the maximal and the minimal elements of $A$ by $\max(A)$ and $\min(A)$ respectively.

A well-ordering, as defined above, exists in ZF on every set. But this is so only if we adopt the Axiom of Choice. And even then it is not at all clear how could we, for example, well-order $\mathcal{P}(\boldsymbol{N})$. On the other hand, well-ordering, defined in this way, exists in AS for every set without the adoption of any additional axiom. Namely, by Definition 2.1 every HF set can be uniquely represented by a string of curly braces. By Definition 2.9 we can lexicographically order these strings and obtain the following result:

**7.1 Theorem** *There exists a linear ordering on every set with regard to which every nonempty subset of this set has the minimal element.*

*Proof.* We shall construct this linear ordering $<$ on the base of the lexicographical ordering of the approximations of sets which we shall denote by $<_{Lex}$. If $x$ and $y$ are two different sets, then let $x < y$ if there exists such $n$ that $F(n, x) = F(n, y)$ and $F(n + 1, x) <_{Lex} F(n + 1, y)$. The relation $<$ linearly orders the universal set $U$ and consequently it linearly orders every subset of $U$. This means that it linearly orders every set since subsets of $U$ are precisely the elements of $U$. (See the proof of Theorem 4.20). We obtain the described linear ordering of some set $A$ in such a way that we linearly order the successors of every node in $\mathbf{T}_A$ from left to right by $<_{Lex}$. If $a_1, a_2 \in A$ and $a_1 \neq a_2$ then $a_1 < a_2$ precisely when the path for $a_1$ in $\mathbf{T}_A$ is more to the left than the path for $a_2$. This implies that $\min(A)$ is represented by the leftmost infinite path in $\mathbf{T}_A$. Hence $\min(A)$ exists and it is an element of $A$. □

**7.2 Definition** The linear ordering from the previous theorem is called the *lexicographical ordering*.

Note that lexicographical ordering exists on every set and not only on HF sets. Therefore, every set in AS can be easily well-ordered with regard to the classical definition of well-ordering. However, in AS this type of ordering does not guarantee a well-ordering as we understand it in the intuitive sense. Namely, to well-order a set $A$ intuitively means to arrange the elements of $A$ in a line in such a way that we can remove $\min(A)$ from $A$ obtaining the set $A \setminus \{\min(A)\}$. Then we can remove the minimal element of $A \setminus \{\min(A)\}$ thus obtaining still smaller set ..., etc. In this way we can remove all the elements of $A$ one by one obtaining at each step a smaller set. In ZF we can always do this, since for every $S \subseteq A$ there exists the set $S \setminus \{\min(S)\}$ by the Axiom of Separation. However, in AS the minimal element of some $S \subseteq A$, with regard to some well-ordering $R$ of $A$, can be an accumulation point in $S$ in $\mathcal{T}_a$. Thus it cannot be removed from $S$ unless we do not previously remove from $S$ infinitely many points which accumulate to $\min(S)$ from above. A necessary condition to carry out the described procedure is that for every $x \in A$ there exists a neighborhood of $x$ in $\mathcal{T}_a$ in which there is no $y \in A$ such that $xRy$. These observation implies that the only sensible definition of well-ordering in AS would be the following one:

**7.3 Definition** A point $x \in A$ is *upward* (*downward*) *isolated* in $A$ with regard to a linear ordering $R$ of $A$ if there exists a neighborhood of $x$ in $\mathcal{T}_a$ in which there is no $y \in A$ such that $xRy$ ( $yRx$). A set $A$ is *well-ordered* by $R$ if

a) $A$ is linearly ordered by $R$
b) every nonempty subset of $A$ has the minimal element with regard to $R$
c) every $x \in A$ is an upward isolated point in $A$ with regard to $R$

If there exists a relation which well-orders $A$ then we say that $A$ *can be well-ordered*, or that $A$ is a *well-orderable* set.

Notice that if $a \in A$ is an upward isolated point in $A$ with regard to some well-ordering $R$, then $a$ is not necessarily an isolated point in $A$ as $a$ may not be downward isolated in $A$. The set $\boldsymbol{N}$ is well-ordered by the lexicographical ordering defined in 7.2, but the set $\mathcal{P}(\boldsymbol{N})$ is not because some of its elements are not upward isolated with regard to lexicographical ordering. In fact $\mathcal{P}(\boldsymbol{N})$ cannot be well-ordered at all as it is shown below.

**7.4 Theorem** *A set which is perfect in* $\mathcal{T}_a$ *cannot be well-ordered.*

*Proof*. Let a set $P$ be perfect in $\mathcal{T}_a$ and let's assume that $<$ is a well-ordering of $P$. Since, by Definition 7.3.c, $\min(P)$ must be an upward isolated point in $P$ with regard to $<$, there exists such neighborhood of $\min(P)$ in $\mathcal{T}_a$ in which there is no $y \in P$ such that $\min(P) < y$. As there is also no $y \in P$, $y < \min(P)$, the point $\min(P)$ is isolated in $P$. But there are no isolated points in $P$ since $P$ is perfect. Hence our assumption that a perfect set can be well-ordered leads to a contradiction. □

**7.5 Lemma** *If a subset of some set A cannot be well-ordered, then A cannot be well-ordered either.*

*Proof*. Let $B \subseteq A$. Every well-ordering of $A$ is an extension of some well-ordering of $B$. So, if there is no well-ordering on $B$ then there is no well-ordering on $A$. □

**7.6 Corollary** *The set* $\mathcal{P}(\boldsymbol{N})$ *cannot be well-ordered.*

*Proof*. The set $\mathcal{P}(\boldsymbol{N})$ contains a perfect subset by Theorem 5.13. So Theorem 7.4 and Lemma 7.5 imply that $\mathcal{P}(\boldsymbol{N})$ cannot be well-ordered. □

The lexicographical ordering of some set $A$ is only one of many possible linear orderings of $A$ since we can linearly order the successors of every node of $\mathbf{T}_A$ in an arbitrary way obtaining thus some linear ordering of $A$. We shall call such ordering an *arrangement* of $A$ since it is obtained by an arrangement of $\mathbf{T}_A$. For every nonempty subset $B$ of every set $A$ there exists min($B$) with regard to any arrangement of $A$. The path for min($B$) is the leftmost infinite path of $\mathbf{T}_B$ in such an arrangement of $\mathbf{T}_A$. Now we shall prove that any well-ordering of some set is obtained by an arrangement of this set.

**7.7 Theorem** *If no arrangement well-orders some set, then this set contains a perfect subset.*

*Proof.* Let no arrangement of a set $A$ well-order $A$. Then there are the following three mutually exclusive possibilities

a) There exists a level in $\mathbf{T}_A$ such that from every one of its nodes develops in $A$ a set which can be well-ordered by some arrangement
b) On every level of $\mathbf{T}_A$ there exists exactly one node from which develops in $A$ a set which is not well-orderable by any arrangement
c) On some level of $\mathbf{T}_A$ there exist at least two nodes from each of which develops in $A$ a set which is not well-orderable by any arrangement

**a)** Let $v_1, v_2, ..., v_k$ be all the nodes on this level, let $A_1, A_2, ..., A_k$ be the sets which develop from these nodes in $A$, and let the arrangements which well-order these sets be $R_1, R_2, ..., R_k$ respectively. Then there exists an arrangement $R$ which is an extension of every $R_i$ for $i = 1, 2, ..., k$ and which well-orders $A$.
1) Obviously $R$ matches with $R_i$ inside $A_i$ for $i = 1, 2, ..., k$. And if $a \in A_i$ and $a' \in A_j$, where $i < j$, then we put $aRa'$. In this way $R$ linearly orders $A$.
2) Let $B$ be a nonempty subset of $A$ and let $m = \min\{i \mid B \cap A_i \neq \emptyset\}$. Then min($B$) with regard to $R$ is equal to min($B \cap A_m$) with regard to $R_m$.
3) If $a \in A$ then $a \in A_i$ for some $i = 1, 2, ..., k$. Let $n$ be the level which satisfies a). If $m > n$ then in the $m$th neighbourhood of $a$ in $A$ exist only points from $A_i$ since $A_1, A_2, ..., A_k$ are mutually separated sets. And since $A_i$ is well-ordered by $R_i$ there exists such $m > n$ that in the $m$th neighbourhood of $a$ in $A$ there is no $a' \in A_i$ such that $aR_ia'$. Hence $a$ is upward isolated in $A$ with regard to $R$.
Items 1), 2) and 3) imply that $A$ is well-ordered by the arrangement $R$. But this contradicts the assumption that no arrangement well-orders $A$. Hence the possibility a) cannot occur.

**b)** Let's denote such node on the $n$th level by $v_n$. From the predecessor of $v_n$ must develop in $A$ a set which is not well-orderable by any arrangement, otherwise, as it is not difficult to see, the set which develops from $v_n$ would also be well-orderable by this arrangement. But this is not the case. Therefore $v_{n-1}$ must be the predecessor of $v_n$ for every $n > 0$, and the sequence $v_0, v_1, v_2, v_3, \ldots$ forms an infinite path in $\mathbf{T}_A$. This path is the developing sequence of some $x \in A$. Let's arrange the successors of every $v_n$ in such a way that $v_{n+1}$ is the rightmost of them. Then the developing sequence of $x$ is the rightmost infinite path in $\mathbf{T}_A$. That is, $x = \max(A)$ with regard to this arrangement.

Now we shall rearrange also all the other nodes of $\mathbf{T}_A$ in such a way that $A$ becomes well-ordered. We shall denote this arrangement of $A$ by $R$. Let $w_{n,1}, w_{n,2}, ... w_{n,k(n)}$ be all the successors of $v_n$ with the exception of $v_{n+1}$. By the assumption made in b) develops in $A$ from every node $w_{n,i}$ a set $A_{n,i}$ which is well-orderable by some arrangement $R_{n,i}$. The path for every $y \in A$, $y \neq x$, separates from the path for $x$ at some level $n$. That is $F(n, y) = F(n, x) = v_n$ and $F(n + 1, y) \neq F(n + 1, x) = v_{n+1}$. Hence $F(n + 1, y) = w_{n,i}$ for some $i$ between 1 and $k(n)$, which means that $y \in A_{n,i}$. If $y, y' \in A_{n,i}$ for some $i$ then let $(y, y') \in R$ if and only if $(y, y') \in R_{n,i}$. In other words, $R_{n,i} \subseteq R$ for every $(n, i)$ where $1 \leq i \leq k(n)$. If, however, $y \in A_{n,i}$ and $y' \in A_{n',i'}$ for some $(n, i) \neq (n', i')$, then let $(y, y') \in R$ if and only if $(n, i) < (n', i')$. Here $(n, i) < (n', i')$ denotes that $n < n'$ or that $n = n'$ and $i < i'$. Since the developing sequence of $x$ is the rightmost infinite path in $\mathbf{T}_A$, we have $(y, x) \in R$ for every $y \in A$, $y \neq x$. Hence $R$ is a total linear ordering of $A$. Every $y \in A$, $y \neq x$, is an upward isolated point in $A$ with regard to $R$ by a similar argument as we have made in item 3) of a). Also $x$ is an upward isolated point in $A$ with regard to $R$ since $x = \max(A)$ with regard to $R$. Let $B$ be a nonempty subset of $A$. Then there exists the smallest pair $(n, i)$ such that $A_{n,i} \cap B \neq \emptyset$. Hence the minimal element of $B$ in $A$ with regard to $R$ is the minimal element of $A_{n,i} \cap B$ in $A_{n,i}$ with regard to $R_{n,i}$. All this implies that the arrangement $R$ well-orders $A$ which contradicts the initial assumption. Hence the possibility b) cannot occur

**c)** Let's denote the sets which develop in $A$ from these two nodes by $A_1$ and $A_2$. Then, by the conclusions we made in discussing a) and b), in each of the developing trees of $A_1$ and $A_2$ there exists a level having two nodes from each of which develops in $A$ a set which is not well-orderable by any arrangement. If we repeat this argument recursively, we can see that $\mathbf{T}_A$ contains a perfect binary tree. This tree is the developing tree of some perfect set $P$ since every $x \in P$ is an accumulation point in $P$. Namely, for every $x \in P$ and every $n$ there exists $y \in P$, $y \neq x$, such that $F(n, x) = F(n, y)$.

Therefore, if $A$ is not well-orderable by any rearrangement of the nodes of $\mathbf{T}_A$, then the only possibility which can occur is c), and consequently $A$ contains a perfect subset. □

**7.8 Theorem** *A set can be well-ordered precisely when it does not contain a perfect subset.*

*Proof.* If a set cannot be well-ordered by any arrangement then, by the previous theorem, it contains a perfect subset. And if this is the case then, by Theorem 7.4 and Lemma 7.5, it cannot be well-ordered at all. And if this is the case then it certainly cannot be well-ordered by any arrangement. This circular argument implies that the statements 'a set cannot be well-ordered' and 'a set contains a perfect subset' are logically equivalent. Hence the statements 'a set can be well-ordered' and 'a set does not contain a perfect subset' are logically equivalent as well. □

So, till now we have proved that the cardinalities of two sets are comparable if at least one of these two sets contains a perfect subset. But what if this is not the case. Then both sets are well-orderable by Theorem 7.8. We shall prove that in this case there always exists an injection from at least one of such two sets into another. These results are classical but we must prove them again, since we have defined well-ordered sets in a new way, and since we must avoid all the conclusions made on the base of the Axiom of Separation which, as we know, does not hold in AS.

**7.9 Definition** Let $A$ be linearly ordered by $<$ and let $S \subseteq A$. If $x \in S$ and $y < x$ imply $y \in S$ for arbitrary $x, y \in A$, then $S$ is an *initial segment* of $A$. If $S$ is a proper subset of $A$, then $S$ is a *proper initial segment* of $A$.

**7.10 Theorem** *Let a set A be well-ordered by < and let S be a proper initial segment of A. Then there exists such* $a \in A$ *that* $S = \{x \mid x \in A \wedge x < a\}$.

*Proof.* By Definition 4.10 and Theorem 4.11 the set $A \setminus S$ exists. Since $S$ is a proper subset of $A$, the set $A \setminus S$ is not empty, and as a subset of a well-orderable set it is well-orderable by Lemma 7.5. Let $a = \min(A \setminus S)$ with regard to $<$. Now, if $x < a$ then $x \notin A \setminus S$. Hence $x \in S$. If $x \geq a$ then $x \notin S$, otherwise $a \in S$ as $S$ is an initial segment of $A$. But this is not possible. Thus $S = \{x \mid x \in A \wedge x < a\}$.□

We shall denote the proper initial segment $S = \{x \mid x \in A$ and $x < a\}$ by $A[a]$. It is evident from the previous proof that $a$ must be an isolated point in $A$ in $\mathcal{T}_a$. If $a$ were an accumulation point, then $a$ would be upward isolated but not downward isolated in $A$, since $A$ is well-ordered. Thus $a$ would be an accumulation point in $S$ as well, and hence $a \in S$, by Theorem 5.3, which is not the case as we have just proved.

**7.11 Definition**

a) Let $A_1$ and $A_2$ be linearly ordered by relations $R_1$ and $R_2$ respectively. An *isomorphism* $f : A_1 \to A_2$ is a bijection which preserves linear ordering. This means, if $x_1, x_2 \in A_1$ and $x_1R_1x_2$, then $f(x_1)R_2f(x_2)$. If there exists an isomorphism between $A_1$ and $A_2$ then $A_1$ and $A_2$ are *isomorphic* to each other.

b) Let a set $A$ be linearly ordered by $<$. The *successor of* $x \in A$ *in* $A$ with regard to $<$ is $x' \in A$, $x < x'$, such that there exists no $y \in A$ such that $x < y < x'$. If $x'$ is the successor of $x$ in $A$ then $x$ is the *predecessor* of $x'$ in $A$.

The concept of successor and predecessor as defined above should not be mixed up with successor and predecessor of a node in a developing tree. The meaning of these words will be always clear from the context. It is not difficult to see that the successor and the predecessor of an element of a linearly ordered set are uniquely determined if they exist. Every element of a well-ordered set $A$ is upward isolated in $A$. Hence every element of $A$, with the exception of $\max(A)$, has its successor in $A$. However, only an isolated point in $A$, with the exception of $\min(A)$, has its predecessor in $A$.

**7.12 Theorem** *At least one of two well-ordered sets is isomorphic to an initial segment of the other.*

*Proof.* Let $A$ and $B$ be well-ordered sets by $<_A$ and $<_B$ respectively, and let $A[x] \sim B[y]$ denote that $A[x]$ is isomorphic to $B[y]$. First we shall prove that there exists the set of precisely all ordered pairs $(x, y)$ where $x \in A$, $y \in B$ and $A[x] \sim B[y]$. Let the predicate $\mathbf{P}(u)$ denote that $u$ is such an ordered pair and let $z \in P$ where $P$ is the smallest set which contains all such pairs. Then, by Theorems 4.7.c and Theorem 4.16, for every $n$ holds

$$F(n + 2, z) = F(n + 2, (x_n, y_n)) = (F(n, x_n), F(n, y_n)) \qquad (1)$$

where $\mathbf{P}((x_n, y_n))$ is true. Thus $F(n + 3, z) = F(n + 3, (x_{n+1}, y_{n+1})) = (F(n + 1, x_{n+1}), F(n + 1, y_{n+1}))$ holds for every $n$. Hence

$$F(n + 2, F(n + 3, z)) = F(n + 2, (F(n + 1, x_{n+1}), F(n + 1, y_{n+1}))) =$$

$$= (F(n, F(n + 1, x_{n+1})), F(n, F(n + 1, y_{n+1}))) \qquad (2).$$

And because of $F(n + 2, F(n + 3, z)) = F(n + 2, z)$ for every $n$, combining (1) and (2) gives

$$F(n, F(n + 1, x_{n+1})) = F(n, x_n) \quad \text{and} \quad F(n, F(n + 1, y_{n+1})) = F(n, y_n).$$

Therefore $F(0, x_0)$, $F(1, x_1)$, $F(2, x_2)$, ... and $F(0, y_0)$, $F(1, y_1)$, $F(2, y_2)$, ... are the developing sequences of some $x' \in A$ and some $y' \in B$ respectively. Thus $z = (x', y') \in A \times B$. If $\mathbf{P}((x', y'))$ is not true then $x'$ must be an upper bound of all those $x \in A$ for which there exists $y \in B$ such that $\mathbf{P}((x, y))$ is true. Namely, if there is such $x'' >_A x'$ that $\mathbf{P}((x'', y''))$ for some $y'' \in B$, then there is an isomorphism $f : A[x''] \to B[y'']$. Hence $A[x'] \sim B[f(x')]$ and $\mathbf{P}((x', f(x')))$ is true. Since $\mathbf{P}((x', y'))$ is not true, we have $y' \neq f(x')$. If $x \in A$ is close enough to $x'$ in $\mathcal{T}_a$ or $x = x'$, then $x <_A x''$ as $x' <_A x''$. Hence $x$ is in the range of $f$. The more $x$ is close to $x'$ the more $f(x)$ is close to $f(x')$ because of the continuity of $f$. Such $f(x)$ cannot be close to $y'$ since $f(x') \neq y'$. Thus in a small enough neighborhood of $(x', y')$ there are no pairs $(x, y)$ for which $\mathbf{P}((x, y))$ is true and so $(x', y') \notin P$. This contradiction shows that $x'$ must be larger then every $x \in A$ for which $\mathbf{P}((x, y))$ is true for some $y \in B$. In a similar way we conclude that $y'$ must be larger then every $y \in B$ for which $\mathbf{P}((x, y))$ is true for some $x \in A$. As in arbitrary small neighborhood of $(x', y')$ there exist pairs $(x, y)$ for which $\mathbf{P}((x, y))$ is true, $x'$ is the precise upper bound of such $x \in A$ and in $y'$ is the precise upper bound of such $y \in B$. Thus there exists at most one pair $(x', y') \in P$ such that $\mathbf{P}((x', y'))$ is not true.

The assignments $x \to y$ and its inverse $y \to x$, where $\mathbf{P}((x, y))$ is true, $x <_A x'$, and $y <_B y'$, are unique otherwise two distinct initial segments of $A$ (or $B$) would be isomorphic to each other which is not possible. The proof for this is the same as in ZF (See, for example, [H&J], p. 105, Corollary 1.5). This assignment is some isomorphism $g$ which contains all pairs $(x, y)$ where $x <_A x'$ and $y = g(x)$. It also contains the pair $(x', y')$ since the $(n+1)$th approximation of $g$ contains the $n$th approximation of $(x', y')$ for every $n$. Namely, the $n$th approximation of $(x', y')$ is the same as the $n$th approximation of $(x_n, y_n) \in g$ since $(x_n, y_n)$ is in the $n$th neighborhood of $(x', y')$. Hence $g$ is bijective on every $x <_A x'$ and also on $x'$ which is mapped by $g$ to $y'$. Therefore $g$ is an isomorphism from $A[x']$ to $B[y']$ and hence $\mathbf{P}((x', y'))$ holds true. This implies that the set of all ordered pairs $(x, y)$, where $x \in A$, $y \in B$ and $A[x] \sim B[y]$, exists.

We shall now prove that for $x'$ and $y'$ discussed above holds either $x' = \max(A)$ or $y' = \max(B)$. If neither of these were true, then there would exist the successor $x''$ of $x'$ in $A$. Similarly there would exist the successor $y''$ of $y'$ in $B$. As $g$ is an isomorphism between $A[x'] \cup x'$ and $B[y'] \cup y'$, it would be an isomorphism from $A[x'']$ to $B[y'']$. Therefore $\mathbf{P}((x'', y''))$ would be true. This, however, is not possible since $x'$ is an upper bound of all $x \in A$ for which $\mathbf{P}((x, y))$ is true. Thus either $x' = \max(A)$ or $y' = \max(B)$, which implies that one of well-ordered sets $A$ and $B$ is always isomorphic to an initial segment of the other. □

This theorem enable us to compare the cardinalities of arbitrary well-order sets.

**7.13 Definition** Let $A$ and $B$ be well-ordered sets. If $A$ is isomorphic to an initial segment of $B$ then we shall write $A \leq_o B$. If $A \leq_o B$ and $B \leq_o A$ then we shall write $A =_o B$. If $A$ is isomorphic to a proper initial segment of $B$ we shall write $A <_o B$.

If $A \leq_o B$ and $B \leq_o A$ then $A$ is isomorphic to $B$. The proof of this is the same as in ZF. (See, for example, [H&J], p. 105, Theorem 1.3.) The proof is based on the fact that no well-ordered set is isomorphic to a proper initial segment of itself. If $A$ and $B$ are finite sets then $A \leq_o B$ is equivalent to $A \leq_c B$.

**7.14 Lemma** *If there is a bijection between a well-ordered set A and a set B, then B can be well-ordered.*

*Proof.* Let a relation $R$ well-order $A$ and let $f$ be a bijection from $A$ to $B$. Then the composition $f^{-1}Rf$ exists by Theorem 6.6. Let's prove that $f^{-1}Rf$ well-orders $B$ by Definition 7.3.

a) The relation $f^{-1}Rf$ linearly orders $B$. Namely, let $b_1, b_2 \in B$, $b_1 \neq b_2$, $f^{-1}(b_1) = a_1$ and $f^{-1}(b_2) = a_2$. Since $f$ is a bijection we have $a_1 \neq a_2$. Without loss of generality we can assume that $a_1Ra_2$. Then $(b_1, a_1) \in f^{-1}$, $(a_1, a_2) \in R$ and $(a_2, b_2) \in f$. Therefore $(b_1, b_2) \in f^{-1}Rf$ and so $f^{-1}Rf$ is total on $B$. Because of the bijection between $A$ and $B$ it must also be asymmetric and transitive on $B$ since $R$ is asymmetric and transitive on $A$.

b) Let $B_1 \subseteq B$ and let $f^{-1}(B_1) = A_1 \subseteq A$. There exists $\min(A_1)$ in $A$ with regard to $R$. If we put in a) that $a_1 = \min(A_1)$ and that $a_2$ is an arbitrary other element of $A_1$, we can conclude that $f(\min(A_1))$ equals to $\min(B_1)$ with regard to $f^{-1}Rf$.

c) Let $b$ be an arbitrary element of $B$ and let $a = f^{-1}(b)$. As $A$ is well-ordered by $R$ there exists such small enough neighbourhood $V(n, a)$ of $a$ (See Definition 5.2) in which there is no point $a'$ such that $(a, a') \in R$. Because of the continuity of $f^{-1}$ there exists such small enough neighbourhood $V(m, b)$ of $b$ that $f^{-1}(V(m, b)) \subseteq V(n, a)$. If in $V(m, b)$ there existed some $b'$ such that $(b, b') \in f^{-1}Rf$ then $f^{-1}(b') \in V(n, a)$ and consequently $(a, f^{-1}(b')) \in R$ by a). But this would be a contradiction. Hence $b$ is an upward isolated point in $B$ with regard to $f^{-1}Rf$. □

**7.15 Theorem** *The cardinalities of two arbitrary sets A and B are always comparable. Precisely one of the following three possibilities is true: $A <_c B$ or $A =_c B$ or $B <_c A$.*

*Proof.* Since the relations $A <_c B$, $A =_c B$ and $B <_c A$ are mutually exclusive it is enough if we prove that at least one of them holds in any case. There are the following tree mutually exclusive possible cases:

a) The sets $A$ and $B$ are well-orderable. Then, by Theorem 7.12, at least one of these two sets is isomorphic to an initial segment of the other. Without loss of generality we can presume that $A$ is isomorphic to an initial segment of $B$. Then there exists an injection from $A$ into $B$. That is $A \leq_c B$. If there exists also an injection from $B$ into $A$, then $A =_c B$, otherwise $A <_c B$.

b) One of the sets $A$ and $B$ is well-orderable and the other is not. Again we can presume, without loss of generality, that $A$ is well-orderable. Theorem 7.8 implies that $B$ contains a perfect subset and $A$ does not. Hence, by Lemma 6.2, there exists an injection from $A$ into $B$. However, no injection exists from $B$ into $A$ since every subset of $A$ is well-orderable by Lemma 7.5 and $B$ cannot be in bijection with a well-orderable set by Lemma 7.14. Thus $A <_c B$.

*c*) None of the sets $A$ and $B$ is well-orderable. Then they both contain a perfect subset by Theorem 7.8. Hence there exists an injection from $A$ into $B$ and an injection from $B$ into $A$ by Lemma 6.2. Therefore $A =_c B$. □

# 8. Determination of Cardinality by Cantor-Bendixson Rank and Cantor-Bendixson Degree

**DERIVED SETS** Every set in AS can be interpreted as a subset of [0, 1] in ZF (See the beginning of section 5). This interpretation is not unique. Every set in AS is closed in the topology of approximations $\mathcal{T}_a$ by Theorem 5.3. This topology is equivalent to the interval topology in ZF. Thus every ZF image of some AS set is closed in the interval topology. Every bijection in AS is a homeomorphism in $\mathcal{T}_a$ by Theorem 5.10. Hence in AS an injection from a set $A$ into a set $B$ exists precisely when a topological embedding of the corresponding images in [0, 1] exists. If the sets $A$ and $B$ are well-orderable in AS then, by Theorem 7.8, neither of them contains a subset which is perfect in $\mathcal{T}_a$. Hence their ZF images do not contain a subset which is perfect in the interval topology. Since these images are closed sets, they are countable in ZF (See, for example, [H&J], p. 190, Theorem 4.5). As it is known an embedding of a closed countable set into another such set depends on Cantor-Bendixson rank and Cantor-Bendixson degree of both sets. All these implies that whether in AS an injection from one well-orderable set into another well-orderable set exists depends on Cantor-Bendixson rank and Cantor-Bendixson degree of the images of these two sets in [0, 1]. If we could transfer this rank and degree into AS then we could with their help determine the cardinality of well-orderable sets in AS. So, let's remind ourselves what these two concepts are.

Let $A$ be some set of real numbers in ZF. Then the *derived* set of $A$, denoted by $A'$, is the set of all accumulation points in $A$ in the interval topology. Let $A^{(0)} = A$, let $A^{(\alpha+1)} = (A^{(\alpha)})'$ for every ordinal $\alpha$, and let $A^{(\alpha)} = \cap_{\beta <_o \alpha} A^{(\beta)}$ for a limit ordinal $\alpha$. If $A$ is countable then there exists such smallest $\alpha$ for which $A^{(\alpha+1)} = \emptyset$. This $\alpha$ is the *Cantor-Bendixson rank* of $A$. In this case $A^{(\alpha)}$ is a finite set of isolated points on the real line. The number of these points is the *Cantor-Bendixson degree* of $A$. Two closed countable sets can be embedded into each other precisely when they have the same Cantor-Bendixson rank and the same Cantor-Bendixson degree. Translated into AS this would mean that two sets have the same cardinality precisely when they have the same Cantor-Bendixson rank and the same Cantor-Bendixson degree. If we want to prove something like this, we must prove the existence of derived sets in AS and the existence of arbitrary large ordinals in AS representing Cantor-Bendixson rank. Since in ZF a derived set exists by the Axiom of Separation which is not valid in AS we must prove the following:

**8.1 Theorem** *The set of all accumulation points of an arbitrary set exists.*

*Proof.* Let $A$ be an arbitrary set. Then the corresponding predicate for $x$ to be an element of the derived set $A'$ is:

$\mathbf{P}(x) \Leftrightarrow \forall n \exists m \exists y_m(\mathbf{N}(n) \wedge \mathbf{N}(m) \wedge n \leq m \wedge y_m \in A \wedge F(m, x) = F(m, y_m) \wedge F(m + 1, x) \neq F(m + 1, y_m))$.

Thus $\mathbf{P}(x)$ expressed in words is: $x$ is an accumulation point in $\mathcal{T}_a$ in $A$. The predicate $\in$ exists by Example 4.6. The predicate **N** is a computable predicate on HF sets. Hence, by Theorem 4.7.b, there exists a set $P$, determined for every $n$ by $F(n + 1, P) = \{F(n, x) \mid \mathbf{P}(x)\}$, for which the following holds. If $x \in P$, then, by Theorem 4.7.c, there exists such $y_n$ for every $n$, not necessarily different from $x$, that $F(n, x) = F(n, y_n)$ and $y_n$ is an accumulation point in $\mathcal{T}_a$ in $A$. Hence $y_n \in A$ for every $n$. Therefore $x$ is either equal to some $y_n$ or it is an accumulation point of the points $y_0$, $y_1$, $y_2$, .... each of which is an element of $A$. In each case $x$ is an accumulation point in $A$ and consequently $\mathbf{P}(x)$ is true. Thus the derived set of $A$ exists by Theorem 4.7.e. □

**VON NEUMANN ORDINALS** As we have already stated, in order to introduce Cantor-Bendixson rank there should exist arbitrary large ordinals in AS. An *ordinal* is in ZF defined as a transitive set which is well-ordered by the relation $\in$. A set is *transitive* if every one of its elements is also its subset. We shall denote these ordinals as *von Neumann ordinals*, or shortly *N-ordinals*, in order to distinguish them from the ordinals we shall define in AS. Namely, we shall soon see that there are not enough N-ordinals in AS for our purposes.

Some N-ordinal is *finite* if it has finitely many elements. Finite N-ordinals can be interpreted as natural numbers. That is, $0 = \emptyset = \{\}$ is a finite N-ordinal, and if $\alpha$ is a finite N-ordinal, then $\alpha \cup \{\alpha\}$ is a finite N-ordinal too. In this way we obtain all finite N-ordinals. Hence finite N-ordinals are hereditarily finite sets. Let $n$ be a natural number expressed in Zermelo form. Let $f(\{\}) = \{\}$ and let $f(\{n\}) = f(n) \cup \{f(n)\}$. Then $f$ is an isomorphism in predicate calculus between natural numbers expressed in Zermelo form and finite N-ordinals. That is, $f$ maps 0 in Zermelo form to N-ordinal 0. And if $f$ maps $n$ in Zermelo form to N-ordinal $\alpha$ then $f$ maps $\{n\}$ to $\alpha \cup \{\alpha\} = \alpha + 1$. This implies that for finite N-ordinals hold Peano Axioms and all the rules concerning the sum of natural numbers since all this holds for natural numbers in Zermelo form (See section 2). Hence a natural number in Zermelo form can be in every expression substituted by the corresponding finite N-ordinal. In particular this is true for the first argument of the function $F(n, A)$ and for the rank of HF sets in the following lemma.

**8.2 Lemma** *If m and n are finite N-ordinals, then* $F(n, m) = \min\{m, n\}$.

*Proof.* For every finite N-ordinal $m$ the rank of $m$ is $m$. This implies $F(n, m) = m$ for $n \geq m$. On the other hand we have $F(0, m) = \{\} = 0$ for every finite N-ordinal $m$. And if $F(n, m) = n$ for some finite N-ordinal $n$ and every $m > n$, then we have

$$F(n + 1, m) = F(n + 1, \{0, 1, 2, ..., m - 1\}) = \{F(n, 0), F(n, 1), \ldots, F(n, m - 1)\} = \{0, 1, ..., n\} = n + 1$$

for every finite N-ordinal $m > n + 1$. Thus $F(n, m) = n$ holds true for every finite N-ordinal $m > n$ by induction on $n$. □

This lemma could also be proved from Example 1.4.1 by the previously constructed isomorphism between natural numbers expressed in Zermelo form and finite N-ordinals.

The smallest infinite N-ordinal is denoted by $\omega$. This N-ordinal is the set of all finite N-ordinals and exists in ZF by the Axiom of Replacement. Since this axiom is not valid in AS we shall prove the existence of $\omega$ in another way. In the following theorem natural numbers will be expressed as finite N-ordinals.

**8.3 Theorem** *There exists the set of all finite N-ordinals, denoted by ω, whose element is also ω. This set is infinite.*

*Proof.* Let the set $\omega$ be determined by $F(n + 1, \omega) = \{F(n, x) \mid$ "$x$ is a finite N-ordinal"$\}$. This set exists by Theorem 4.7.b and Theorem 4.5 because the predicate "to be a finite N-ordinal" is computable on HF sets. Let $x \in \omega$. Then for every finite N-ordinal $n$ there exists such finite N-ordinal $y_n$ that $F(n, x) = F(n, y_n)$, which by Lemma 8.2 means that $F(n, x) = \min\{n, y_n\}$. Therefore $F(n, x)$ is a finite N-ordinal not larger than $n$. If $F(n, x) <_o n$ for some $n$, and $r$ denotes the rank of a hereditarily finite set then $r(F(n, x)) < n$ and hence $F(n, x) = x$ by Theorem 3.15. In this case $x$ is a finite N-ordinal not larger than $n$. The remaining possibility is that $F(n, x) = n$ for every $n$. This implies

$$F(n + 1, x) = n + 1 = \{i \mid i \le n\} = \{F(n, i) \mid \text{“}i \text{ is a finite N-ordinal”}\}$$

Since this is true for every $n$ we have $x = \omega$. Therefore $\omega$ is the only element of $\omega$ which is not a finite N-ordinal, similarly as $\infty$ is the only element of $\boldsymbol{N}$ which is not a natural number.

The map which maps every finite N-ordinal $n$ to $n + 1$ and $\omega$ to $\omega$, is a bijection from $\omega$ to $\omega \setminus \{0\}$. The set $\omega \setminus \{0\}$ exists by Theorem 4.11, and the described map exists because it is a homeomorphism in $\mathcal{T}_a$. As 0 is a HF set it is an isolated point in $\omega$. Hence the set $\omega \setminus \{0\}$ is a proper subset of $\omega$. Therefore $\omega$ is infinite. □

Hence $\omega$ in AS is not equal to $\omega$ in ZF, since in AS holds $\omega \in \omega$ and in ZF this does not hold. Therefore $\omega$ in AS is not well-ordered by $\in$. The N-ordinal which follows $\omega$ in ZF is $\omega \cup \{\omega\}$. This N-ordinal is the *successor* of $\omega$ and is in ZF different from $\omega$. But in AS we have $\omega \cup \{\omega\} = \omega$. Hence the N-ordinal $\omega$ is the successor of itself in AS. Consequently, $\omega$ is the largest N-ordinal in AS and at the same time it is the set of all N-ordinals in AS. So there are very few N-ordinals in AS. However, this has at least one positive side effect. Although the set of all N-ordinals exists, the following holds:

**8.4 Theorem** *The argument for Burali-Forti paradox fails on N-ordinals.*

*Proof.* This paradox appears in the unaxiomatized set theory, developed by Cantor and Frege, where for every property **P** there exists the set of all objects having the property **P**. In this theory the set of all N-ordinals is well-ordered and transitive and is therefore an N-ordinal. Let's denote this set by $\Omega$. The successor of $\Omega$ is an N-ordinal as well and is equal to $\Omega \cup \{\Omega\}$. However, since $\Omega$ is the set of all N-ordinals, we have $\Omega \cup \{\Omega\} \in \Omega$. And, because $\Omega$ is an N-ordinal it is transitive. Therefore we have $\Omega \cup \{\Omega\} \subseteq \Omega$. Thus $\Omega \in \Omega$, which contradicts the fact that N-ordinals are well-ordered by $\in$. In AS the set of all N-ordinals exists but there is no paradoxical situation. This set is equal to $\omega$. It is transitive but it is not well-ordered by $\in$ since $\omega \in \omega$. □

Now we shall prove that in AS the set of all N-ordinals is in bijection with the set of all natural numbers represented in Zermelo form.

**8.5 Theorem** *The set ω is in bijection with the set* ***N****.*

*Proof.* Every element of $\omega$, with the exception of $\omega$ itself, is a finite N-ordinal, which in turn is a HF set. Hence, by Lemma 6.5 and Theorem 5.12, the only accumulation point in $\omega$ in $\mathcal{T}_a$ is $\omega$. Since $\infty$ is the only accumulation point in $\boldsymbol{N}$ in $\mathcal{T}_a$, there exists a bijection between $\boldsymbol{N}$ and $\omega$, that is, a homeomorphism in $\mathcal{T}_a$ between $\boldsymbol{N}$ and $\omega$ which maps every natural number $n$ in Zermelo form to the N-ordinal $n$ and maps $\infty$ to $\omega$. □

**ORDINAL ARITHMETIC OF VON NEUMANN ORDINALS** As the set of all N-ordinals is so very reduced in AS, the ordinal arithmetic of N-ordinals in AS is very simple. Every N-ordinal, with the exception of $\omega$, is in fact a natural number. For this reason the rules for addition, multiplication and exponentiation for finite N-ordinals are the same as the corresponding rules for natural numbers. The sum and the product of finite N-ordinals are associative and commutative operations connected by the distributive law. And when $\omega$ is involved in the computation, we must take into account that $\omega$ is a *limit* N-ordinal. That is, there is no N-ordinal $\alpha <_o \omega$ such that there is no N-ordinal between $\alpha$ and $\omega$. Knowing all this we apply the following rules of the ordinal arithmetic:

| | sum | product | exponentiation |
|---|---|---|---|
| for any $\alpha$ : | $\alpha + 0 = \alpha$ | $\alpha \cdot 0 = 0$ | $\alpha^0 = 1$ |
| for any $\beta$ : | $\alpha + (\beta + 1) = (\alpha + \beta) + 1$ | $\alpha \cdot (\beta + 1) = \alpha \cdot \beta + \alpha$ | $\alpha^{\beta+1} = \alpha^{\beta} \cdot \alpha$ |
| for a limit $\beta$ : | $\alpha + \beta = \max\{\alpha + \gamma \mid \gamma <_o \beta\}$ | $\alpha \cdot \beta = \max\{\alpha \cdot \gamma \mid \gamma <_o \beta\}$ | $\alpha^{\beta} = \max\{\alpha^{\gamma} \mid \gamma <_o \beta\}$ |

The relation $\alpha <_o \beta$ holds for finite N-ordinals $\alpha$ and $\beta$ precisely when it holds for the corresponding natural numbers. And for every finite N-ordinal $\alpha$ holds $\alpha <_o \omega$. It is not difficult to see that every set of N-ordinals is linearly ordered, that the maximal element of such set always exists, and that this element is equal to $\omega$ if the set is infinite. Now, as $\omega + 1 = \omega \cup \{\omega\} = \omega$, we have, by induction, $\omega + n = \omega$ for every $n <_o \omega$. This implies that the ordinal arithmetic gives in the cases in which $\omega$ is involved the following results:

| | |
|---|---|
| $n + \omega = \omega + n = \omega + \omega = \omega$ | for any finite N-ordinal $n$ |
| $n \cdot \omega = \omega \cdot n = \omega \cdot \omega = \omega$ | for any finite N-ordinal $n > 0$ |
| $\omega^n = \omega$ | for any finite N-ordinal $n > 0$ |
| $n^{\omega} = \omega$ | for any finite N-ordinal $n > 1$ |

We can see that here $\omega$ has the same role as has $\infty$ in the arithmetic of the set $\boldsymbol{N}$.

**ORDINALS IN AS** According to the above conclusions we cannot introduce Cantor-Bendixson rank as an N-ordinal since there are not enough N-ordinals in AS for this purpose. Even more. There does not seem to exist a natural representative set of any well-ordered set in AS which would be the same for isomorphic well-ordered sets. An idea might be to define ordinals as equivalence classes where the elements of each class are all well-ordered sets isomorphic to each other. But such classes overlap in the same way as the classes described in Remark 6.9. Therefore it seems that the only sensible way is to proclaim an arbitrary well-ordered set is an ordinal. Before we define the concept of ordinal we must be also aware of the fact that the same set can represent different ordinals if we well-order it in different ways.

**8.6 Definition** Let a set $A$ be well-ordered by a relation $R$. Then the ordered pair $(A, R)$ is an *ordinal*. If $R$ is clear from the context then we simply say that $A$ is an *ordinal*. A *successor* of an ordinal $\alpha$ is every ordinal $\beta$ such that $\beta \setminus \{\max(\beta)\} =_o \alpha$. If $\beta$ is a successor of $\alpha$ then $\alpha$ is a *predecessor* of $\beta$. Let $\alpha$ be an ordinal different from $\emptyset$. If there exists a *predecessor* of $\alpha$ then $\alpha$ is a *successor ordinal* otherwise $\alpha$ is a *limit ordinal*.

Successors and predecessors of an ordinal are not uniquely determined. All ordinals which are isomorphic to some successor of an ordinal $\alpha$ are successors of $\alpha$. The same holds for predecessors of an ordinal.

Now we can define the derived set $A^{(\alpha)}$ for arbitrary set $A$ and arbitrary ordinal $\alpha$.

**8.7 Definition** Let $A$ be some set. Then $A^{(\emptyset)} = A$. If $\alpha$ is a successor ordinal then $A^{(\alpha)} = (A^{(p(\alpha))})'$ where $p(\alpha)$ is a predecessor of $\alpha$ and $X'$ denotes the derived set of $X$. If, however, $\alpha$ is a limit ordinal, then $A^{(\alpha)} = \{x \mid x \in A^{(\beta)}$ for every $\beta <_o \alpha\}$.

To make this definition sensible we must prove that $A^{(\alpha)}$ exists for every set $A$ and every ordinal $\alpha$. In the proof we shall apply transfinite induction on all ordinals. So we must first prove that transfinite induction can be applied in AS. When we use transfinite induction for proving that every ordinal has some property, we first prove that the ordinal Ø has this property. Then we prove that an ordinal has this property if all smaller ordinals, with regard to $<_o$, have this property.

**8.8 Lemma** *If a property* **P** *is true for at least one element of some well-ordered set A, then there exists the minimal element of A for which* **P**$(x)$ *is true.*

*Proof.* By Theorem 4.7.b there exists a set $P$, determined by $F(n + 1, P) = \{F(n, x) \mid \mathbf{P}(x) \wedge x \in A\}$ for every $n$. This set contains all elements of $A$ for which **P** is true and their accumulation points in $\mathcal{T}_a$. Therefore $P$ is a subset of the well-ordered set $A$, and as such is well-ordered as well. Hence the minimal element of $P$ exists and is an isolated point in $P$. Thus **P**(min($P$)) must be true. So min($P$) is the minimal element of $A$ for which **P**$(x)$ is true. □

**8.9 Theorem** *Transfinite induction can be applied in* AS.

*Proof.* Let some property **P** hold true for Ø, and let the fact that **P** holds true for every ordinal $\beta <_o \alpha$ imply that **P** holds true also for the ordinal $\alpha$. Then **P** holds true for every ordinal. If this were not so, then there would exist such ordinal $(\alpha, <)$ that **P**$(\alpha)$ would not be true. Every proper initial segment of $\alpha$ would be an ordinal $\alpha[x]$ determined by some $x \in \alpha$. (See the comment following the proof of Theorem 7.10.) Therefore, by Lemma 8.8, there would exist the minimal element $x'$ of $\alpha$ for which **P**$(\alpha[x'])$ would not be true. This, however, is not possible. Namely, **P**$(\alpha[\min(\alpha)])$ is true since $\alpha[\min(\alpha)] = \emptyset$. And **P**$(\alpha[y])$ is true for every $y \in \alpha$ , $y < x'$. Therefore, by assumption, **P** must be true also for $\alpha[x']$. □

**8.10 Theorem** *The set* $A^{(\alpha)}$ *exists for every set A and every ordinal* $\alpha$.

*Proof.* Let $A$ be an arbitrary set and let $A^{(\beta)}$ exist for every ordinal $\beta <_o \alpha$. Let first $\alpha$ be a successor ordinal and let $p(\alpha)$ be its predecessor. Then $A^{(p(\alpha))}$ exists by assumption and $A^{(\alpha)}$ exists by Theorem 8.1. Let now $\alpha$ be a limit ordinal and let **P**$(y, A, \alpha)$ be the predicate for an $y$ to be an element of $A^{(\alpha)}$. Then we have **P**$(y, A, \alpha) \Leftrightarrow \forall\beta( y \in A^{(\beta)} \Leftrightarrow \beta <_o \alpha)$. By Theorem 4.7.d there exists the smallest set $P$ which contains every $y$ such that **P**$(y, A, \alpha)$ is true. Let $x \in P$. Then, by Theorem 4.7.c, for every $n$ holds $F(n, x) = F(n, y_n)$ where $y_n \in A^{(\beta)}$ for $\beta <_o \alpha$. Therefore $F(n, x) = F(n, y_n) \in F(n + 1, A^{(\beta)})$ for every $n$ for every $\beta <_o \alpha$. Hence $x \in A^{(\beta)}$ for every $\beta <_o \alpha$ and so **P**$(x, A, \alpha)$ is true. Thus $A^{(\alpha)}$ exists by Theorem 4.7.e. □

Now we shall prove by transfinite induction on $\alpha$ that obtaining the set $A^{(\alpha)}$ is a legitimate operation.

**8.11 Theorem** *If* $\alpha =_o \beta$ *then* $A^{(\alpha)} = A^{(\beta)}$.

*Proof.* The theorem is true for $\alpha = \emptyset$. Let the theorem be true for every ordinal $\gamma <_o \alpha$ and let first $\alpha$ be a successor ordinal. Then $A^{(\alpha)} = (A^{(p(\alpha))})'$ where $p(\alpha)$ is a predecessor of $\alpha$. And since $\beta$ being isomorphic to $\alpha$ is a successor ordinal as well, we also have $A^{(\beta)} = (A^{(p(\beta))})'$ where $p(\beta)$ is a predecessor of $\beta$. Since $p(\alpha) =_o p(\beta)$ we have $A^{(p(\alpha))} = A^{(p(\beta))}$ by assumption. This implies

$$A^{(\alpha)} = (A^{(p(\alpha))})' = (A^{(p(\beta))})' = A^{(\beta)}.$$

Let now $\alpha$ be a limit ordinal. Then $\alpha =_o \beta$ implies that $\beta$ is a limit ordinal too, and that $\gamma <_o \alpha$ holds precisely when $\gamma <_o \beta$ holds. Hence some $x$ is an element of every $A^{(\gamma)}$ for $\gamma <_o \alpha$ precisely when this $x$ is an element of every $A^{(\gamma)}$ for $\gamma <_o \beta$. Therefore $A^{(\alpha)} = A^{(\beta)}$. So, by transfinite induction on $\alpha$, the theorem is true for any ordinal $\alpha$. □

To introduce Cantor-Bendixson rank into AS we need to prove two additional facts.

**8.12 Lemma** *Let A and B be two sets. Then* $(A \cup B)^{(\alpha)} = A^{(\alpha)} \cup B^{(\alpha)}$ *for every ordinal* $\alpha$.

*Proof.* Let $A$ and $B$ be two sets. First we shall prove that $(A \cup B)' = A' \cup B'$. Since $A, B \subseteq A \cup B$, every accumulation point in $A$ or in $B$ is an accumulation point in $A \cup B$. Thus $A' \cup B' \subseteq (A \cup B)'$. And every accumulation point in $A \cup B$ is an accumulation point either in $A$ or in $B$. Consequently we have $(A \cup B)' \subseteq A' \cup B'$. Now we shall prove that the theorem holds for arbitrary ordinal $\alpha$. The theorem certainly holds for $\alpha = \emptyset$. Let first $\alpha$ be a successor ordinal and let the theorem hold for its predecessor $p(\alpha)$. Then the conclusion we made in the beginning of this proof implies

$$(A \cup B)^{(\alpha)} = ((A \cup B)^{(p(\alpha))})' = (A^{(p(\alpha))} \cup B^{(p(\alpha))})' = (A^{(p(\alpha))})' \cup (B^{(p(\alpha))})' = A^{(\alpha)} \cup B^{(\alpha)}.$$

Let now $\alpha$ be a limit ordinal and let $A^{(\beta)} \cup B^{(\beta)} = (A \cup B)^{(\beta)}$ hold for every $\beta <_o \alpha$. If $x \in A^{(\alpha)} \cup B^{(\alpha)}$, then $x \in A^{(\beta)}$ for $\beta <_o \alpha$, or $x \in B^{(\beta)}$ for $\beta <_o \alpha$. Then $x \in A^{(\beta)} \cup B^{(\beta)}$ for $\beta <_o \alpha$ which by assumption means that $x \in (A \cup B)^{(\beta)}$ for $\beta <_o \alpha$. Then $x \in (A \cup B)^{(\alpha)}$ by definition. If, however, $x \notin A^{(\alpha)} \cup B^{(\alpha)}$ then there must exist such $\beta_A <_o \alpha$ that $x \notin A^{(\beta)}$ for $\beta >_o \beta_A$ and there exists such $\beta_B <_o \alpha$ that $x \notin B^{(\beta)}$ for $\beta >_o \beta_B$. So, if $\beta_A <_o \beta <_o \alpha$ and $\beta_B <_o \beta <_o \alpha$, then $x \notin A^{(\beta)} \cup B^{(\beta)} = (A \cup B)^{(\beta)}$. Therefore $x \notin (A \cup B)^{(\alpha)}$. Hence for every $x$ holds $x \in (A \cup B)^{(\alpha)}$ precisely when $x \in A^{(\alpha)} \cup B^{(\alpha)}$. Thus $(A \cup B)^{(\alpha)} = A^{(\alpha)} \cup B^{(\alpha)}$. □

For every closed countable set $A$ exists in ZF such countable ordinal $\alpha$ that $A^{(\alpha)} = \emptyset$ ([H&J] , p. 192, Theorem 4.10). We shall prove the analogical statement in AS. That is, for every well-ordered set $A$ exists such ordinal $\alpha$ that $A^{(\alpha)} = \emptyset$.

**8.13 Theorem** *Let A be a well-orderable set. Then* $(A^{(A)})' = \emptyset$.

*Proof.* The theorem holds for $A = \emptyset$. Let $A$ be a nonempty ordinal and let the theorem hold for every ordinal $B <_o A$. If $A$ is a successor ordinal, then $p(A) = A \setminus \{\max(A)\}$ is a predecessor of $A$ and $(p(A)^{(p(A))})' = \emptyset$ by assumption. Hence, by Lemma 8.12, we have

$$A^{(A)} = (p(A) \cup \{\max(A)\})^{(A)} = p(A)^{(A)} \cup \{\max(A)\}^{(A)} = (p(A)^{(p(A))})' \cup \emptyset = \emptyset.$$

So the theorem holds for every successor ordinal $A$. Let now $A$ be a limit ordinal. If there exists some $x \in A^{(A)}$, $x \neq \max(A)$, then $x = \max(B_0)$ for some proper initial segment $B_0$ of $A$. By definition of $A^{(A)}$ we have $x \in A^{(B)}$ for every $B <_o A$. Hence

$$x \in A^{(B)} = (B_0 \cup (A \setminus B_0))^{(B)} = B_0^{(B)} \cup (A \setminus B_0)^{(B)} \qquad (1)$$

for every $B <_o A$. As $x$ is upward isolated in $A$ we have $B_0 \cap (A \setminus B_0) = \emptyset$. Since $x = \max(B_0) \in B_0$ we have $x \notin A \setminus B_0$. Hence (1) implies that $x \in B_0^{(B)}$ for every $B <_o A$. Thus also $x \in B_0^{(B_0)}$ which contradicts the assumption that $B^{(B)} = \emptyset$ for every $B <_o A$. So, either $A^{(A)} = \emptyset$ or $A^{(A)} = \{\max(A)\}$ which implies that $(A^{(A)})' = \emptyset$. □

This theorem enables us to define Cantor-Bendixson rank of an arbitrary well-ordered set $A$.

**8.14 Definition** *Cantor-Bendixson rank* of a well-ordered set $A$ is the smallest ordinal $\alpha$ with regard to $<_o$ such that $(A^{(\alpha)})' = \emptyset$. *Cantor-Bendixson degree* of a well-ordered set $A$ is the set $A^{(\alpha)}$ where $\alpha$ is a Cantor-Bendixson rank of $A$. We shall denote Cantor-Bendixson rank and Cantor-Bendixson degree of $A$ by $r_{CB}(A)$ and $d_{CB}(A)$ respectively.

**8.15 Theorem**

*a) Cantor-Bendixson rank exists for every well-ordered set.*

*b) All Cantor-Bendixson ranks of a well-ordered set are isomorphic to each other.*

*c) Cantor-Bendixson degree of every well-ordered set exists and is a uniquely determined finite set. If* $A \neq \emptyset$ *then* $d_{CB}(A) \neq \emptyset$

*Proof.*
a) Let $\mathbf{P}(x) \Leftrightarrow x \in A \wedge \alpha = A[x\,] \cup \{x\} \wedge (A^{(\alpha)})' = \emptyset$. If $x = \max(A)$ then $\mathbf{P}(x)$ is true by Theorem 8.13. Hence, by Lemma 8.8, there exists the minimal $x \in A$ for which $\mathbf{P}(x)$ is true. This $x$ is a Cantor-Bendixson rank of $A$.

b) Let $A$ be a well-orderable set and let $\alpha$ and $\beta$ be two Cantor-Bendixson ranks of $A$. By Theorem 7.12 one of $\alpha$ and $\beta$ is isomorphic to an initial segment of the other. Let $\alpha'$ be the initial segment of $\beta$ which is isomorphic to $\alpha$. Then $A^{(\alpha')} = A^{(\alpha)} \neq \emptyset$ by Theorem 8.11, and $A^{(\gamma)} = \emptyset$ for every $\gamma >_o \alpha'$. Since $\alpha' \leq_o \beta$ and since $A^{(\beta)} \neq \emptyset$ this is possible only if $\alpha' = \beta$ and hence $\alpha =_o \beta$.

c) If $\alpha$ and $\beta$ are two Cantor-Bendixson ranks of a well-ordered set $A$ then $\alpha =_o \beta$ by b). Theorem 8.11 implies $A^{(\alpha)} = A^{(\beta)} = d_{CB}(A)$. Hence $d_{CB}(A)$ is uniquely determined and exists since Cantor-Bendixson rank of $A$ exists by a). The set $d_{CB}(A)$ is finite. If it were infinite it would have at least one accumulation point. Then $\emptyset \neq (d_{CB}(A))' = (A^{(\alpha)})'$ which is not true as $\alpha =_o r_{CB}(A)$ and so $(A^{(\alpha)})' = \emptyset$. If $A \neq \emptyset$ and $\alpha = r_{CB}(A)$ then $d_{CB}(A) = A^{(\alpha)} \neq \emptyset$ by Definition 8.14. □

To determine the cardinality of sets by Cantor-Bendixson rank and Cantor-Bendixson degree we also need to state some additional lemmas.

**8.16 Lemma** *Let f be an injection from a set A to a set B and let $A_i \subseteq A$ for every $i \in J$. Then $x \in A_i$ for every $i \in J$ precisely when $f(x) \in f(A_i)$ for every $i \in J$.*

This lemma is proved in the same way as in ZF. Therefore we omit its proof.

**8.17 Lemma** *If f is a bijection between A and B, then the restriction of f to $A^{(\alpha)}$ is a bijection between $A^{(\alpha)}$ and $B^{(\alpha)}$ for every ordinal $\alpha$.*

*Proof.* We shall use transfinite induction. The lemma obviously holds for $\alpha = \emptyset$. Let the lemma hold for every $\beta <_o \alpha$. Let first $\alpha$ have its predecessor $p(\alpha)$. The restriction of $f$ to $A^{(p(\alpha))}$ is by assumption a bijection from $A^{(p(\alpha))}$ to $B^{(p(\alpha))}$. By Theorem 5.10 and Lemma 5.11, the restriction of this bijection to $A^{(\alpha)}$ is a bijection between $A^{(\alpha)}$ and $B^{(\alpha)}$. For this reason the lemma holds for $\alpha$ as well. Let now $\alpha$ be a limit ordinal. As $f(A^{(\beta)}) = B^{(\beta)}$ for every $\beta <_o \alpha$, Lemma 8.16 implies that some $x$ is an element of every $A^{(\beta)}$ for $\beta <_o \alpha$ precisely when $f(x)$ is an element of every $B^{(\beta)}$ for $\beta <_o \alpha$. This fact and Definition 8.7 imply that the restriction of $f$ to $A^{(\alpha)}$ is a bijection from $A^{(\alpha)}$ to $B^{(\alpha)}$. Thus we have proved the lemma completely. □

**8.18 Lemma** *If $A \subseteq B$, then $A^{(\alpha)} \subseteq B^{(\alpha)}$ for every ordinal $\alpha$.*

*Proof.* Again we shall use transfinite induction. The lemma obviously holds for $\alpha = \emptyset$. Let the lemma hold for every $\beta <_o \alpha$. So, if $p(\alpha)$ is a predecessor of $\alpha$, then $A^{(p(\alpha))} \subseteq B^{(p(\alpha))}$. If $x$ is an accumulation point in $A^{(p(\alpha))}$, then $x$ is an accumulation point in $B^{(p(\alpha))}$. Therefore $A^{(\alpha)} \subseteq B^{(\alpha)}$ and the lemma holds for $\alpha$. Let now $\alpha$ be a limit ordinal. By Definition 8.7 we have $A^{(\alpha)} \subseteq A^{(\beta)}$ for every $\beta <_o \alpha$, and by assumption we have $A^{(\beta)} \subseteq B^{(\beta)}$ for every $\beta <_o \alpha$. This implies $A^{(\alpha)} \subseteq B^{(\beta)}$ for every $\beta <_o \alpha$. Hence, by Definition 8.7, we have $A^{(\alpha)} \subseteq B^{(\alpha)}$. Thus the lemma is proved completely. □

**8.19 Lemma** *A well-ordered nonempty set A can be well-ordered in such a way that* $\max(d_{CB}(A)) = \max(A)$.

*Proof.* Let $(A, <)$ be a well-ordered nonempty set and let $m = \max(d_{CB}(A))$. Let's define $\leq$ and $>$ as usual. Then the sets $A_1 = \{x \mid x \in A \wedge x \leq m\}$ and $A_2 = \{x \mid x \in A \wedge x > m\}$ exist. Namely, $A_1$ and $A_2$ defined in this way contain all their accumulation points. The point $m$ and has its successor $m'$ in $(A, <)$. Therefore $m = \max(A_1)$ and $m' = \min(A_2)$. Now we define the relation $<_m$ which is inside $A_1$ and $A_2$ the same as $<$. If, however, $a_1 \in A_1$ and $a_2 \in A_2$ then let $a_2 <_m a_1$. In other words, we move $A_2$ in front of $A_1$. Hence $m = \max(d_{CB}(A)) = \max(A)$ with regard to $<_m$. □

Now, at last, we can classify well-orderable sets by their Cantor-Bendixson rank and degree.

**8.20 Theorem** *Let A and B be well-orderable sets. Then $A =_c B$ precisely when $r_{CB}(A) =_o r_{CB}(B)$ and $d_{CB}(A) =_o d_{CB}(B)$.*

*Proof.* If $A =_c B$ then there exists an injection $f$ from $A$ into $B$. Therefore there exists a bijection between $A$ and $f(A)$. So the restriction of $f$ to $A^{(\alpha)}$ is by Lemma 8.17 a bijection from $A^{(\alpha)}$ to $(f(A))^{(\alpha)}$ for every ordinal $\alpha$. Since $f(A) \subseteq B$, this restriction is by Lemma 8.18 an injection from $A^{(\alpha)}$ into $B^{(\alpha)}$ for every ordinal $\alpha$. This implies $r_{CB}(A) \leq_o r_{CB}(B)$. Since there also exists an injection from $B$ into $A$ we conclude in the same way that $r_{CB}(B) \leq_o r_{CB}(A)$. Hence $r_{CB}(A) =_o r_{CB}(B)$. Therefore, if $\alpha =_o r_{CB}(A)$ then $A^{(\alpha)} = d_{CB}(A)$ and $B^{(\alpha)} = d_{CB}(B)$. The sets $d_{CB}(A)$ and $d_{CB}(B)$ are by Theorem 8.15.c finite, and from each of them exists an injection into the other one. This is possible only if $d_{CB}(A) =_o d_{CB}(B)$.

Let now $r_{CB}(A) =_o r_{CB}(B)$ and $d_{CB}(A) =_o d_{CB}(B)$. By Lemma 8.19 there exists a well-ordering $<_A$ of $A$ such that $\max(A) = \max(d_{CB}(A))$ with regard to $<_A$. There also exists such well-ordering $<_B$ of $B$. At least one of $A$ and $B$ is isomorphic to an initial segment of the other by Theorem 7.12. Without loss of generality we can assume that $A$ is isomorphic to an initial segment $S$ of $B$. Hence there is a bijection between $A^{(\alpha)}$ and $S^{(\alpha)}$ for every ordinal $\alpha$ by Lemma 8.17. Let $\alpha =_o r_{CB}(A) =_o r_{CB}(B)$. The sets $S^{(\alpha)}$ and $B^{(\alpha)}$ are finite having $d_{CB}(A)$ and $d_{CB}(B)$ elements respectively. We have $d_{CB}(A) =_o S^{(\alpha)} \subseteq B^{(\alpha)} = d_{CB}(B)$ for every $\alpha$ by Lemma 8.18. If $S$ were a proper initial segment of $B$, then $\max(B) \notin S$ and therefore $\max(B) \notin S^{(\alpha)} = d_{CB}(S)$. On the other hand, $\max(B) = \max(d_{CB}(B)) \in d_{CB}(B)$ by Lemma 8.19. This implies that $d_{CB}(S)$ is a proper subset of the finite set $d_{CB}(B)$. And since $d_{CB}(S) =_o d_{CB}(A)$ we can conclude that $d_{CB}(A) <_o d_{CB}(B)$. But this is in contradiction with the assumption that $d_{CB}(A) =_o d_{CB}(B)$. Hence $A =_o S = B$ which implies $A =_c B$. □

**8.21 Theorem** *Let A and B be well-orderable sets. Then $A <_c B$ precisely when $r_{CB}(A) <_o r_{CB}(B)$ or $r_{CB}(A) =_o r_{CB}(B)$ and $d_{CB}(A) <_o d_{CB}(B)$.*

*Proof.* If $A <_c B$ is true then there exists an injection $f$ from $A$ into $B$. Hence $A =_c f(A)$. Then Theorem 8.20 implies $r_{CB}(A) =_o r_{CB}(f(A))$ in $d_{CB}(A) =_o d_{CB}(f(A))$. Since $f(A) \subseteq B$ we have $r_{CB}(A) \leq_o r_{CB}(B)$. If $r_{CB}(A) =_o r_{CB}(B)$ then $d_{CB}(A) \leq_o d_{CB}(B)$ since $f(A) \subseteq B$. But $d_{CB}(A) =_o d_{CB}(B)$ cannot hold in this case otherwise $A =_c B$ by Theorem 8.20 which, however, is not true. If we resume then $A <_c B$ implies either $r_{CB}(A) <_o r_{CB}(B)$ or $r_{CB}(A) =_o r_{CB}(B)$ and $d_{CB}(A) <_o d_{CB}(B)$.

If, however, $A <_c B$ does not hold then $B \leq_c A$ by Theorem 7.15. If $B <_c A$ then the previous argument implies either $r_{CB}(B) <_o r_{CB}(A)$ or $r_{CB}(B) =_o r_{CB}(A)$ and $d_{CB}(B) <_o d_{CB}(A)$. And if $B =_c A$ then Theorem 8.20 implies $r_{CB}(B) =_o r_{CB}(A)$ and $d_{CB}(B) =_o d_{CB}(A)$. Therefore, if $A <_c B$ is not true then we have either $r_{CB}(B) <_o r_{CB}(A)$ or $r_{CB}(B) =_o r_{CB}(A)$ and $d_{CB}(B) \leq_o d_{CB}(A)$. In other words, if $A <_c B$ is not true then the following statement is not true: $r_{CB}(A) <_o r_{CB}(B)$ or $r_{CB}(A) =_o r_{CB}(B)$ and $d_{CB}(A) <_o d_{CB}(B)$.

If we resume, $A <_c B$ is true precisely when the statement $r_{CB}(A) <_o r_{CB}(B)$ or $r_{CB}(A) =_o r_{CB}(B)$ and $d_{CB}(A) <_o d_{CB}(B)$ is true. □

# 9. The Sum of Sets and its Cardinality

We have seen that the ordinal arithmetic of von Neumann ordinals is trivial in AS. However, we can define the not so trivial arithmetic based on ordinals defined in the previous section. This arithmetic is analogical to the cardinal arithmetic in ZF and we shall refer to it in the sequel by this name, although this is not an arithmetic of cardinals but an arithmetic of arbitrary well-ordered sets. As we shall see cardinal arithmetic in AS is pretty different from cardinal arithmetic in ZF.

As we shall prove in section 6, several isomorphic cardinals belong to the same set. Hence we shall not define the sum and the product of cardinals as this definition complicates proofs. Instead we shall define the sum and the product of arbitrary sets. The sum of a family $\{A_j \mid j \in J\}$ is in ZF defined as the union $\cup_{j\in J} A_j$ of mutually disjoint sets $A_j$. We shall not restrict ourselves only to the families of mutually disjoint sets but shall replace every set $A_j$ in a family with the set $\{j\} \times A_j$. This set is in bijection with $A_j$ and is mutually disjoint with every set $\{k\} \times A_k$ where $k \in J$ and $j \neq k$. Another reason for such approach is that the sum of a well-ordered family of well-ordered sets is a well-ordered set since $j <_J k$, where $<_J$ is a well-ordering of $J$, implies that every element of $\{j\} \times A_j$ can be treated as being smaller than every element of $\{k\} \times A_k$, which, in general, does not hold for the elements of $A_j$ and $A_k$ even if $A_j$ and $A_k$ are disjoint.

**9.1 Definition** Let $\mathcal{F}$ be a family of sets and let $J$ be an arbitrary set, called the *set if indexes*. Then every function $I : J \rightarrow \mathcal{F}$ is an *indexed system of sets*. Thus $I = \{(j, A_j) \mid j \in J \wedge A_j \in \mathcal{F}\}$. The *sum* of $I$ is the set $\cup\{\{j\} \times A_j \mid (j, A_j) \in I\}$. We denote this sum $\cup_{j\in J}\{j\} \times A_j$ or $\sum_{j\in J} A_j$. If $J = \{1, 2, ... n\}$ then we write it $A_1 + A_2 + \ldots + A_n$.

**9.2 Remark** The sum of sets is here defined for an arbitrary indexed system of sets, while the sum of sets is in ZF defined only for such indexed system of sets where the sets $A_j$ are mutually disjoint. If for some indexed system of sets both sums can be defined, then these two sums have the same cardinality. Namely, in this case $\cup_{j\in J} A_j =_c \cup_{j\in J}\{j\} \times A_j$ because the map $a \rightarrow (j, a)$, where $a$ can be any element of $A_j$ for any $j \in J$, is a bijection from $\cup_{j\in J} A_j$ to $\cup_{j\in J}\{j\} \times A_j$.

Of course, Definition 9.1 presumes that the sum of $I$ exists for every indexed system of sets $I$. Therefore we must first prove that this really is the case.

**9.3 Lemma** *For every $n \geq 1$ holds $F(n + 2, A \times B) = F(n, A) \times F(n, B)$.*

*Proof.* By Theorem 4.17 we have for every $n \geq 1$:

$$F(n + 3, A \times B) = \{(F(n, a), F(n, b)) \mid (a, b) \in A \times B\} =$$

$$\{F(n, a) \mid a \in A\} \times \{F(n, b) \mid b \in B\} = F(n + 1, A) \times F(n + 1, B).$$

The replacement of $n + 1$ by $n$ where, of course, $n \geq 1$ proves the theorem. □

**9.4 Theorem** *The sum of every indexed system of sets $I = \{(j, A_j) \mid j \in J \wedge A_j \in \mathcal{F}\}$ exists.*

*Proof.* The sum of $I$ is the union of the set $\{\{j\} \times A_j \mid (j, A_j) \in I\}$. Hence we must first prove that this set exists. The corresponding predicate for $x$ to be an element of this set is

$$\mathbf{P}(x) \Leftrightarrow \exists j\, \exists A_j (j \in J \wedge (j, A_j) \in I \wedge x = \{j\} \times A_j).$$

This predicate exists by Theorem 4.5. By Theorem 4.7.d there exists the smallest set $P$ which contains all $x$ determined by $\mathbf{P}$. That is, $F(n + 1, P) = \{F(n, x) \mid \mathbf{P}(x)\}$. If $x \in P$ then, by Theorem 4.7.c, for every $n$ there exists such set $\{j(n)\} \times A_{j(n)}$, where $(j(n), A_{j(n)}) \in I$, that $F(n, x) = F(n, \{j(n)\} \times A_{j(n)})$. Among the pairs $(j(n), A_{j(n)})$ there are either only finitely many different, or they have at least one accumulation point in $\mathcal{T}_a$, and this point is, of course, an element of $I$. In each case there exists such pair $(k, A_k) \in I$ that for infinitely many $n$ holds $F(n, (k, A_k)) = F(n, (j(n), A_{j(n)}))$. Then, by Theorem 4.16, we have $(F(n - 2, k), F(n - 2, A_k)) = (F(n - 2, j(n)), F(n - 2, A_{j(n)}))$ for infinitely many $n$. Hence $F(n - 2, k) = F(n - 2, j(n))$ and $F(n - 2, A_k) = F(n - 2, A_{j(n)})$ holds for infinitely many $n$. This fact and Lemma 9.3 imply that

$$F(n, \{k\} \times A_k) = \{F(n-3, k)\} \times F(n-2, A_k) = \{F(n-3, j(n))\} \times F(n-2, A_{j(n)}) = F(n, x)$$

holds for infinitely many $n$. Now, if two sets have the same $n$th approximation then these two sets have the same $m$th approximation for every $m < n$. Therefore, since $F(n, \{k\} \times A_k) = F(n, x)$ holds for infinitely many $n$, it must hold for every $n$. Hence $x = \{k\} \times A_k$. Since $\mathbf{P}(\{k\} \times A_k)$ is true, $\mathbf{P}(x)$ is true as well. Thus the set $\{\{j\} \times A_j \mid (j, A_j) \in I\}$ exists by Theorem 4.7.e. And the sum of $I$ being the union of this set exists by Theorem 4.9. □

It is not difficult to construct bijections showing that the sum of sets is associative and commutative with regard to $=_c$. That is $(A_1 + A_2) + A_3 =_c A_1 + (A_2 + A_3)$ and $A_1 + A_2 =_c A_2 + A_1$. And it is not difficult to see that the sum of two finite sets $A_1$ and $A_2$ corresponds to the sum of corresponding natural numbers, as described in Definition 2.3, and to the sum of corresponding finite N-ordinals in the ordinal arithmetic.

Let's show that sum of sets is a legitimate operation, which means that the cardinality of a sum is determined by the cardinalities of its summands.

**9.5 Theorem** *Let* $\{(j, A_j) \mid j \in J, A_j \in \mathcal{F}\}$ *and* $\{(j, B_j) \mid j \in J, B_j \in \mathcal{G}\}$ *be two indexed systems of sets where* $A_j =_c B_j$ *for every* $j \in J$*. Then* $\sum_{j\in J} A_j =_c \sum_{j\in J} B_j$*.*

*Proof.* Since $A_j =_c B_j$ for every $j \in J$ there exists an injection $f_j : A_j \to B_j$ for every $j \in J$. First we shall prove that for every $j \in J$ there exists the set of all ordered pairs $\{((j, x), (j, y)) \mid (x, y) \in f_j\}$. Let $P_j$ be the smallest set which contains all these pairs. The set $P_j$ exists by Theorem 4.7.d. Let $z \in P_j$. Then $F(n, z) = F(n, ((j, x_n), (j, y_n)))$ for every $n$ where $(x_n, y_n) \in f_j$. Among the pairs $(x_n, y_n)$ there are either only finitely many different, or they have at least one accumulation point in $\mathcal{T}_a$, and this point is, of course, an element of $f_j$. In each case there exists such pair $(x, y) \in f_j$ that for infinitely many $n$ holds $F(n, z) = F(n, ((j, x), (j, y)))$. This implies that $z = ((j, x), (j, y)) \in P_j$ which, by Theorem 4.7.e, means that the elements of $P_j$ are only the sets of the form $((j, x), (j, y))$ where $(x, y) \in f_j$. It is not difficult to see that $P_j$ is an injection from $\{j\} \times A_j$ into $\{j\} \times B_j$. There also exists the set of all these $P_j$. That is, there exists the set $\{P_j \mid j \in J\}$. This can be proved in a similar way as the existence of the set $\{\{j\} \times A_j \mid j \in J \wedge A_j \in \mathcal{F}\}$ (See the proof of Theorem 9.4). Thus there exists $\cup\{P_j \mid j \in J\}$. It is not difficult to see that this set is an injection from $\cup\{\{j\} \times A_j \mid j \in J\}$ into $\cup\{\{j\} \times B_j \mid j \in J\}$. This implies $\sum_{j\in J} A_j \leq_c \sum_{j\in J} B_j$. As there also exist injections $g_j : B_j \to A_j$ for every $j \in J$ we prove in the same way as before that $\sum_{j\in J} B_j \leq_c \sum_{j\in J} A_j$. Hence $\sum_{j\in J} A_j =_c \sum_{j\in J} B_j$. □

An important question is when the sum of sets is well-orderable.

**9.6 Theorem** *The sum of an indexed system of sets* $\{(j, A_j) \mid j \in J, A_j \in \mathcal{F}\}$ *is well-orderable precisely when* $J$ *is well-ordered and when every* $A_j$ *is well-ordered.*

*Proof.* Let every $A_j$ and $J$ be well-ordered sets. First we arrange the nodes of $\mathbf{T}_J$ in such a way that the order of these nodes determines a well-ordering < of $J$. This can be done because $J$ is well-orderable (See the proof of Theorem 7.8). Since every $j \in J$ is an upward isolated point in the well-ordering < of $J$, there exists such $n$ that the node representing $F(n, j)$ is in the arrangement of $\mathbf{T}_J$ more to the left than the node representing $F(n, j')$ if $j < j' \in J$. In other words, no infinite path in $\mathbf{T}_J$ separates from the path for $j$ to the right in the $(n+i)$th node of this path for any $i$. As also every $A_j$ is well-ordered, the same holds for the paths of the elements of every $A_j$ in a well-ordering arrangement $<_j$ of the developing tree of $A_j$. The sum $\cup_{j\in J} \{j\} \times A_j$ consists of ordered pairs $(j, a_{j,i})$ where $j \in J$ and $a_{j,i} \in A_j$. We well-order this sum as follows. The $n$th approximation of $(j, a_{j,i})$ is an ordered pair of the $(n-2)$th approximations of $j$ and $a_{j,i}$. Therefore, by the previous argument, there must exist such $n$ that the $(n+i)$th approximation of $(j, a_{j,i})$ is smaller than the $(n+i)$th approximation of $(j', a_{j,i}')$ for any $i$ if $j < j'$, or if $j = j'$ and $a_{j,i} <_j a_{j,i}'$. So every pair $(j, a_{j,i})$ is an upward isolated point in the described linear ordering of the sum $\sum_{j\in J} A_j$. Hence, if $J$ is a well-orderable set and if all $A_j$, $j \in J$, are well-orderable sets, then $\cup_{j\in J} \{ j\} \times A_j$ is a well-orderable set.

The set $\cup_{j\in J} \{ j\} \times A_j$ is not well-orderable in any other case. If some $A_j$ is not well-orderable then there exists a perfect subset $P$ of $A_j$, and it is not difficult to see that $\{ j\} \times P$ is a perfect set as well. Since $\{ j\} \times P$ is a subset of the sum $\cup_{j\in J} \{ j\} \times A_j$, this sum cannot be well-ordered. If, however, $J$ is not well-orderable then $\cup_{j\in J} \{ j\} \times A_j$ is again not well-orderable which we prove as follows. In this case there must exist some perfect subset $P$ of $J$ and there exists the injection $f : J \to \cup_{j\in J} \{ j\} \times A_j$ which maps every $j \in J$ to the ordered pair $( j, \min(A_j))$ where $\min(A_j)$ is the minimal element of $A_j$ with regard to the lexicographical ordering of $A_j$ (See Definition 7.2). The function $f$ exists since there exists its approximation $F(n, f)$ for every $n$. Namely, for every $n \geq 4$ we have

$$F(n, f) = F(n, ( j, ( j, \min(A_j)))) = (F(n-2, j), F(n-2, ( j, \min(A_j)))) =$$

$$(F(n-2, j), (F(n-4, j), F(n-4, \min(A_j)))).$$

The HF set from the last row exists for every $n \geq 4$, for every $j$ and for every $A_j$. Therefore $f$ exists and it is not difficult to see that $f$ is an injection. Hence the $f$ image of the perfect set $P$, which is a subset of $J$, is a perfect subset of $\cup_{j\in J} \{ j\} \times A_j$ since $f$ being an injection maps accumulation points in $J$ to accumulation points. Therefore $\cup_{j\in J} \{ j\} \times A_j$ is not well-orderable. □

If $A_1$ and $A_2$ are infinite sets in ZFC and $A_1 \leq_c A_2$ then $A_1 + A_2 =_c A_2$. This is true in AS only to a limited extent as Theorem 9.8 will show. First we prove one lemma.

**9.7 Lemma** *Let $A_1 =_c A_2$. Then $A_1^{(\alpha)} =_c A_2^{(\alpha)}$ for every ordinal $\alpha$.*

*Proof.* As $A_1 =_c A_2$ there exists an injection $f_1 : A_1 \to A_2$. Hence, by Lemma 8.17, there exists a bijection between $A_1^{(\alpha)}$ and $(f_1(A_1))^{(\alpha)}$ for every ordinal $\alpha$. Since $f_1(A_1) \subseteq A_2$ we have $(f_1(A_1))^{(\alpha)} \subseteq A_2^{(\alpha)}$ by Lemma 8.18 for every ordinal $\alpha$. This two facts imply that for every ordinal $\alpha$ there exists an injection from $A_1^{(\alpha)}$ into $A_2^{(\alpha)}$. Since $A_1 =_c A_2$ there exists also an injection $f_2 : A_2 \to A_1$. So we conclude in the same way that there exists an injection from $A_2^{(\alpha)}$ into $A_1^{(\alpha)}$ for every ordinal $\alpha$. Thus $A_1^{(\alpha)} =_c A_2^{(\alpha)}$ for every ordinal $\alpha$. □

**9.8 Theorem** *Let $A_1 \leq_c A_2$ where $A_1$ and $A_2$ are well-ordered. If $r_{CB}(A_1) <_o r_{CB}(A_2)$, then $A_1 + A_2 =_c A_2$. If, however, $r_{CB}(A_1) =_o r_{CB}(A_2)$, then $A_1 + A_2 >_c A_2$.*

*Proof.* We shall determine the cardinality of $A_1 + A_2$ by Cantor-Bendixson rank and degree. Let now $A_1 \leq_c A_2$ where $A_1$ and $A_2$ are well-ordered sets. Then, by Lemma 8.12, we have

$$(A_1 + A_2)^{(\alpha)} = (\{1\} \times A_1 \cup \{2\} \times A_2)^{(\alpha)} = (\{1\} \times A_1)^{(\alpha)} \cup (\{2\} \times A_2)^{(\alpha)} \qquad (1)$$

For every sets $A$ and $a$ and for every ordinal $\alpha$ holds $(\{a\} \times A)^{(\alpha)} =_c A^{(\alpha)}$ by Lemma 9.7. And since certainly $A^{(\alpha)} =_c \{a\} \times A^{(\alpha)}$ is true for every $A$ , $a$, and $\alpha$, we have

$$(\{1\} \times A_1)^{(\alpha)} \cup (\{2\} \times A_2)^{(\alpha)} =_c \{1\} \times A_1^{(\alpha)} \cup \{2\} \times A_2^{(\alpha)} = A_1^{(\alpha)} + A_2^{(\alpha)} \qquad (2)$$

(1) and (2) together give

$$(A_1 + A_2)^{(\alpha)} =_c (A_1)^{(\alpha)} + (A_2)^{(\alpha)} \qquad (3)$$

for every $\alpha$. Since $A_1 \leq_c A_2$ there exists an injection $f : A_1 \to A_2$ which is a bijection between $A_1$ and $f(A_1)$. And since $f(A_1) \subseteq A_2$ we have $A_1^{(\alpha)} =_c f(A_1)^{(\alpha)} \subseteq A_2^{(\alpha)}$ for every ordinal $\alpha$ by Lemmas 9.7 and 8.18. Hence $A_1^{(\alpha)} \leq_c A_2^{(\alpha)}$ for every ordinal $\alpha$ which implies $r_{CB}(A_1) \leq_o r_{CB}(A_2)$. Let first $r_{CB}(A_1) <_o r_{CB}(A_2) = r$. Then (3) implies

$$(A_1 + A_2)^{(r)} =_c A_1^{(r)} + A_2^{(r)} = \emptyset + A_2^{(r)} = A_2^{(r)}.$$

This implies $r = r_{CB}(A_2) =_o r_{CB}(A_1 + A_2)$ and $d_{CB}(A_2) = A_2^{(r)} =_c (A_1 + A_2)^{(r)} = d_{CB}(A_1 + A_2)$. So, by Theorem 8.20, we have $A_1 + A_2 =_c A_2$. Let now $r_{CB}(A_1) =_o r_{CB}(A_2) = r$. Then $A_1^{(r)}$ and $A_2^{(r)}$ are finite sets $d_{CB}(A_1)$ and $d_{CB}(A_2)$ respectively and we have

$$d_{CB}(A_1) + d_{CB}(A_2) = A_1^{(r)} + A_2^{(r)} =_c (A_1 + A_2)^{(r)} = d_{CB}(A_1 + A_2)$$

Since $d_{CB}(A_1)$ is not the empty set, we have $d_{CB}(A_1 + A_2) >_o d_{CB}(A_2)$. Hence $A_1 + A_2 >_c A_2$ by Theorem 8.21. □

# 10. The Product of Sets and its Cardinality

In ZF the product of an indexed system $I = \{(j, A_j) \mid j \in J, Aj \in \mathcal{F}\}$ is the set of all those functions $f$ whose domain is $J$ and for which holds $f(j) \in A_j$ for every $j \in J$. In AS this set does in general not exist. Let, for example, $I = \{(j, A_j) \mid j \in \boldsymbol{N}, A_j = \{0, 1\}\}$. Let $k$ be a natural number and let $R_k$ be the set of ordered pairs $(j, m) \in \boldsymbol{N} \times \{0, 1\}$ for which the following holds: if $j = 2i$ for some $i < k$ then $(j, 0) \in R_k$ otherwise $(j, 1) \in R_k$. The set $R_k$ is a function $f$ described above since it maps $j \in \boldsymbol{N}, j < \infty$, uniquely to 0 or to 1 and it maps $\infty$ to 1. Let now $R$ be the set of all ordered pairs of the form $(2n, 0)$ and $(2n + 1, 1)$ where $n \in \boldsymbol{N}$. Then $R$ is in ZF a function $f$ described above and hence $R$ is in ZF an element of the product of $I$. However, $R$ is not a function in AS since here both $(\infty, 0)$ and $(\infty, 1)$ are elements or $R$. Nevertheless $R$ is an element of the smallest set in AS which contains all functions $f$ described above. Namely, it is not difficult to see that $R$ is an accumulation point of the functions $R_k$ in $\mathcal{T}_a$. According to these findings and to Theorem 4.7 we shall state the following definition.

**10.1 Definition** A *choice function on* an indexed system of sets $I = \{(j, A_j) \mid j \in J, Aj \in \mathcal{F}\}$ is a function $f$ with domain $J$ such that for every $j \in J$ holds $f(j) \in A_j$. The *product* of $I$ is the smallest set to which belong all choice functions on $I$. We denote the product of $I$ by $\prod_{j\in J} A_j$. If $J = \{1, 2, ..., n\}$ then we write the product of $I$ also as $A_1 \cdot A_2 \cdot \ldots \cdot A_n$.

The first thing we should do is to prove that the product of sets is legitimate. That is, we should prove that the cardinality of the product is determined by the cardinalities of its factors. However the proof of this does not appear to be a simple one. It seems that the easiest way to carry it out is to determine first when the product of sets is well-orderable. We shall see that this is true precisely when only finitely many of its factors are greater than 1. In this case the product of sets can be interpreted as Cartesian product which, as we shall see, is well-orderable precisely when all of its factors are well-orderable. So let us first introduce Cartesian product of a finite number of factors. This one can be defined in AS in the same way as in ZF.

**10.2 Definition** Every $x \in U$ is an ordered 1-*tuple*. If $x$ is an ordered $n$-tuple and $y \in U$, then the ordered pair $(x, y)$ is an ordered $(n + 1)$-*tuple*. Every set is *the Cartesian product of one set*. If a set $P$ is the Cartesian product of $n$ sets and $A$ is an arbitrary set, then the set $P \times A$ is the *Cartesian product of* $n + 1$ *sets*. The Cartesian product of the sets $A_1, A_2, \ldots, A_n$ is written as $A_1 \times A_2 \times \ldots \times A_n$.

The Cartesian product $A_1 \times A_2 \times \ldots \times A_n$ exists for arbitrary sets $A_1, A_2, \ldots, A_n$. This can be proved by induction on $n$ since the Cartesian product of two sets exists by Theorem 4.17. Cartesian product is associative and distributive with regard to the sum of sets defined in the previous section. That is, $(A \times B) \times C = A \times (B \times C)$ and $A \times (B + C) = A \times B + A \times C$. These facts can be proved in the same way as in ZF. The map $(a, b) \to (b, a)$, where $a \in A$ and $b \in B$, is a bijection between $A \times B$ and $B \times A$. This fact can also be proved in the same way as in ZF. Therefore $A \times B =_c B \times A$ although $A \times B \neq B \times A$ in general.

Let's now show that the product of finitely many sets, as defined in 10.1, can be interpreted as Cartesian product.

**10.3 Theorem** *The set* $A_1 \cdot A_2 \cdot \ldots \cdot A_n$ *is in bijection with the set* $A_1 \times A_2 \times \ldots \times A_n$.

*Proof.* First we shall prove that every element of the set $A_1 \cdot A_2$ is a map $f$ with domain $\{1, 2\}$ such that $f(1) \in A_1$ and $f(2) \in A_2$. Namely, let $z \in A_1 \cdot A_2$. If $z$ is not one of the described functions $f$, then $z$ is an accumulation point of such functions in $\mathcal{T}_a$ by Theorem 4.7.c. Hence $F(n, z) = F(n, f_n)$ holds for every $n$ where $f_n = \{(1, a_{n,1}), (2, a_{n,2})\}$, $a_{n,1} \in A_1$ and $a_{n,2} \in A_2$. Therefore

$$F(n + 3, z) = F(n + 3, f_n) = F(n + 3, \{(1, a_{n,1}), (2, a_{n,2})\}) = \{F(n + 2, (1, a_{n,1})), F(n + 2, (2, a_{n,2}))\} =$$

$$\{(F(n, 1), F(n, a_{n,1})), (F(n, 2), F(n, a_{n,2}))\} = \{(1, F(n, a_{n,1})), (2, F(n, a_{n,2}))\}$$

holds for every large enough $n$ such that $F(n, 2) = 2$. So $F(0, a_{0,1})$, $F(1, a_{1,1})$, $F(2, a_{2,1})$, ... is the developing sequence of some $a_1 \in A_1$. Namely, the previous argument and the equations, together with the fact that $F(n, (F(n + 1, z)) = F(n, z)$ holds for every $n$, imply that $F(n, F(n + 1, a_{n+1,1}) = F(n, a_{n,1})$ holds for every $n$. In a similar way we can conclude that $F(0, a_{0,2})$, $F(1, a_{1,2})$, $F(2, a_{2,2})$, ... is the developing sequence of some $a_2 \in A_2$. For this reason is $F(3, f_0)$, $F(4, f_1)$, $F(5, f_2)$, ... the developing sequence of the set $\{(1, a_1),(2, a_2)\}$. As this developing sequence is also the developing sequence of $z$ we have $z = \{(1, a_1),(2, a_2)\}$. Therefore $z$ is a function which maps 1 to some element of $A_1$ and 2 to some element of $A_2$.

Let's now prove that the map $(a_1, a_2) \to \{(1, a_1), (2, a_2)\}$, where $a_1$ and $a_2$ are arbitrary elements of $A_1$ and $A_2$ respectively, is a bijection between $A_1 \times A_2$ and $A_1 \cdot A_2$. This map is total on $A_1 \times A_2$. It is injective, since $(a_1, a_2) \neq (a_1', a_2')$ implies $\{(1, a_1), (2, a_2)\} \neq \{(1, a_1'), (2, a_2')\}$. It is surjective, since we have just proved that every element of $A_1 \cdot A_2$ is of the form $\{(1, a_1), (2, a_2)\}$, where $a_1 \in A_1$ and $a_2 \in A_2$, and so every element of $A_1 \cdot A_2$ has the inverse picture in $A_1 \times A_2$. The map is also a function since the same pair $(a_1, a_2)$ from $A_1 \times A_2$ cannot be mapped to different elements of $A_1 \cdot A_2$.

From the above arguments we prove by induction on $n$ that for every $n$ the product $A_1 \cdot A_2 \cdot \ldots \cdot A_n$ is the set of all functions which map $j$ to some element of $A_j$ for $1 \leq j \leq n$ and consequently we prove that there exists a bijection from $A_1 \times A_2 \times \ldots \times A_n$ to $A_1 \cdot A_2 \cdot \ldots \cdot A_n$. □

This theorem and the properties of Cartesian product imply $(A \cdot B) \cdot C =_c A \cdot (B \cdot C)$ and $A \cdot B =_c B \cdot A$ for arbitrary sets $A$, $B$ and $C$. Also $A \cdot (B + C) = A \cdot B + A \cdot C$. Namely, because of the properties of Cartesian product, Theorem 9.5, and Theorem 10.3 we have

$$A \cdot (B + \mathrm{C}) =_c A \times (B + \mathrm{C}) = A \times B + A \times C =_c A \cdot B + A \cdot C$$

It is also not difficult to see that finite product of finite sets matches with the product of the corresponding finite N-ordinals in the ordinal arithmetic.

In a cardinal arithmetic we are interested only in the cardinality of sets. As we have proved, finite product is in bijection with the corresponding Cartesian product which is easier to imagine and to handle. Hence, in order to determine the cardinality of a finite product we shall from now on investigate only Cartesian products instead of finite products. First we shall examine when the Cartesian product is well-orderable..

**10.4 Theorem** *If the sets A and B are well-orderable then $A \times B$ can be well-ordered.*

*Proof.* As $A$ is well-orderable, Theorem 7.8 implies that there exists an arrangement of $\mathbf{T}_A$ which well-orders $A$. Let $v_1$ and $v_2$ be different nodes on the same level of $\mathbf{T}_A$. If $v_1$ is to the left of $v_2$ in this arrangement we shall write $v_1 <_A v_2$. It is not difficult to see that $<_A$ is a linear ordering of the nodes on any level of $\mathbf{T}_A$. Such an ordering uniquely determines the arrangement of $\mathbf{T}_A$. If $x \in A$ then the approximation $F(n, x)$ determines a unique node from $\mathbf{T}_A$. This node is on the $(n+1)$th level of $\mathbf{T}_A$. Hence $<_A$ is a linear ordering of the $n$th approximations of the elements of $A$ for any $n$.

As $B$ is well-orderable there exists an arrangement $<_B$ of $\mathbf{T}_B$ which well-orders $B$. From $<_A$ and $<_B$ we shall construct an arrangement $<_{A\times B}$ of $\mathbf{T}_{A\times B}$ which will well-order $A \times B$. By Theorem 4.16 we have $F(n + 2, (a, b)) = (F(n, a), F(n, b))$ for any $a$, $b$ and $n$. Let $a, a' \in A$ and $b, b' \in B$. We shall define that $F(n + 2, (a', b')) <_{A\times B} F(n + 2, (a, b))$ holds if $F(n, a') <_A F(n, a)$ or if $F(n, a') = F(n, a)$ and $F(n, b') <_B F(n, b)$. It is not difficult to see that $<_{A\times B}$ linearly orders the successors of every node of $\mathbf{T}_{A\times B}$. Therefore $<_{A\times B}$ is an arrangement of $\mathbf{T}_{A\times B}$. Every arrangement of the developing tree of some set $X$ linearly orders $X$. Every subset $Y$ of $X$ has its minimal element with regard to any arrangement of $\mathbf{T}_X$. This element is represented by the utmost left infinite path in $\mathbf{T}_Y$. To prove that $<_{A\times B}$ well-orders $A \times B$ we must only show that every point from $A \times B$ is upward isolated in $\mathcal{T}_a$ with regard to $<_{A\times B}$.

Let $(a, b) \in A \times B$. We shall prove that there exists such $n$ that $F(n, (a', b')) = F(n, (a, b))$ implies $F(n + m, (a', b')) <_{A\times B} F(n + m, (a, b))$ for every $m > 1$ if only $(a', b') \in A \times B$ and $(a', b') \neq (a, b)$. This will prove that $(a, b)$ is an upward isolated point in $\mathcal{T}_a$ in $A \times B$ with regard to $<_{A\times B}$. As $<_A$ well-orders $A$ there exists such $n_1$ that $F(n_1, a') = F(n_1, a)$ implies $F(n_1 + m, a') <_A F(n_1 + m, a)$ for $m > 1$. And because $<_B$ well-orders $B$ there exists such $n_2$ that $F(n_2, b') = F(n_2, b)$ implies $F(n_2 + m, b') <_A F(n_2 + m, b)$ for $m > 1$. Let $n_3 = \max \{n_1, n_2\}$. If $F(n_3 + 2, (a', b')) = F(n_3 + 2, (a, b))$ then $F(n_3, a') = F(n_3, a)$ and $F(n_3, b') = F(n_3, b)$. Since $(a', b') \neq (a, b)$ it holds either $a' \neq a$ or $b' \neq b$. If $a' \neq a$ then $F(n_3 + m, a') <_A F(n_3 + m, a)$ for $m > 1$. If, however, $b' \neq b$ then $F(n_3 + m, b') <_A F(n_3 + m, b)$ for $m > 1$. In any case $F(n_3 + 2, (a', b')) <_{A\times B} F(n_3 + 2, (a, b))$. This implies that $(a, b)$ is an upward isolated point in $A \times B$ with regard to $<_{A\times B}$. Since $(a, b)$ can be any element of $A \times B$, every element from $A \times B$ is an upward isolated point in $\mathcal{T}_a$ with regard to $>_{A\times B}$. Therefore $>_{A\times B}$ is an arrangement of $\mathbf{T}_{A\times B}$ which well-orders $A \times B$. □

**10.5 Theorem** *The Cartesian product $A_1 \times A_2 \times \ldots \times A_k$ can be well-ordered precisely when the sets $A_1, A_2, \ldots, A_k$ are well-orderable.*

*Proof.* Using induction on $k$ we can prove from Theorem 10.4 and Definition 10.2 that the Cartesian product $A_1 \times A_2 \times \ldots \times A_k$ is well-orderable if every one of the sets $A_1, A_2, \ldots, A_k$ is well-orderable. If, however, one of the sets $A_1, A_2, ..., A_n$ is not well-orderable, then $A_1 \times A_2 \times \ldots \times A_k$ is not well-orderable either. Let $A_i$ not be well-orderable for some $i \leq k$ and let $m_j, j \leq k$, be the minimal elements of $A_j$ with regard to the lexicographical ordering of $A_j$. Then the set $(m_1, ..., m_{i-1}) \times A_i \times (m_{i+1}, ...., m_k)$ is a subset of $A_1 \times A_2 \times \ldots \times A_k$. It is not difficult to see that this subset is in bijection with $A_i$ and is therefore not well-orderable by Lemma 7.14. Hence, the Cartesian product $A_1 \times A_2 \times \ldots \times A_k$ contains a non-well-orderable subset and for this reason it is, by Lemma 7.5, not well-orderable. □

Before we determine precise conditions for a product of sets to be well-orderable we need to prove some additional facts.

**10.6 Theorem** *Every infinite set contains a countable subset.*

*Proof.* Let $A$ be an infinite set. By Lemma 6.5 there exists an accumulation point $x \in A$ in $\mathcal{T}_a$. Therefore there exist $n_0$ and $y_0 \in A$ such that $F(n_0, x) = F(n_0, y_0)$ and $F(n_0 + 1, x) \neq F(n_0 + 1, y_0)$. And there exist $n_1 > n_0$ and $y_1 \in A$ such that $F(n_1, x) = F(n_1, y_1)$ and $F(n_1 + 1, x) \neq F(n_1 + 1, y_1)$ ..., etc. By this construction we obtain points $y_0, y_1, y_2, \ldots$ of $A$ all being different from each other. The only accumulation point of these points in $\mathcal{T}_a$ is $x$. Therefore, by Theorem 4.7.d, the smallest set to which belong $y_0$, $y_1$, $y_2$, ... is a subset of $A$ which, beside the points $y_0$, $y_1$, $y_2$, ..., contains only the point $x$. This subset is in bijection with $\boldsymbol{N}$. Such a bijection is, for example, the function which maps $i$ to $y_i$ for every $i$, and maps $\infty$ to $x$. □

**10.7 Lemma** *$F_n(A \cup B) = F_n(A) \cup F_n(B)$ for every set A and B and every n.*

*Proof.* The lemma holds trivially for $n = 0$, since $F_0(A) = \emptyset$ for every $A$. If $n > 0$ we have

$$F_{n+1}(A \cup B) = \{F_n(x) \mid x \in A \cup B\} = \{F_n(x) \mid x \in A \text{ or } x \in B\} =$$

$$= \{F_n(x) \mid x \in A\} \cup \{F_n(x) \mid x \in B\} = F_{n+1}(A) \cup F_{n+1}(B).$$ □

**10.8 Lemma** *Let $I = \{(j, A_j) \mid j \in J, A_j \in \mathcal{F}\}$ and $I' = \{(j, A_j) \mid j \in J', A_j \in \mathcal{F}\}$ be two indexed systems of sets where $J' \subseteq J$. If the elements of $\prod_{j\in J'} A_j$ are only choice functions on $I'$, and if $A_j =_c 1$ for all $j \in J, j \notin J'$, then $\prod_{j\in J} A_j$ is in bijection with $\prod_{j\in J'} A_j$.*

*Proof.* Let's denote $\prod_{j\in J} A_j$ by $P$ and $\prod_{j\in J'} A_j$ by $P'$. Let $x \in P$. Then for every $n$ there exists a choice function $f_n$ on $I$ (See Definition 10.1) such that $F(n, x) = F(n, f_n)$. The function $f_n$ can be expressed as $f_n = f_c \cup f_n'$ where $f_c$ is a choice function on $I$ with domain $J \setminus J'$, while $f_n'$ is a choice function on $I'$ with domain $J'$, and $f_c$ and $f_n'$ coincide on the intersection of their domains. The proof of this is the same as in ZF. (The intersection of $J \setminus J'$ and $J'$ is not necessarily empty! By Theorem 5.4 these two sets have in common their boundary points in $\mathcal{T}_a$.) The choice function $f_c$ is for every $j \in J, j \notin J'$, constant since $A_j =_c 1$ for every such $j$. Lemma 10.7 implies

$$F(n, x) = F(n, f_n) = F(n, f_c \cup f_n') = F(n, f_c) \cup F(n, f_n') \qquad (1)$$

The sequence $F(0, f_0')$, $F(1, f_1')$, $F(2, f_2')$, ... is the developing sequence of some $f' \in P'$. Hence (1) and Lemma 10.7 imply $F(n, x) = F(n, f_c) \cup F(n, f') = F(n, f_c \cup f')$, which means $x = f_c \cup f'$. Every element of $P'$ is by assumption a choice function on $I'$. Hence $f'$ is a choice function on $I'$. Therefore $x$ being the union of $f_c$ and $f'$ is a choice function on $I$. So the elements of $P$ are only choice functions on $I$.

Let's now prove that $P$ and $P'$ are in bijection. Let $f_1, f_2 \in P$. Then, as we have just proved, we have

$$f_1 = f_1' \cup f_c \qquad (2)$$

and

$$f_2 = f_2' \cup f_c \qquad (3)$$

where $f_1', f_2' \in P'$. If $f_1' \neq f_2'$ then we can assume, without loss of generality, that there exists at least one pair $(j', a_j') \in f_1'$ which is not an element of $f_2'$. This pair is in $f_1$ because of (2). This pair cannot be in $f_c$ otherwise (3) would imply that $f_2$ is not a function. Namely, $f_1'$ and $f_c$ would in this case map $j'$ to $a_j'$, while $f_2'$ would map $j'$ to something else since $(j', a_j') \notin f_2'$ but $j'$ is in the domain of $f_2'$. Therefore $(j', a_j') \notin f_2$. Hence $f_1 \neq f_2$. Let now $f_1 \neq f_2$. Then (2) and (3) imply $f_1' \neq f_2'$. Thus different elements of $P$ correspond by (2) and (3) to different elements of $P'$ and vice versa. □

In the following theorems natural numbers will be represented as N-ordinals, that is, in von Neumann form and not in Zermelo form because a natural number $n$ in Zermelo form is a set with only one element while the ordinal $n$ is a set with $n$ elements and can thus represent the cardinality of the sets having $n$ elements.

**10.9 Theorem** *The product* $\prod_{j\in J} A_j$ *is well-orderable precisely when* $A_j >_c 1$ *holds for finitely many sets* $A_j$ *and when every* $A_j$, $j \in J$, *is well-ordered.*

*Proof.* If one of the sets $A_j$ is empty then $\prod_{j\in J} A_j$ is empty and as such is trivially well-orderable. Otherwise $A_j \geq_c 1$ for every $j \in J$. Let only for finitely many sets $A_j$ hold $A_j >_c 1$, and let these sets be $A_{j(1)}, A_{j(2)}, \ldots, A_{j(n)}$. Then $A_j =_c 1$ for every $j \neq j(i)$, $i \leq n$, and so $\prod_{j\in J} A_j$ is in bijection with $\prod_{1\leq i\leq n} A_{j(i)}$ by Lemma 10.8. The product $\prod_{1\leq i\leq n} A_{j(i)}$ is in bijection with the product $A_{j(1)} \times A_{j(2)} \times \ldots \times A_{j(n)}$ by Theorem 10.3. So, if every set $A_{j(i)}$, $i \leq n$, is well-orderable then the product $\prod_{j\in J} A_j$ is well-orderable by Theorems 10.3 and 10.5 and Lemma 7.14. Let some set $A_{j(i)}$, where $i \leq n$, not be well-orderable. The set $(a_1, a_2, ..., a_{i-1}) \times A_{j(i)} \times (a_{i+1}, ... a_n)$, $a_k \in A_{j(k)}$ for $k \leq n$ and $k \neq j(i)$, is in bijection with $A_{j(i)}$ and is thus not well-orderable by Lemma 7.14. Since this set is a subset of the set $A_{j(1)} \times A_{j(2)} \times \ldots \times A_{j(n)}$ this set is not well-orderable. And since this set is in bijection with the set $\prod_{j\in J} A_j$, the set $\prod_{j\in J} A_j$ is not well-orderable.

Let now for infinitely many sets $A_j$ hold $A_j >_c 1$. In other words let there exist an infinite subset $J'$ of $J$ such that $A_j >_c 1$ for every $j \in J'$. Thus $A_j$ has at least two elements, say $a_{j,0}$ and $a_{j,1}$. This holds for every $j \in J'$. On the other hand every $A_j$, $j \in J \setminus J'$, has at least one element, say $b_j$, as in the opposite case $\prod_{j\in J} A_j$ would be empty. Every function $f$ which maps $j \in J'$ to $a_{j,0}$ or to $a_{j,1}$, and maps $j$ to $b_j$ for $j \in J$, $j \notin J'$, is an element of $\prod_{j\in J} A_j$. By Theorem 10.6 there exists a countable $J'' \subseteq J'$, such that $A_j >_c 1$ for every $j \in J''$. Let $j_1, j_2, j_3$ ... be all elements of $J''$. All $j_i$ for $i > k$ have the same $k$th approximation $F(k, j_i)$ which, together with $k$th approximation of $a_{i,0}$ or $a_{i,1}$, constitutes the $(k + 2)$th approximation of $f$ since $F(k + 2, f) \subseteq \{(F(k, j_i), F(k, a_{i,m})) \mid m = 0 \vee m = 1\}$. And since there are only finitely many different approximations of the elements of $F(k + 2, f)$, it is not recognizable from $F(k + 2, f)$ for infinitely many $j_i$ how $f$ maps them. Hence, by increasing $k$, we can at a definite $i$th approximation $F(i, f)$ choose arbitrarily whether $f$ maps $j_i \in J''$, $i > k$, to $a_{i,0}$ or to $a_{i,1}$. And, by increasing $k$ further, this choice is possible infinitely many times again and again. In this way we obtain a perfect binary tree whose infinite paths are the developing sequences of the described functions $f$. This tree is the developing tree of some set which is perfect in $\mathcal{T}_a$ and is a subset of $\prod_{j\in J} A_j$. Consequently $\prod_{j\in J} A_j$ is not well-orderable by Theorem 7.4 and Lemma 7.5. □

This theorem enables us to show at last that the product of sets is legitimate. That is, it enables us to show that the cardinality of the product is determined by the cardinalities of its factors.

**10.10 Theorem** *If* $I = \{(j, A_j) \mid j \in J, A_j \in \mathcal{F}\}$ *and* $I' = \{(j, A_j') \mid j \in J, A_j' \in \mathcal{F}'\}$ *are such indexed systems of sets that* $A_j =_c A_j'$ *for every* $j \in J$ *then* $\prod_{j\in J} A_j =_c \prod_{j\in J} A_j'$.

*Proof.* Let $I$ and $I'$ be indexed systems of sets described in the theorem. Let only for finitely many sets $A_j$ of $I$ hold $A_j >_c 1$, and let these sets be $A_1, A_2, \ldots A_n$. Then, by Lemma 10.8 and by Theorem 10.3, it must be true that $\prod_{j\in J} A_j =_c A_1 \times A_2 \times \ldots \times A_n$. And because of $A_j =_c A_j'$ for every $j \in J$, it must also hold $\prod_{j\in J} A_j' =_c A'_1 \times A'_2 \times \ldots \times A'_n$. In the same way as in ZF we prove that $A_1 =_c A_1'$ and $A_2 =_c A_2'$ imply $A_1 \times A_2 =_c A_1' \times A_2'$. Thus, by induction on $n$, we have $A_1 \times A_2 \times \ldots \times A_n =_c A'_1 \times A'_2 \times \ldots \times A'_n$. All these equations taken together give $\prod_{j\in J} A_j =_c \prod_{j\in J} A_j'$.

If, however, for infinitely many sets $A_j$ of $I$ hold $A_j >_c 1$ then $\prod_{j\in J} A_j$ and $\prod_{j\in J} A_j'$ contain a perfect subset by Theorem 10.9. Hence again $\prod_{j\in J} A_j =_c \prod_{j\in J} A_j'$. □

Now we can also define exponentiation for arbitrary sets. The expression $A^B$ denotes in ZF the set of all functions from $B$ to $A$. Of course, in AS the set of all these functions does in general not exist if $B$ is infinite as we have already shown for the case where $A = \{0, 1\}$ and $B = N$. For this reason we must set the following definition:

**10.11 Definition** The set $A^B$ is the smallest set whose elements are all functions $f : B \rightarrow A$.

The set $A^B$ is a special case of the product of an indexed system $\{(j, A_j) \mid j \in J, A_j \in \mathcal{F}\}$, where $J = B$ and $A_j = A$ for every $j \in J$. This product is legitimate by Theorem 10.10. So $A =_c A'$ and $B =_c B'$ imply $A^B =_c A'^{B'}$. If $B$ is a finite set of $n$ elements, then $A^B$ is, by Theorem 10.3, in bijection with the product $A \times A \times \ldots \times A$ which has $n$ factors.

Exponentiation of finite sets matches with the exponentiation of corresponding finite N-ordinals in the ordinal arithmetic. In ZFC it is pretty difficult to determine the cardinality of the set $A^B$ when $B$ is infinite. This problem has not yet been solved in general ([H&J], p. 164). In AS the exponentiation with infinite sets is trivial. Let $B$ be an infinite set. If $A =_c 0$ then $A^B =_c 0$. If $A =_c 1$ then $A^B =_c 1$. And if $A \geq_c 2$ then $A^B$ is not well-orderable by Theorem 10.9. Hence $A^B =_c \mathcal{P}(\omega)$.

Let's now take a look at the cardinalities of some sets which are expressed by sums and products. We shall see that cardinal arithmetic in AS is pretty different from cardinal arithmetic in ZFC.

**10.12 Lemma** *For any well-ordered set A and any finite N-ordinal* $n > 0$ *holds* $r_{CB}(n \times A) =_o r_{CB}(A)$ *and* $d_{CB}(n \times A) =_o n \times d_{CB}(A)$

*Proof.* As $n$ is a N-ordinal different from 0 we have by Lemma 8.12

$$(n \times A)^{(\alpha)} = (\cup_{1\leq i<n} \{i\} \times A)^{(\alpha)} = \cup_{1\leq i<n} (\{i\} \times A)^{(\alpha)} \qquad (1)$$

for every ordinal $\alpha$. By Lemma 9.7 we have $(\{i\} \times A)^{(\alpha)} =_c A^{(\alpha)}$ for every ordinal $\alpha$ and every $i$. Therefore, if $\alpha =_o r_{CB}(A)$ then $(\{i\} \times A)^{(\alpha)} = d_{CB}(\{i\} \times A)$ is a finite nonempty set for every $i$. Consequently the set $\cup_{1\leq i<n} (\{i\} \times A)^{(\alpha)}$ is finite, and (1) implies that the set $(n \times A)^{(\alpha)}$ is finite as well. Hence $r_{CB}(n \times A) =_o r_{CB}(A)$. Since the sets $\{i\} \times A$ are mutually disjoint for any $i$, the sets $(\{i\} \times A)^{(\alpha)}$ are mutually disjoint too. Therefore $\alpha =_o r_{CB}(A)$ implies

$$d_{CB}(n \times A) = (n \times A)^{(\alpha)} = \cup_{1\leq i<n} (\{i\} \times A)^{(\alpha)} = \cup_{1\leq i<n} d_{CB}(\{i\} \times A) =_c n \times d_{CB}(A)$$

since $\{i\} \times A =_c A$ for every $i$. Thus we have proved the lemma. □

**10.13 Theorem** *For any well-ordered set* $A \neq \emptyset$ *and any finite N-ordinal n holds* $n \times A <_c (n + 1) \times A$.

*Proof.* Obviously $n \times A \leq_c (n + 1) \times A$ since $n \times A \subseteq (n + 1) \times A$. But $n \times A =_c (n + 1) \times A$ does not hold. Namely, we have $r_{CB}(n \times A) =_o r_{CB}(A) =_o r_{CB}((n + 1) \times A)$ by Lemma 10.12. Since $d_{CB}(A) >_o 0$ if $A \neq \emptyset$, we have $d_{CB}(n \times A) =_o n \times d_{CB}(A) <_o (n + 1) \times d_{CB}(A) =_o d_{CB}((n + 1) \times A)$ by Lemma 10.12. Therefore we have $n \times A <_c (n + 1) \times A$ by Theorem 8.21. □

If the sets $A_1$ and $A_2$ are infinite, then in ZFC we have $A_1 + A_2 =_c A_1 \times A_2$. The AS cardinal arithmetic is not so trivial if both $A_1$ and $A_2$ are well-orderable.

**10.14 Theorem** *If the sets* $A_1$ *and* $A_2$ *are infinite and well-orderable, then* $A_1 + A_2 <_c A_1 \times A_2$.

*Proof.* Since $A_1$ and $A_2$ are infinite, they cannot be of smaller cardinality then $\omega$ since $\omega$ as a countable set is of the smallest infinite cardinality by Theorem 10.6. Therefore we can assume, without loss of generality that $\omega \leq_c A_1 \leq_c A_2$. Then $A_1 + A_2 \leq_c A_2 + A_2 =_c 2 \times A_2 <_c 3 \times A_2 \leq_c \omega \times A_2 \leq_c A_1 \times A_2$ by Theorem 10.13. □

In ZFC is in force *König's Theorem* which and asserts the following. If $\{(j, A_j) \mid j \in J, A_j \in \mathcal{F}\}$ and $\{(j, B_j) \mid j \in J, B_j \in \mathcal{G}\}$ are indexed systems of sets such that $A_j <_c B_j$ for every $j \in J$, then we have $\sum_{j\in J} A_j <_c \prod_{j\in J} B_j$. König's Theorem is true in AS only under certain conditions.

**10.15 Theorem** *König's Theorem holds if and only if the set of indexes J is well-orderable.*

*Proof.* Obviously every $A_j, j \in J$, must be a well-orderable set. If there exists some non-well-orderable $A_j$ then $A_j <_c B_j$ cannot hold and the conditions of König's Theorem are not fulfilled. And, if $J$ is not well-orderable, then $\sum_{j\in J} A_j$ is not well-orderable either (Theorem 9.6), and König's Theorem does not hold. Let for infinitely many $j \in J$ hold $A_j >_c 0$. Then for infinitely many $B_j$ holds $B_j >_c 1$ and consequently $\prod_{j\in J} B_j =_c \mathcal{P}(\omega)$ by Theorem 10.9, while $\sum_{j\in J} A_j$ is a well-orderable set by Theorem 9.6. Therefore König's Theorem holds in this case. König's Theorem holds also if $A_j >_c 0$ only for finitely many $j \in J$. In this case the smallest set $M$, whose elements are all such $j$, is finite. Thus $M$ has no accumulation points in $\mathcal{T}_a$. Hence every element of $M$ is such $j$. If $M = \emptyset$ then $\sum_{j\in J} A_j =_c 0$ while $\prod_{j\in J} B_j >_c 0$ and König's Theorem holds. If $M \neq \emptyset$ then let $A_m$ be the greatest $A_j$ for $j \in M$ with regard to $<_c$. Because $B_j =_c 1$ for every $j \in J, j \notin M$, we have $\prod_{j\in J} B_j =_c \prod_{j\in M} B_j$ by Lemma 10.8. Further we have $\prod_{j\in M} B_j \geq_c 2^{M\setminus\{m\}} \times B_m$ because $B_j \geq_c 2$ for $j \in M \setminus \{m\}$. Here is $2^{M\setminus\{m\}}$ the exponentiation of the N-ordinal 2 by the set $M \setminus \{m\}$. Then we have $2^{M\setminus\{m\}} \times B_m >_c 2^{M\setminus\{m\}} \times A_m$. Namely, it certainly holds $2^{M\setminus\{m\}} \times B_m \geq_c 2^{M\setminus\{m\}} \times A_m$ but we cannot have $2^{M\setminus\{m\}} \times B_m =_c 2^{M\setminus\{m\}} \times A_m$ since Theorem 8.20 and Lemma 10.12 would then imply

$$r_{CB}(A_m) =_o r_{CB}(2^{M\setminus\{m\}} \times A_m) =_o r_{CB}(2^{M\setminus\{m\}} \times B_m) =_o r_{CB}(B_m) \qquad (1)$$

$$2^{M\setminus\{m\}} \times d_{CB}(A_m) =_o d_{CB}(2^{M\setminus\{m\}} \times A_m) =_o d_{CB}(2^{M\setminus\{m\}} \times B_m) =_o 2^{M\setminus\{m\}} \times d_{CB}(B_m) \qquad (2)$$

From (1) we obtain $r_{CB}(A_m) =_o r_{CB}(B_m)$ and from (2) we obtain $d_{CB}(A_m) =_o d_{CB}(B_m)$. Thus $A_m =_c B_m$ by Theorem 8.20 which is not true by assumption. Further we have $2^{M\setminus\{m\}} \times A_m \geq_c M \times A_m$ since $2^{n-1} \geq_c n$ for every finite N-ordinal $n > 0$ which can be proved by induction on $n$. And $M \times A_m \geq_c \sum_{j\in M} A_j$ since $A_m$ has the greatest cardinality among the sets $A_j$ for $j \in M$. Finally we have $\sum_{j\in M} A_j =_c \sum_{j\in J} A_j$ since $A_j =_c 0$ for every $j \in J, j \notin M$. If we resume:

$$\prod_{j\in J} B_j =_c \prod_{j\in M} B_j \geq_c 2^{M\setminus\{m\}} \times B_m >_c 2^{M\setminus\{m\}} \times A_m \geq_c M \times A_m \geq_c \sum_{j\in M} A_j =_c \sum_{j\in J} A_j$$

In short $\prod_{j\in J} B_j >_c \sum_{j\in J} A_j$. This proves the theorem completely. □

# 11. The Universe of all Cardinals

Till now we have developed the cardinal arithmetic in AS without explicit definition of the concept of cardinal since this concept was not needed for our purposes. Now we can no longer avoid this concept. In ZF a cardinal is an ordinal which is not in bijection with any of its proper initial segments. But in AS every set is not well-orderable. However, for every set $A$ there should exist such cardinal (number) $\kappa$ that $A =_c \kappa$. Taking into account that all sets which are not ordinals are of the same cardinality we shall define cardinals as follows.

**11.1 Definition** A *cardinal* is

1) an ordinal $\kappa$ such that $\kappa =_c S$ does not hold for any proper initial segment $S$ of $\kappa$
2) a set which cannot be well-ordered.

**11.2 Theorem** *For every set A there exists such cardinal* $\kappa$ *that* $A =_c \kappa$.

*Proof.* If $A = \emptyset$ then $A =_c \emptyset$ where $\emptyset$ is a cardinal. Let now $A$ be a nonempty well-orderable set. Let $\mathbf{P}(x) \Leftrightarrow x \in A \wedge A[x] \cup \{x\} =_c A$. Certainly $\mathbf{P}(\max(A))$ holds true. Hence, by Lemma 8.8, there exists the minimal $x' \in A$ for which $\mathbf{P}$ holds true. Then $A\,[x'] \cup \{x'\}$ is by Definition 11.1 a cardinal. Let finally $A$ be a non-well-orderable set. Then $A$ is a cardinal by Definition 11.1 and the equality $A =_c A$ proves the theorem completely. □

The sets $\boldsymbol{N}$ and $\omega$ have the same structure (See the proof of Theorem 8.5). Hence the proof that $\mathcal{P}(\omega)$ contains a perfect subset is a natural adaptation of the proof that $\mathcal{P}(\boldsymbol{N})$ contains a perfect subset by Theorem 5.13. Thus $A =_c \mathcal{P}(\omega)$ holds for every non-well-orderable set $A$. In the sequel we shall see that all interesting well-orderable cardinals can be represented by the sets $n \times \omega$ and $\omega^n$ where $n$ is a finite N-ordinal. For this reason $\mathcal{P}(\omega)$ will be from now on the standard representative of the non-well-orderable cardinals.

Let's show that all cardinals representing the same well-orderable set are isomorphic to each other.

**11.3 Theorem** *If $\kappa_1$ and $\kappa_2$ are well-ordered cardinals and $\kappa_1 =_c \kappa_2$ then $\kappa_1 =_o \kappa_2$.*

*Proof.* Let $\kappa_1$ and $\kappa_2$ be well-ordered cardinals for which holds $\kappa_1 =_c \kappa_2$. By Theorem 7.12 one of them is isomorphic to an initial segment of the other. Without loss of generality we can presume that $\kappa_1$ is isomorphic to an initial segment $S$ of $\kappa_2$. That is, $\kappa_1 =_o S$. Therefore $\kappa_2 =_c S$. Since $\kappa_2$ is a cardinal, $S$ cannot be its proper initial segment. Hence $S = \kappa_2$ which implies $\kappa_1 =_o \kappa_2$. □

The cardinality of the continuum is the cardinality of the set $\mathcal{P}(\omega)$. The continuum problem is to determine how do non-isomorphic cardinals follow each other with regard to $<_c$ and where in this sequence is the cardinal $\mathcal{P}(\omega)$. The second part of this question is easy to answer. The cardinal $\mathcal{P}(\omega)$ is at the very end of this sequence since there are no greater cardinals with regard to $<_c$. The concepts of sum and product of sets, which we have developed in the previous two sections, enable us to answer also the first part of the question, which implicitly includes the question "how many" mutually non-isomorphic cardinals precede $\mathcal{P}(\omega)$.

**11.4 Theorem** *The set $n \times \kappa$ is a cardinal for any well-ordered cardinal $\kappa$ and any finite N-ordinal $n$.*

*Proof.* Let $<_\kappa$ be a well-ordering of $\kappa$ such that $\max(\kappa) = \max(d_{CB}(\kappa))$. This is always possible by Lemma 8.19. First we shall well-order the product $n \times \kappa$ in the following way. If $(i, a), (i', a') \in n \times \kappa$ then let $(i, a) < (i', a')$ if $i <_o i'$ or if $i = i'$ and $a <_\kappa a'$ where As $n$ is a finite set it is not difficult to see that the described linear ordering $<$ is a well-ordering of $n \times \kappa$. By proving that $\max(n \times \kappa)$, with regard to $<$, is an element of $d_{CB}(n \times \kappa)$ we prove that $n \times \kappa$ is a cardinal. Namely, let $S$ be a proper initial segment of $n \times \kappa$ . Certainly $r_{CB}(S) \leq_o r_{CB}(n \times \kappa)$. If $r_{CB}(S) <_o r_{CB}(n \times \kappa)$ then $S <_c n \times \kappa$ by Theorem 8.21. And if $r_{CB}(S) =_o r_{CB}(n \times \kappa)$ then $d_{CB}(S) \subseteq d_{CB}(n \times \kappa)$ by Definition 8.14 and Lemma 8.18. As $d_{CB}(n \times \kappa)$ is a finite set whose element is by assumption $\max(n \times \kappa)$ which, of course, cannot be an element of $d_{CB}(S)$, we have $d_{CB}(S) <_c d_{CB}(n \times \kappa)$. And again Theorem 8.21 implies $S <_c n \times \kappa$.

So it remains to prove that $\max(n \times \kappa) \in d_{CB}(n \times \kappa)$. The set $n \times \kappa$ is the union of mutually disjoint sets $\{i\} \times \kappa$ where $i <_o n$ since $n$ is a N-ordinal. Therefore we have $\max(n \times \kappa) = \max(\{n-1\} \times \kappa)$ with regard to $<$. The set $\{n-1\} \times \kappa$ is a cardinal since the map $a \to (n-1, a)$, where $a$ is an arbitrary element of $\kappa$, is a $<$ isomorphism between $\kappa$ and $\{n-1\} \times \kappa$. For this reason $\max(\kappa) = \max(d_{CB}(\kappa))$ imply $\max(\{n-1\} \times \kappa) = \max(d_{CB}(\{n-1\} \times \kappa)) \in d_{CB}(\{n-1\} \times \kappa)$. Therefore the following holds $\max(n \times \kappa) = \max(\{n-1\} \times \kappa) \in d_{CB}(\{n-1\} \times \kappa) \subseteq d_{CB}(n \times \kappa)$. In short $\max(n \times \kappa) \in d_{CB}(n \times \kappa)$. □

**11.5 Theorem** *If $d_{CB}(A) =_o 1$ holds for some well-ordered set $A$, then there is no set between $n \times A$ and $(n+1) \times A$ with regard to $<_c$ for any finite N-ordinal $n > 0$.*

*Proof.* If $d_{CB}(A) =_o 1$, then $d_{CB}(n \times A) =_o n$ and $d_{CB}((n+1) \times A) =_o n+1$ by Lemma 10.12. Let $n \times A <_c B <_c (n+1) \times A$ for some set $B$ and some $n >_o 0$. Since $r_{CB}(n \times A) =_o r_{CB}((n+1) \times A) =_o r_{CB}(A)$ by Lemma 10.12, we have $r_{CB}(B) =_o r_{CB}(A)$ and $n <_o d_{CB}(B) <_o n+1$ by Theorem 8.21. The last inequality, of course, cannot hold. □

**11.6 Remark** If $d_{CB}(A) \neq_o 1$, then the previous theorem does not hold. Let, for example, $A =_c n \times A'$ where $n >_o 1$. Then $d_{CB}(A) = d_{CB}(n \times A') = n \times d_{CB}(A') >_o 1$ by Lemma 10.12, and by Theorem 10.13 we have $A =_c n \times A' <_c (n + 1) \times A' <_c (n + 2) \times A' <_c ... <_c (2n - 1) \times A' <_c (2n) \times A' =_c 2 \times A$.

**11.7 Definition** A cardinal $\kappa'$ is a *successor of a cardinal* $\kappa$ if $\kappa <_c \kappa'$ and if there is no cardinal $\kappa''$ such that $\kappa <_c \kappa'' <_c \kappa'$. In this case the cardinal $\kappa$ is a *predecessor* of the cardinal $\kappa'$. If a nonempty cardinal does not have a predecessor, then it is a *limit cardinal*. A cardinal $\kappa$ *follows immediately after the cardinals* $\kappa_0, \kappa_1, \kappa_2, \kappa_3, ...$ if $\kappa_i <_c \kappa$ for every $i$ and if for every cardinal $\kappa'$ for which holds $\kappa_i <_c \kappa'$ for every $i$ also holds $\kappa <_c \kappa'$.

A cardinal can have several successors. All these are isomorphic to each other. The same holds for the predecessors of a cardinal. This definition and Theorems 11.4 and 11.5 imply the following. If $\kappa$ is some cardinal for which $d_{CB}(\kappa) =_c 1$ then $\kappa$ is a predecessor of the cardinal $2 \times \kappa$ which is a predecessor of the cardinal $3 \times \kappa$ ..., etc. And now we shall prove that $\omega \times \kappa$ is a cardinal which follows immediately after the cardinals $\kappa, 2 \times \kappa, 3 \times \kappa, ....$

**11.8 Theorem** *If $\kappa$ is a well-ordered cardinal then $\omega \times \kappa$ is a limit well-ordered cardinal which follows immediately after the cardinals* $\kappa, 2 \times \kappa, 3 \times \kappa, ....$

*Proof.* Let $P$ be some perfect set whose developing tree is a perfect binary tree. Such set is, for example, the perfect subset of $\mathcal{P}(\boldsymbol{N})$ constructed in the proof of Theorem 5.13. Let $v_0, v_1, v_2, v_3, ...$ be the consequent nodes on the utmost right infinite path of $\mathbf{T}_P$ where $v_0$ is the root of $\mathbf{T}_P$. Since $\mathbf{T}_P$ is a perfect binary tree each node of $\mathbf{T}_P$ has precisely two successors. Let $v_i'$ be the left successor of $v_i$ (the right successor is $v_{i+1}$). From every $v_i'$ develops in $P$ a perfect set $P_i$. The sets $P_0, P_1, P_2, ...$ are mutually separated by Definition 6.1.b. The cardinal $\kappa$ can be embedded into each perfect set $P_i$ by Lemma 6.2. Let $\kappa_i$ be the embedding of $\kappa$ into $P_i$. Since the sets $P_i$ are mutually separated, the sets $\kappa_i$ are mutually separated as well. Let's preserve only those infinite paths from $\mathbf{T}_P$ which represent an element of some $\kappa_i$. In this way we obtain the developing tree of some set $A \subseteq P$. Since $\kappa$ is well-ordered we can, by Theorem 7.7, well-order each $\kappa_i$ by an arrangement of the developing tree of $\kappa_i$. And since $\kappa_i$ are mutually separated subsets of $A$ we obtain in this way a well-ordering of $A$. Let's denote this well-ordering by <. Then we have $i \times \kappa =_c \kappa_1 + \kappa_2 + \kappa_3 + ... + \kappa_i \leq_c A$ for every $i$. Therefore $i \times \kappa <_c (i + 1) \times \kappa \leq_c A$ for every $i$.

Let's presume that there exists such set $A'$ that $i \times \kappa <_c A' <_c A$ for every $i$. Then $A'$ is well-orderable and must be < isomorphic to a proper initial segment $S$ of $A$. Hence $\max(S) < \max(A)$. Let $<_A$ be the order of the nodes of $\mathbf{T}_A$ which determines the well-ordering < of $A$. There must exist such $n$ that $F(n, \max(S)) = F(n, \max(A))$ and $F(n + 1, \max(S)) <_A F(n + 1, \max(A))$. The way in which we have embedded the cardinal $\kappa$ into each $P_i$ implies that the path for $\max(S)$ in $\mathbf{T}_A$ is in the developing tree of $\kappa_n$. Therefore $A' =_o S \leq_c \kappa_1 + \kappa_2 + ... + \kappa_n =_c n \times \kappa <_c (n + 1) \times \kappa$. This contradicts the assumption that $i \times \kappa <_c A'$ for every $i$. Hence a cardinal representing $A$ follows immediately after the cardinals $\kappa, 2 \times \kappa, 3 \times \kappa, .....$ And $A$ is in fact a cardinal. Namely, $A$ cannot be in bijection with any of its proper initial segments $S$ since, by the previous argument, we have $S \leq_c n \times \kappa <_c A$ for some $n <_o \omega$. And $A$ is also a limit cardinal since if $\kappa' <_c A$ holds for some cardinal $\kappa'$ then, by the previous argument, there exists such $n$ that $\kappa' <_c n \times \kappa <_c A$. Therefore the cardinal $A$ has no predecessor. □

**11.9 Theorem** *An infinite well-ordered cardinal $\kappa$ is a limit cardinal precisely when* $d_{CB}(\kappa) =_o 1$.

*Proof.* Let $\kappa$ be an infinite well-ordered cardinal such that $d_{CB}(\kappa) =_o 1$ and let $\kappa'$ be a cardinal such that $\kappa' <_c \kappa$. Then $r_{CB}(\kappa') <_o r_{CB}(\kappa)$ by Theorems 8.21 and 8.15.c. Theorem 11.4 implies that $2 \times \kappa'$ is a cardinal as well, and Lemma 10.12 implies that $r_{CB}(2 \times \kappa') =_o r_{CB}(\kappa') <_c r_{CB}(\kappa)$. Hence $\kappa' <_c 2 \times \kappa' <_c \kappa$ by Theorem 8.21. This implies that $\kappa$ is a limit cardinal.

Let now $\kappa$ be a limit cardinal. If $d_{CB}(\kappa) =_o n > 1$, then $d_{CB}(\kappa) = \{a_1, a_2, ..., a_n\}$. Without loss of generality we can assume that $a_1 <_\kappa a_2 <_\kappa ..., <_\kappa a_n$ where $<_\kappa$ is a well-ordering of $\kappa$. The initial segment $\kappa'$ of $\kappa$ whose maximal element with regard to $<_\kappa$ is $a_1$ is a cardinal. Namely, if $\alpha = r_{CB}(\kappa)$ then $a_1 \in \kappa'^{(\alpha)}$ while for every proper initial segment $S$ of $\kappa'$ we have $S^{(\alpha)} = \emptyset$. Hence $r_{CB}(S) <_o r_{CB}(\kappa_i)$ and therefore $S <_c \kappa_i$ by Theorem 8.21. Thus $r_{CB}(\kappa') =_o r_{CB}(\kappa)$ and $d_{CB}(\kappa') =_o 1$. Therefore Theorem 8.20 and Lemma 10.12 imply $\kappa =_c n \times \kappa'$. By Theorem 11.5 there are no cardinals between the cardinals $(n - 1) \times \kappa'$ and $n \times \kappa'$ with regard to $<_c$. Thus $(n - 1) \times \kappa'$ is a predecessor of $\kappa$. But this contradicts the assumption that $\kappa$ is a limit cardinal. □

Now we can conclude the following. The smallest cardinal is 0. A successor of the cardinal 0 is the cardinal 1. A successor of 1 is 2 ..., etc. Theorem 11.8 and $\omega \times 1 =_c \omega$ imply that $\omega$ is a limit cardinal which follows immediately after the cardinals 0, 1, 2, .... Therefore $d_{CB}(\omega) =_o 1$ by Theorem 11.9. Theorems 11.4 and 11.5 imply that $\omega$ is followed by the cardinals $2 \times \omega$, $3 \times \omega$, … etc. Theorem 11.8 implies that a cardinal which follows immediately after these cardinals is $\omega \times \omega$ or $\omega^2$. (Here $\omega^2$ does not denote the exponentiation in the ordinal arithmetic of N-ordinals, where $\omega^2 = \omega$, but the exponentiation in the cardinal arithmetic!). Now Theorems 11.4, 11.5 and 11.8 imply that $\omega^2$ is a limit cardinal followed by the cardinals $2 \times \omega^2$, $3 \times \omega^2$, …. A cardinal which follows immediately after these cardinals is $\omega \times \omega^2$ or $\omega^3$. This cardinal is followed by the cardinals $2 \times \omega^3$, $3 \times \omega^3$, …., etc. Every of the cardinals stated so far can be expressed in the form $m \times \omega^n$ where $m$ and $n$ are finite N-ordinals. It is not difficult to conclude that a cardinal which follows immediately after all cardinals of this form follows also immediately after all the cardinals 1, $\omega$, $\omega^2$, $\omega^3$, .... A natural conjecture would be that $\omega^\omega$ is such a cardinal. However, Theorem 10.9 implies that $\omega^\omega$ is not a well-ordered set. But we shall prove now that there exists a well-ordered cardinal which follows immediately after the cardinals 1, $\omega$, $\omega^2$, $\omega^3$, ....

**11.10 Theorem** *Let $\kappa_0 <_c \kappa_1 <_c \kappa_2 <_c ...$ be well-ordered cardinals such that $\kappa_i + \kappa_{i+1} =_c \kappa_{i+1}$ holds for every i. Then there exists a well-ordered cardinal which follows immediately after $\kappa_0, \kappa_1, \kappa_2, ....$ This cardinal is a limit cardinal.*

*Proof*. We prove this theorem in the same way as Theorem 11.8 with the difference that instead of embedding the same cardinal $\kappa$ into each perfect set $P_i$ we embed the cardinal $\kappa_i$ into $P_i$ obtaining thus $\kappa_i' \subseteq P_i$ for every $i$. We preserve only those infinite paths from $\mathbf{T}_P$ which represent an element of some $\kappa_i'$. In this way we obtain the developing tree of some set $A \subseteq P$. Since each of the cardinals $\kappa_i$ is well-ordered we can, by Theorem 7.7, well-order each $\kappa_i'$ by an arrangement of its developing tree. And since the sets $\kappa_i'$ are mutually separated subsets of $A$ we obtain in this way a well-ordering $<$ of $A$. Besides, $\kappa_i \leq_c A$ holds for every $i$ and hence $\kappa_i <_c \kappa_{i+1} \leq_c A$ for every $i$.

We show in the same way as in the proof of Theorem 11.8 that there does not exist a set $A'$ such that $\kappa_i <_c A' <_c A$ for every $i$. Namely, this would imply that $A' \leq_c \kappa_0 + \kappa_1 + \kappa_2 + ... + \kappa_n =_c \kappa_n <_c \kappa_{n+1}$ holds for some $n <_o \omega$ which contradicts the assumption that $\kappa_i <_c A'$ for every $i$. Hence a cardinal representing $A$ follows immediately after the cardinals $\kappa_0, \kappa_1, \kappa_2, \kappa_3, ....$ And $A$ is itself a cardinal since for any of its proper initial segments $S$ holds $S \leq_c \kappa_n <_c A$ for some $n <_o \omega$. It is also a limit cardinal since if for some cardinal $\kappa$ holds $\kappa <_c A$ then there exists such $n$ that $\kappa <_c \kappa_n <_c A$. (The details of this sketched proof are the same as in the proof of Theorem 11.8.) □

Theorems 11.4, 11.5, 11.8 and 11.9 imply that $d_{CB}(\omega^n) =_o 1$ for every $n$ and that $\omega^n <_c \omega^{n+1}$ for every $n$. Hence $r_{CB}(\omega^n) <_o r_{CB}(\omega^{n+1})$ for every $n$ by Theorem 8.21. Consequently $\omega^n + \omega^{n+1} =_c \omega^{n+1}$ by Theorem 9.8. So the cardinals 1, $\omega$, $\omega^2$, $\omega^3$, .... meet the conditions of Theorem 11.10. Therefore there exists a well-ordered limit cardinal, say $\kappa_1$, which follows immediately after the cardinals 1, $\omega$, $\omega^2$, $\omega^3$, .... Theorem 11.9 implies that $d_{CB}(\kappa_1) =_o 1$, and Theorems 11.4, 11.5 and 11.8 imply that $\kappa_1$ is followed by the cardinals

$2 \times \kappa_1$, $3 \times \kappa_1$, ..., $\omega \times \kappa_1$, $2 \times \omega \times \kappa_1$, $3 \times \omega \times \kappa_1$, ..., $\omega^2 \times \kappa_1$, $2 \times \omega^2 \times \kappa_1$, $3 \times \omega^2 \times \kappa_1$,..., $\omega^3 \times \kappa_1$, ....

The cardinal which follows immediately after all these cardinals is the same as the cardinal which follows immediately after the cardinals $\kappa_1$, $\omega \times \kappa_1$, $\omega^2 \times \kappa_1$, $\omega^3 \times \kappa_1$, .... This sequence meets the conditions of Theorem 11.10 which we conclude in the same way as for the sequence 1, $\omega$, $\omega^2$, $\omega^3$, .... Hence, by Theorem 11.10, there is some well-ordered cardinal $\kappa_2$ which follows immediately after all the cardinals $\kappa_1$, $\omega \times \kappa_1$, $\omega^2 \times \kappa_1$, $\omega^3 \times \kappa_1$, .... Now we repeat the argument and find out that $\kappa_2$ is followed by the cardinals $2 \times \kappa_2$, $3 \times \kappa_2$, ..., $\omega \times \kappa_2$, $2 \times \omega \times \kappa_2$, ..., $\omega^2 \times \kappa_2$, $2 \times \omega^2 \times \kappa_2$, ..., $\omega^3 \times \kappa_2$, ..., and that these cardinals are followed by a well-ordered cardinal $\kappa_3$. And because the sequence of the cardinals $\kappa_1 <_c \kappa_2 <_c \kappa_3 <_c$ ... satisfies the conditions of Theorem 11.10, these cardinals are followed by some well-ordered cardinal, say $\kappa'$..., etc without an end. Let $\kappa$ be as large well-ordered cardinal as we please. There will always exist larger well-ordered cardinals $2 \times \kappa$, $3 \times \kappa$, .... We must not forget that all these cardinals are countable sets in ZF.

**COFINALITY OF CARDINALS** The following definition of cofinality is based on the definition of cofinality in [G]. The definition of cofinality in [H&J] (p. 163, Definition 2.7) gives in AS trivial and therefore uninteresting results.

**11.11 Definition** The *cofinality* of a well-ordered cardinal $\kappa$ is the shortest initial segment $\lambda$ of $\kappa$ for which there exists such indexed system of sets $A_j <_c \kappa$ for $j \in \lambda$ that $\kappa =_c \sum_{j\in\lambda} A_j$. If $\lambda$ is the cofinality of $\kappa$ then we write $\mathrm{cf}(\kappa) = \lambda$.

By this definition the cofinality of the cardinals of cardinality 0, 1 or $\mathcal{P}(\omega)$ does not exist.

**11.12 Theorem** *If* $\mathrm{cf}(\kappa)$ *exists then* $\mathrm{cf}(\kappa) \leq_o \omega$.

*Proof.* By Theorems 7.7 and 7.8 the nodes of the developing tree of a well-ordered cardinal $\kappa$ can be arranged in such a way that their order determines a well-ordering of $\kappa$. Let the nodes on the utmost right infinite path of $\mathbf{T}_\kappa$, which is arranged in such a way, be $v_0$ (the root), $v_1$, $v_2$, $v_3$, ... Let from the node $v_0$ develop a set $A_0$, from the node $v_1$ a set $A_1$, ... etc. The sets $A_0 \setminus A_1$, $A_1 \setminus A_2$, $A_2 \setminus A_3$, ... are all smaller than $\kappa$ with regard to $<_c$ since $A_i \setminus A_{i+1} \leq_c A_0 \setminus A_{i+1}$ and $A_0 \setminus A_{i+1}$ is a proper initial segment of the cardinal $\kappa$ for every $i$. Since the sets $A_i \setminus A_{i+1}$ are also mutually disjoint and their union is equal to $\kappa$, we have $\kappa =_c \sum_{i\in\omega} A_i$. Hence $\mathrm{cf}(\kappa) \leq_o \omega$. □

If the cofinality of a well-ordered cardinal exists then it is isomorphic either to 2 or to $\omega$.

**11.13 Theorem** *Let* $\kappa >_c 1$ *be a well-ordered cardinal. If* $d_{CB}(\kappa) >_c 1$ *then* $\mathrm{cf}(\kappa) =_o 2$, *and if* $d_{CB}(\kappa) =_c 1$ *then* $\mathrm{cf}(\kappa) =_o \omega$.

*Proof.* We arrange the nodes of $\mathbf{T}_\kappa$ in such a way that their order determines a well-ordering $<$ of $\kappa$. Let $\alpha =_o r_{CB}(\kappa)$ and let $n =_o d_{CB}(\kappa)$. Then $\kappa^{(\alpha)} = \{m_1, m_2, ..., m_n\}$. Let $\kappa_1 = \{x \mid x \in \kappa \wedge x \leq m_1\}$ and let $\kappa_i = \{x \mid x \in \kappa \wedge m_{i-1} < x \leq m_i\}$ for $2 \leq i \leq n$. These sets exist since it is not difficult to see that they contain all their accumulation points in $\mathcal{T}_a$. They are mutually disjoint and $\kappa = \kappa_1 \cup \kappa_2 \cup \ldots \cup \kappa_n$. Every $\kappa_i$ is a cardinal since $m_i \in \kappa_i^{(\alpha)}$ while for every proper initial segment $S$ of $\kappa_i$ we have $S^{(\alpha)} = \emptyset$. Hence $S <_c \kappa_i$ by Theorem 8.21. We have $r_{CB}(\kappa_i) =_o \alpha =_o r_{CB}(\kappa)$ and $d_{CB}(\kappa_i) =_o 1$ for $1 \leq_o i \leq_o n$. Thus $\kappa$ can be expressed as a sum of $n$ cardinals which are all smaller than $\kappa$ if $d_{CB}(\kappa) >_o 1$. In this case the set $\kappa' = \kappa \setminus \kappa_n$ is equal to the nonempty set $\kappa_1 \cup \kappa_2 \cup \ldots \cup \kappa_{n-1}$ since no element of $\kappa'$ is an element of $\kappa_n$. And $\kappa'$ is a cardinal. Namely, $r_{CB}(\kappa') =_o \alpha$ and $d_{CB}(\kappa') = n - 1$ while for every proper initial segment $S$ of $\kappa'$ we have either $r_{CB}(S) <_o \alpha$ or $r_{CB}(S) =_o \alpha$ and $d_{CB}(S) <_o n - 1$ since $m_{n-1}$ being the maximal element of $\kappa'$ is not an element of $S$. Thus $S <_c \kappa'$ by Theorem 8.21. Therefore $\kappa$ is of the same cardinality as the sum of two smaller cardinals $\kappa'$ and $\kappa_n$ and hence $\mathrm{cf}(\kappa) =_o 2$.

If, however, $d_{CB}(\kappa) =_o 1$, then we cannot have $\kappa =_c \kappa_1 + \kappa_2 + ... + \kappa_n$ for some natural number $n > 1$ where $\kappa_1$, $\kappa_2$, ..., $\kappa_n$ are cardinals being all smaller than $\kappa$. Let's say that this would hold and let $\kappa_{\max}$ be the greatest cardinal, with regard to $<_c$, among the cardinals $\kappa_i$ for $i \leq n$. As $\kappa_{\max} <_c \kappa$ and $d_{CB}(\kappa) =_o 1$, we should have $r_{CB}(\kappa_{\max}) <_o r_{CB}(\kappa)$ by Theorem 8.21. This, Theorem 8.20 and Lemma 10.12 would then imply

$$r_{CB}(\kappa) =_o r_{CB}(\kappa_1 + \kappa_2 + ... + \kappa_n) \leq_o r_{CB}(n \times \kappa_{\max}) =_o r_{CB}(\kappa_{\max}) <_o r_{CB}(\kappa)$$

which is a contradiction. Hence $\text{cf}(\kappa)$ cannot be isomorphic to any finite N-ordinal and so $\text{cf}(\kappa) =_o \omega$ by Theorem 11.12. □

**11.14 Corollary** *An infinite well-ordered cardinal $\kappa$ is a limit cardinal precisely when* $\text{cf}(\kappa) =_o \omega$.

*Proof.* Theorems 11.13 and 11.9. □

## REGULAR, SINGULAR AND INACCESSIBLE CARDINALS

**11.15 Definition** A cardinal $\kappa$ is *regular* if it is of greater cardinality than any sum of less than $\kappa$ cardinals of smaller cardinality than $\kappa$. A cardinal which is not regular is *singular*.

**11.16 Theorem** *A cardinal is regular precisely when it is of cardinality* 0, 1, 2, $\omega$ *or* $\mathcal{P}(\omega)$.

*Proof.* The cardinals 0, 1 and 2 are obviously regular. Every other finite cardinal is singular since for such a cardinal holds $\kappa =_c \kappa \setminus \min(\kappa) + 1$ where $\kappa \setminus \min(\kappa)$ and 1 are cardinals being both of smaller cardinality than $\kappa$. The cardinal $\omega$ is regular because $\omega$, as a representative of the smallest infinite cardinals with regard to $<_c$, is of greater cardinality than any finite sum of finite cardinals. By Theorem 11.12 the cofinality of any well-ordered cardinal $\kappa$ is isomorphic at most to $\omega$. In other words, any well-ordered cardinal is of the same cardinality as a sum of at most $\omega$ smaller cardinals. Thus every well-ordered cardinal, which is greater than $\omega$, is singular. Any cardinal which is of smaller cardinality than $\mathcal{P}(\omega)$ is well-orderable. But cardinals of cardinality $\mathcal{P}(\omega)$ are not well-orderable and are therefore, by Theorem 9.6, of greater cardinality than any well-ordered sum of well-ordered cardinals. Therefore any cardinal of cardinality $\mathcal{P}(\omega)$ is a regular cardinal. □

**11.17 Definition** An infinite cardinal $\kappa$ is a *strong limit* cardinal if for every cardinal $\kappa' <_c \kappa$ holds $\mathcal{P}(\kappa') <_c \kappa$. A cardinal is *inaccessible* if it is regular, limit and not countable. A cardinal is *strongly inaccessible* if it is strong limit, regular and not countable.

**11.18 Theorem**

*a) A cardinal $\kappa$ is a strong limit cardinal if and only if $\kappa =_c \omega$.*
*b) A cardinal $\kappa$ is inaccessible if and only if $\kappa =_c \mathcal{P}(\omega)$.*
*c) Strongly inaccessible cardinals do not exist.*

*Proof.*

a) The cardinals of cardinality $\omega$ are the smallest strong limit cardinals. If, however, a cardinal $\kappa$ is of greater cardinality than $\omega$, and if we put $\kappa' = \omega$, then $\mathcal{P}(\kappa') <_c \kappa$ does not hold. For this reason no cardinal of greater cardinality than $\omega$ is a strong limit cardinal.

b) As finite cardinals and $\omega$ are countable cardinals, Theorem 11.16 implies that only the cardinals of cardinality $\mathcal{P}(\omega)$ can be inaccessible. The cardinals of this cardinality are the only not well-orderable cardinals. Hence they are limit cardinals by Definition 11.7 and Theorem 11.4. Therefore they are inaccessible.

c) A simple consequence of a). □

# 12. The Size of Infinite Sets in AS

**12.1 Observation** In this section we shall make some intuitive observations in the area of infinity in AS. After everything we have found out till now it would be difficult to accept cardinality as a proper measure of the sizes of infinite sets in AS. Namely, cardinality of sets is based on injections and injections in AS do not preserve only the size of sets but also their topological structure in $\mathcal{T}_a$. Since $\mathcal{T}_a$ in AS is equivalent to the interval topology in ZF, every set in AS can be represented as a set of real numbers in ZF. And since every set in AS is closed in $\mathcal{T}_a$, the corresponding set of real numbers is closed in the interval topology in ZF. If an infinite set of real numbers is closed in the interval topology in ZF then it is either countable or it contains a perfect subset of cardinality $\mathcal{P}(\omega)$. This and Theorem 7.8 imply that every well-orderable set in AS is countable in ZF. In the previous section we have found out that there is a rich variety of well-orderable infinite cardinals in AS. However all these cardinals are countable in ZF and are therefore of equal size in ZF. For this reason it is unlikely that the view would be accepted by which well-orderable infinite sets in AS would be considered as being of different infinite size. Thus the only remaining question is whether such sets are of smaller size then those of cardinality $\mathcal{P}(\omega)$. It turns out that also this distinction of infinite sizes is highly questionable. Namely, the only argument in the favor of this distinction is the following one. In ZF we can prove, by diagonal procedure, that it is not possible to construct a bijection between $\boldsymbol{N}$ and $\mathcal{P}(\boldsymbol{N})$, since in every such attempt we can find an element of $\mathcal{P}(\boldsymbol{N})$ which cannot be paired with any $n \in \boldsymbol{N}$. However, the argument that an infinite set is of larger size than another infinite set because at every attempt of the construction of a bijection between these two sets there remains one element of the first set which is not paired with any element of the second set, is not a very convincing one. Even more. The following argument intuitively shows that $\boldsymbol{N}$ and $\mathcal{P}(\boldsymbol{N})$ have the same number of elements. Let's take a complete binary tree of depth $n$. Let $N_n$ be the set (in ZF) of its internal nodes and let $P_n$ be the set of its paths that start in the root and continue to the bottom of the tree. Then the number of the elements of $N_n$ is $2^n - 1$ and the number of the elements of $P_n$ is $2^n$. Thus the set $N_n$ has only one element less than the set $P_n$ for every finite $n$. Consequently, the most natural conclusion would be that the infinite set $N_\omega$ has only one element less than the infinite set $P_\omega$. However, it is possible to put $N_\omega$ in bijection with $\boldsymbol{N}$ and $P_\omega$ in bijection with $\mathcal{P}(\boldsymbol{N})$. Namely on each of the countably many levels of a complete infinite binary tree there are only finitely many nodes, and every element of $\mathcal{P}(\boldsymbol{N})$ can be uniquely interpreted as an infinite path which starts in the root of this tree. From this point of view the infinite set $\boldsymbol{N}$ has only one element less than the infinite set $\mathcal{P}(\boldsymbol{N})$. In other words, these two sets have the same number of elements.

In AS the view that $\boldsymbol{N}$ and $\mathcal{P}(\boldsymbol{N})$ are of equal size is even more natural. Here injections are homeomorphic embeddings in $\mathcal{T}_a$. So saying that for every injection which maps $\boldsymbol{N}$ into $\mathcal{P}(\boldsymbol{N})$ there always exists an unpaired point in $\mathcal{P}(\boldsymbol{N})$ is the same as saying that for every homeomorphic embedding of a line into a circle there always exists an unpaired point in the circle. This certainly does not imply that a circle has more points than a line.

Taking all these considerations into account the only sensible point of view would be that in AS all infinite sets are of equal size. In other words, all these sets have the same number of elements and their size can be therefore represented by the number $\infty$ or by the ordinal $\omega$. As we know $\infty$ and $\omega$ represent the same number in AS since they are the only accumulation points of the isomorphic sets $\boldsymbol{N}$ and $\omega$ respectively. On the other hand it is clear that a proper set size measure must assign the natural number $m$ to a finite set having $m$ elements. As it is known such set size measure is the *counting measure*. It can be defined in AS in a natural way by approximations of sets compared by the relation $\leq_c$. Namely, by Theorem 3.21, the nodes on the $n$th level of the developing tree of a set $A$ represent the elements of the $(n + 1)$th approximation of $A$. As every one of these nodes has at least one successor in $\mathbf{T}_A$, we have $F(n, A) \leq_c F(n + 1, A)$ for any $n$. And further we have $F(n, A) \leq_c F(\infty, A) = A$ for any $n$ by Theorem 5.15. Thus $A$ has at least as many elements as any one of its approximations. Let $A$ and $B$ be two sets. If for every natural number $i$ there exists such natural number $j$ that $F(i, A) \leq_c F(j, B)$ then it is intuitively clear that $A$ cannot have more elements than $B$. We shall express these observations in the following way:

**12.2 Definition** A set $A$ is *not of greater size than* a set $B$ if for every natural number $i$ there exists such natural number $j$ that $F(i, A) \leq_c F(j, B)$. We shall denote this by $A \leq_q B$ where the subscript $_q$ symbolizes quantity of the elements of a set. If $A \leq_q B$ and $B \leq_q A$, then the sets $A$ and $B$ are of *equal size*, which we shall denote by $A =_q B$. If $A \leq_q B$ but not $B \leq_q A$, then $A$ is *of smaller size* than $B$ which we shall denote by $A <_q B$.

Now we shall prove some basic properties of $\leq_q$ and $=_q$ in order to make these two relations legitimate.

**12.3 Theorem** *The relation $\leq_q$ is an ordering .*

*Proof.* If $A \leq_q B$ does not hold true, then for some natural number $i$ there exists no such natural number $j$ that $F(i, A) \leq_c F(j, B)$. Let's denote this $i$ by $i_0$. Then, by Theorem 7.15, we have $F(i_0, A) >_c F(j, B)$ for every $j$. Therefore for every $j$ there exists such $i$ that $F(j, B) \leq_c F(i, A)$ as one such $i$ is for example $i_0$. This implies $B \leq_q A$. Hence arbitrary two sets are always comparable by size and so $\leq_q$ is a total relation on $U$. This relation is already by definition antisymmetric with regard to $=_q$ which is an equivalence relation as we shall show in the sequel (Theorem 12.4). And $\leq_q$ is also transitive. Namely, if $A \leq_q B$ and $B \leq_q C$ then for every $i$ there exists such $j$ that $F(i, A) \leq_c F(j, B)$ and for this $j$ there exists such $k$ that $F(j, B) \leq_c F(k, C)$. Hence for every $i$ there exists such $k$ that $F(i, A) \leq_c F(k, C)$ which implies $A \leq_q C$. □

**12.4 Theorem** *The relation $A =_q B$ is an equivalence relation.*

*Proof.* The relation is reflexive. Namely, $F(i, A) \leq_c F(i, A)$ for every $i$. Hence $A \leq_q A$ and consequently $A =_q A$. The relation is symmetric already by its definition. The relation is transitive. If $A =_q B$ and $B =_q C$, then $A \leq_q B$ and $B \leq_q C$. Hence $A \leq_q C$ by Theorem 12.3. The symmetry of $=_q$ implies that we have also $C =_q B$ and $B =_q A$. From this we obtain $C \leq_q A$ by the same argument as we have obtained $A \leq_q C$. Consequently, $A =_q C$. □

Let's now show that $=_q$ really is the counting measure.

**12.5 Theorem** *For every set $A$ there exists a unique $n \in \omega$ such that $A =_q n$. If $n <_o \omega$ then $A$ is finite, and if $n = \omega$ then $A$ is infinite.*

*Proof.* If there exists such finite N-ordinal $i$ that $F(i, A) =_c F(i + j, A)$ for every finite N-ordinal $j$, then $A$ is finite. Namely, let $i_0$ be the smallest finite N-ordinal such that $F(i_0 + 1, A) =_c F(i_0 + 1 + j, A)$ for every $j$. Let $v_0, v_1, v_2, ..., v_n$ be all the nodes on the $i_0$th level of the tree $\mathbf{T}_A$. By Theorem 3.21 every of these nodes represents one element of the set $F(i_0 +1, A)$. Since $F(i_0 +1, A) =_c F(i_0 + 2, A)$, every of these nodes has precisely one successor on the level $i_0 + 1$. And since $F(i_0 +1, A) =_c F(i_0 + 3, A)$, every of these successors has precisely one successor on the level $i_0 + 2$ ..., etc. Therefore through every of the nodes on the $i_0$th level passes precisely one infinite path of $\mathbf{T}_A$. Hence $\mathbf{T}_A$ has $n$ infinite paths and so $A$ has $n$ elements. Thus $A =_c n$ for a precisely determined N-ordinal $n$.

If, however, there is no such $i$ that $F(i, A) =_c F(i + j, A)$ for every $j$, then for every $i$ there exists such $j$ that $F(i, A) <_c F(i + j, A)$. This means that the cardinality of $F(i, A)$ can be greater than the cardinality of an arbitrary finite N-ordinal. The isomorphism between the sets $\boldsymbol{N}$ and $\omega$ constructed in the subsection of von Neumann ordinals, Lemma 8.2, and Definition 4.21 imply that $F(j, \omega) = j$ for every finite N-ordinal $j$. Hence for every $j$ there exists such $i$ that $F(j, \omega) = j \leq_c F(i, A)$. That is, $\omega \leq_q A$ and $A$ must be infinite. On the other hand there certainly exists such $k$ that $F(i, A) \leq_c k = F(k, \omega)$ holds for every $i$. Therefore $A \leq_q \omega$ and thus $A =_q \omega$. □

From the philosophical point of view the relation $=_q$ abstracts more properties of sets than the relation $=_c$, and this relation abstracts more properties of sets than the relation $=_o$. This follows from the fact that $A =_o B$ implies $A =_c B$ and that $A =_c B$ implies $A =_q B$.

If we measure the sizes of sets with the relations $=_q$ and $<_q$ instead of $=_c$ and $<_c$ then also the definition of infinity becomes more natural and simpler.

**12.6 Definition** A set is *infinite* if there exists no set of greater size.